\documentclass[a4paper,11pt,reqno]{amsart}
\usepackage{latexsym,amscd,amssymb,url}
\usepackage{extarrows}
\usepackage{verbatim}
\usepackage{amssymb}
\usepackage{soul}
\usepackage{dsfont}

\newcommand*{\longhookrightarrow}{\ensuremath{\lhook\joinrel\relbar\joinrel\rightarrow}}

\usepackage[all,cmtip]{xy}  
\usepackage{graphicx}
\usepackage{xcolor}
\usepackage{enumitem} 
\definecolor{dblue}{rgb}{0,0,.6}
\usepackage[colorlinks=true, linkcolor=dblue, citecolor=dblue, filecolor = dblue, menucolor = dblue, urlcolor = dblue]{hyperref}

\usepackage{mathtools}
\usepackage{mathrsfs}

\numberwithin{equation}{section}

\newtheorem{theorem}{Theorem}[section]

\theoremstyle{plain}

\newtheorem{corollary}[theorem]{Corollary}

\newtheorem{definition}[theorem]{Definition}

\newtheorem{lemma}[theorem]{Lemma}
\newtheorem{notation}[theorem]{Notation}

\newtheorem{proposition}[theorem]{Proposition}
\newtheorem{claim}[theorem]{Claim}

\theoremstyle{definition}
\newtheorem{remark}[theorem]{Remark}
\newtheorem{remarks}[theorem]{Remarks}

\newtheorem{claim-sub}[theoremsub]{Claim}

\newcommand{\Z}{\mathbb Z}
\newcommand{\Q}{\mathbb Q}
\newcommand{\A}{\mathbb A}
\newcommand{\G}{\mathbb G}

\newcommand{\C}{\mathbb C}

\newcommand{\R}{\mathbb R}

\newcommand{\CP}{\mathbb P}
\newcommand{\OO}{\mathcal O}

\let\P\relax 
\newcommand{\P}{\mathbb P}
\newcommand{\bb}[1]{{\mathbb{#1}}}

\newcommand{\ol}[1]{{\overline{#1}}}
\newcommand{\mr}[1]{{\mathscr{#1}}}
\newcommand{\mf}[1]{{\mathfrak{#1}}}
\newcommand{\tn}[1]{{\textnormal{#1}}}
\let\rm\relax 
\newcommand{\rm}[1]{{\mathrm{#1}}}

\newcommand{\im}{\operatorname{im}}
\newcommand{\Ima}{\operatorname{Im}}

\newcommand{\Hom}{\operatorname{Hom}}
\newcommand{\Pic}{\operatorname{Pic}}

\newcommand{\End}{\operatorname{End}}

\newcommand{\Aut}{\operatorname{Aut}}

\newcommand{\id}{\operatorname{id}}

\newcommand{\GL}{\operatorname{GL}}

\newcommand{\pr}{\operatorname{pr}}
\newcommand{\NS}{\operatorname{NS}}
\newcommand{\Sing}{\operatorname{Sing}}

\newcommand{\Ext}{\operatorname{Ext}}

\newcommand{\Ker}{\operatorname{Ker}}

\newcommand{\ca}[1]{{\mathcal{#1}}}
\newcommand{\wt}[1]{{\widetilde{#1}}}

\newcommand{\dashedlongrightarrow}{\xymatrix@1@=15pt{\ar@{-->}[r]&}}
\renewcommand{\longrightarrow}{\xymatrix@1@=15pt{\ar[r]&}}
\renewcommand{\mapsto}{\xymatrix@1@=15pt{\ar@{|->}[r]&}}
\renewcommand{\twoheadrightarrow}{\xymatrix@1@=15pt{\ar@{->>}[r]&}}
\newcommand{\hooklongrightarrow}{\xymatrix@1@=15pt{\ar@{^(->}[r]&}}
\newcommand{\congpf}{\xymatrix@1@=15pt{\ar[r]^-\sim&}}
\renewcommand{\cong}{\simeq}

\newcommand{\p}[1]{\left(#1\right)}

\newcommand{\set}[1]{\left\{ #1 \right\}}
\newcommand{\va}[1]{\left| #1 \right|}

\usepackage[backend=bibtex,style=alphabetic, sorting=nyt]{biblatex}

\begin{document}    

\title[Abelian varieties with no power isogenous to a Jacobian]{Abelian varieties with no power isogenous to a Jacobian}

\author{Olivier de Gaay Fortman} 

\address{Department of Mathematics, Utrecht University,  Budapestlaan 6, 3584 CD Utrecht, The Netherlands.}
\email{a.o.d.degaayfortman@uu.nl}

\author{Stefan Schreieder} 
\address{Institute of Algebraic Geometry, Leibniz University Hannover, Welfengarten 1, 30167 Hannover, Germany.}
\email{schreieder@math.uni-hannover.de}

\date{\today}



 \begin{abstract} 
 Let $X$ be a curve of genus $\geq 4$ that is very general or very general hyperelliptic.
 We classify all the ways in which a power $(JX)^k$ of the Jacobian of $X$ can be isogenous to a product of Jacobians of curves.
 As an application, we show that if $A$ is a very general principally polarized abelian variety of dimension $\geq 4$ or the intermediate Jacobian of a very general cubic threefold, then no power $A^k$ is isogenous to a product of Jacobians of curves.
This confirms various cases of the Coleman--Oort conjecture.
We further deduce from our results some progress on the question whether the integral Hodge conjecture fails for $A$ as above.
 \end{abstract}

\maketitle 

\tableofcontents

 \newpage
 
\section{Introduction} \label{section:introduction}

In this paper we work over the field of complex numbers. 
For a positive integer $g$, let $\ca M_g$ be the moduli space of smooth projective connected curves of genus $g$. If $C$ is a smooth projective curve, let $JC$ denote its Jacobian. 
Our main result is the following theorem. 

\begin{theorem}
\label{theorem:maintheorem1:hyperelliptic} 
Let $Z \subset \ca M_g$ be an irreducible subvariety containing the hyperelliptic locus and let $X$ be a curve that defines a very general point in $Z$. 
Assume that there is an isogeny
$$
\varphi\colon (JX)^k\longrightarrow JC_1 \times \cdots \times JC_n
$$
for some positive integers $k$ and $n$ and some smooth projective connected curves $C_1,\dots ,C_n$ of positive genus. If $g \geq 4$, then $k = n$ and there is an isomorphism $C_i \cong X$ for each $i \in \set{1, \dotsc, n}$. 
\end{theorem}

Theorem \ref{theorem:maintheorem1:hyperelliptic} generalizes in various directions the main result of Bardelli and Pirola in \cite{bardellipirola-curvesofgenusg}, which says that  
the Jacobian $JX$ of a very general curve of genus $g\geq 4$ is not isogenous to the Jacobian of any other curve $C$ with $C\not \cong X$. 

The case $k=1$ of Theorem \ref{theorem:maintheorem1:hyperelliptic} is due to Naranjo and Pirola, see \cite[Theorem 1.1]{naranjopirola2018}. 
By proving Theorem \ref{theorem:maintheorem1:hyperelliptic},  
we fix a gap in their proof (cf.\ Remark \ref{rem:naranjo-pirola:2}) and generalize their theorem to arbitrary powers $(JX)^k$ with $k\geq 1$.  
For $k\geq 2$,  
additional difficulties appear and the proof requires new ingredients, most notably work of Kneser \cite{kneser-quadratisch} on the classification of integral inner product spaces of dimension at most 16,  a generalization of a theorem of Lu and Zuo \cite{luzuo-shimuracurves} on Shimura curves in the Torelli locus, and a recent result of Lazarsfeld and Martin \cite{lazarsfeld2023measures} that grew out from their study of various measures of irrationality. 

We will use 
Theorem \ref{theorem:maintheorem1:hyperelliptic} to prove the following.
 
\begin{theorem} \label{thm:cor:maintheorem2}
Let $A$ be either the intermediate Jacobian of a very general cubic threefold, or a very general principally polarized abelian variety of dimension $g \geq 4$. Then there exists no integer $k \geq 1$ such that $A^k$ is isogenous to a product of Jacobians of smooth projective curves. 
\end{theorem}

Theorem \ref{thm:cor:maintheorem2} implies 
that no power of a very general principally polarized abelian variety of dimension $g \geq 4$ is isogenous to the Jacobian of a smooth projective curve. 
In particular, this proves an instance of the Coleman--Oort conjecture \cite[Expectation 4.2]{moonenoort-specialsubvarieties}, which predicts that for $g\geq 8$, no positive-dimensional special subvariety $Z \subset \ca A_g$ is generically contained in the Torelli locus.
More precisely, Theorem \ref{thm:cor:maintheorem2} has the following consequence:

\begin{corollary} \label{cor:coleman-oort}
Let $g = hk$ with $h \geq 4$ and $k\geq 2$. Let $Z \subset \ca A_g$ be a subvariety such that the general element of $Z$ is isogenous to the $k$-th power of a general principally polarized abelian variety of dimension $h$. Then $Z \subset \ca A_g$ is a special subvariety that satisfies the Coleman--Oort conjecture. In particular, the generic point of $Z$ does not lie in the Torelli locus.
\end{corollary}

Notice that the union of all subvarieties $Z \subset \ca A_g$ as in the above corollary is stable under Hecke translation, hence dense in $\ca A_g$ for the euclidean topology. 

The assertion of Theorem \ref{thm:cor:maintheorem2} that concerns intermediate Jacobians of cubic threefolds seems to go beyond what is predicted by the Coleman--Oort conjecture. 

Another consequence of Theorem \ref{theorem:maintheorem1:hyperelliptic} is the following result.

\begin{corollary} \label{corollary:IHCconsequence}
Let $A$ be either the intermediate Jacobian of a very general cubic threefold, or a very general principally polarized abelian variety of dimension $g \geq 4$. Let $A_1$ be an abelian variety isogenous to a power of $A$, and let $A_2$ be an abelian variety with $\Hom(A,A_2) = 0$. Then $A_1 \times A_2$ is not isomorphic (as unpolarized abelian varieties) to a product of Jacobians of curves.
\end{corollary}

In the above corollary, we cannot exclude that $A_1\times A_2$ is isogenous to a product of Jacobians. In fact, for any abelian variety $A$, a sufficiently general complete intersection curve $C\subset A$ has the property that $JC$ is isogenous to $A\times A_2$ for some abelian variety $A_2$ with $\Hom(A,A_2)=0$.

It is a famous open problem to decide whether the integral Hodge conjecture for curve classes holds on any principally polarized abelian variety. 
This is partly motivated by \cite{voisin-universalCHgroup}, where it is shown that a smooth cubic threefold is not stably rational, if the minimal class of its intermediate Jacobian is not algebraic.
The question on the stable rationality of cubic threefolds is in turn open since the work of Clemens--Griffiths in the nineteen-seventies \cite{clemensgriffiths-cubicthreefold}.

It is shown in \cite{beckmann-degaayfortman,voisin2022cycle} that a principally polarized abelian variety $A$ satisfies the integral Hodge conjecture for curve classes if and only if there is an abelian variety $B$ such that $A\times B$ is isomorphic to a product of Jacobians of curves. If $\dim(A) \geq 4$ and $A$ is very general or the intermediate Jacobian of a very general cubic threefold,  
then Corollary \ref{corollary:IHCconsequence} excludes several possibilities for the abelian variety $B$, as follows.
Since $A$ is simple, we know that $B$ is an extension $0\to B_1\to B\to B_2\to 0$ of an abelian variety $B_2$ with $\Hom(A,B_2)=0$ (i.e.\ $B_2$ does not have $A$ as isogeny factor) by an abelian variety $B_1$ that is isogenous to a power of $A$. 
Corollary \ref{corollary:IHCconsequence} shows that this extension cannot be split: 

\begin{corollary}\label{cor:IHCconsequence-2}
Let $A$ be either the intermediate Jacobian of a very general cubic threefold, or a very general principally polarized abelian variety of dimension $g \geq 4$. 
Let $B$ be an abelian variety.
Assume that $B=B_1\times B_2$, where $B_1$ is isogenous to a power of $A$ and $B_2$ does not have $A$ as an isogeny factor.
Then $A\times B$ is not  isomorphic (as unpolarized abelian varieties) to a product of Jacobians of curves.  
\end{corollary}

The above corollary implies for instance that $A\times B$ is not isomorphic to a product of Jacobians whenever $A$ is as in the corollary and $B$ is a product of simple abelian varieties.

A natural strategy to prove the integral Hodge conjecture for curve classes on an abelian variety $A$ is to construct isogenies of coprime degrees from $k$-th powers of $A$ to Jacobians of curves. 
This approach for $k=1$ allowed Voisin to prove the property for special families of intermediate Jacobians of cubic threefolds in \cite{voisin-universalCHgroup}. 
Similarly, Beckmann and the first named author used this approach in \cite{beckmann-degaayfortman} to prove that the locus of principally polarized abelian varieties that satisfy the integral Hodge conjecture for curve classes is dense in moduli. 
Theorem \ref{thm:cor:maintheorem2} above shows that such an approach does not work for a very general abelian variety of dimension at least four, nor for the intermediate Jacobian of a very general cubic threefold.

\begin{remark}
A simple dimension count shows that there are complex abelian varieties of dimension $\geq 4$ that are not isogenous to the Jacobian of a curve.
While such an argument cannot work over countable fields, the statement remains true over $\overline \Q$ by work of Chai--Oort \cite{chai-oort} and Tsimerman \cite{tsimerman}, and, via a different method, by Masser--Zannier \cite{masser-zannier}.  
Our paper shows that there are abelian varieties $A$ over $\C$ such that no power of $A$ is isogenous to the Jacobian of a curve, see Theorem \ref{thm:cor:maintheorem2}. 
It is reasonable to ask whether such examples exist over $\overline \Q$ as well; first results in this direction are proven by Chen, Lu and Zuo in \cite{MR4332480}.   
\end{remark} 

As aforementioned,  an ingredient in the proofs of Theorems \ref{theorem:maintheorem1:hyperelliptic} and \ref{thm:cor:maintheorem2} is a generalization of a result of Lu and Zuo on Shimura curves in the Torelli locus. 
To be precise, note that \cite[Theorem A]{luzuo-shimuracurves} implies that for an elliptic curve $E$ with transcendental $j$-invariant, and for $g\geq 12$, the $g$-th power $E^g$ yields a point in $\mathcal A_g$ that is not in the same Hecke orbit as the Jacobian of a smooth projective connected curve. 
In other words, there is no smooth projective connected curve $C$ for which there exists an isogeny $E^g \to JC$ that respects the natural polarizations up to a positive integer multiple
(cf.\ Remark \ref{remark:be-aware}). 
In Appendix \ref{appendix:B} to this paper, we show how to deduce from the results of Lu and Zuo in \cite{luzuo-shimuracurves} the stronger statement, where the compatibility assumption on the polarizations is dropped; see Theorem \ref{theorem:verygeneral-powers-elliptic-isogenies} for the precise statement.

\subsection{Outline of the argument} \label{section:outline-of-the-argument}

\subsubsection{Theorem \ref{theorem:maintheorem1:hyperelliptic} implies Theorem \ref{thm:cor:maintheorem2}.} \label{subsec:proofs-intro}
Let us first explain how to deduce Theorem \ref{thm:cor:maintheorem2} from Theorem \ref{theorem:maintheorem1:hyperelliptic}.
To this end, let $A$ be either the intermediate Jacobian of a very general cubic threefold or a very general principally polarized abelian variety of dimension $\geq 4$.
We assume for a contradiction that there is a smooth projective curve $C$ and an isogeny $f \colon A^k \to JC$ for some $k\geq 1$.
 
Our assumptions on $A$ ensure by \cite{collino-fanofundamentalgroup} that $A$ specializes to $JX$, where $X$ is a very general hyperelliptic curve of genus $g=\dim(A)$. This yields a specialization of $C$ to a compact type curve $C_0$ and an isogeny $f_0 \colon (JX)^k \to JC_0$. 
Applying Theorem \ref{theorem:maintheorem1:hyperelliptic} to $f_0$, we see that there exists an isomorphism $g_0 \colon JC_0 \cong (JX)^k$. The composition $g_0 \circ f_0$ is an isogeny $(JX)^k \to (JX)^k$, given by a matrix $\rm{M}_k(\Z)$ with non-zero determinant as $\End(JX) =\Z$. By an idea from \cite{bardellipirola-curvesofgenusg} (generalized in Lemma \ref{lemma:isogenyVHS}), we deduce that $JC \cong A^k$ as unpolarized abelian varieties. 

By the above step, we are reduced to the case where $f\colon A^k \to JC$ is an isomorphism of complex tori.
Using this isomorphism, the canonical principal polarization on $JC$ induces an indecomposable principal polarization on $A^k$.
We will study all principal polarizations on $A^k$ in Section \ref{subsection:polarizations-on-powers} and see that such polarizations correspond to indecomposable integral inner product spaces, that is, indecomposable pairs  $(\Z^k,\alpha)$ where $\alpha \in \mathrm{M}_k(\Z)$ is symmetric and the bilinear form it defines on $\Z^k$ is positive definite and unimodular. Kneser classified such pairs for $k\leq 16$, see \cite{kneser-quadratisch}: there is one such space  in dimension $k=1,8,12,14,15$ and there are two for $k=16$.
 (As an aside, we point out that the number of such spaces grows exponentially with $k$, e.g.\ there are more than $10^{51}$ such spaces for $k=40$, and there is no classification for large $k$, see \cite[p.\ 28, Remark 1]{milnorhusemoller}.)

The aforementioned classification implies that the isomorphism $f \colon A^k \xrightarrow{\sim}JC$ is actually an isomorphism of principally polarized abelian varieties $f \colon (A^k, \alpha) \xrightarrow{\sim} (JC, \Theta_C)$, where $(A^k,\alpha)$ denotes the principally polarized abelian variety associated to $A^k$ and some indecomposable integral inner product space $(\Z^k,\alpha)$. 
To conclude the argument, we specialize $A$ to a product $E\times B$ where $E$ is a very general elliptic curve and $B$ is a principally polarized abelian variety with $\Hom(E,B)=0$. 
We find that $(A^k, \alpha)$ specializes to the product of principally polarized abelian varieties $(E^k\times B^k,\alpha)\cong (E^k,\alpha)\times (B^k,\alpha)$, and this principally polarized abelian variety is isomorphic to the Jacobian of the specialization of the curve $C$. 
Since $(\Z^k,\alpha)$ is indecomposable, so is the principally polarized abelian variety $(E^k,\alpha)$, which must thus be isomorphic to the Jacobian of an irreducible curve. 
In particular, as $E$ is very general, we get $k\leq 11$ in view of Theorem \ref{theorem:verygeneral-powers-elliptic-isogenies} in Appendix \ref{appendix:B}, which generalizes work of Lu and Zuo \cite{luzuo-shimuracurves}.

At this point Kneser's classification can be applied and we are reduced to the case where either $k=1$, or $k=8$ and $\alpha$ is induced by the $E_8$-lattice.
The latter is ruled out by comparing the automorphism group of the $E_8$-lattice with the automorphism group of the Jacobian of a smooth curve of genus eight. 
We thus arrive at $k=1$. 
This yields an isomorphism $f \colon A\cong JC$ of unpolarized abelian varieties, which has to respect the polarizations on both sides (because $A$ is very general, hence has Picard rank one). 
We have finally arrived at a contradiction because the principally polarized abelian variety $A$ is not isomorphic to a Jacobian by our assumptions (when $A$ is the intermediate Jacobian of a smooth cubic threefold, this follows from \cite{clemensgriffiths-cubicthreefold}).

\subsubsection{Sketch of the proof of Theorem \ref{theorem:maintheorem1:hyperelliptic}.} 
This is the technical heart of the paper.  
The assumptions that $Z\subset \mathcal M_g$ contains the hyperelliptic locus $\ca H_g \subset \ca M_g$ and that $[X]\in Z$ is very general 
quickly reduce the proof of Theorem \ref{theorem:maintheorem1:hyperelliptic} to the case where $n=1$ and $Z = \ca H_g$. 
Thus, we may assume $X$ is a very general hyperelliptic curve and there exists an isogeny $\varphi \colon (JX)^k \to JC$ for some smooth projective connected curve $C$.
We then have to show $k=1$ and $C\cong X$. For this, the idea is to split the proof into two steps, and prove that:
\begin{enumerate} [label=(\roman*)]
\item there exists an isomorphism of unpolarized abelian varieties $(JX)^k \cong JC$; and\label{item:(i)}
\item  if $\varphi \colon (JX)^k \cong JC$ is an isomorphism, then $k = 1$ and $C \cong X$. \label{item:(ii)}
\end{enumerate}

Let us first sketch how to prove item \ref{item:(i)}. Define $H \coloneqq H^1(JX,\Z)$. The isogeny $\varphi \colon (JX)^k \to JC$ induces an embedding $\varphi^\ast \colon H^1(JC,\Z) \hookrightarrow H^{\oplus k}$ whose image we denote by $M \subset H^{\oplus k}$. We remark that to prove item \ref{item:(i)}, it suffices to prove that:
\begin{enumerate} [label=(\roman*')]
\item $M = \alpha \cdot H^{\oplus k}$ for some $\alpha \in \rm{M}_k(\Z)$ with non-zero determinant.\label{item:(i')}
\end{enumerate}
At the core of the proof of item \ref{item:(i')} lies a carefully chosen degeneration of  
 $X$ to four different one-nodal hyperelliptic curves $X_{01}, X_{02}, X_{03}$ and $X_{04}$. These degenerations induce degenerations of $C$ to nodal curves $C_{01}, \dotsc, C_{04}$. In fact, for each $i \in \set{1,2,3,4}$, we deform the curve $X_{01}$ while keeping its normalization $\wt X_{0i}$ fixed, which moves the extension class on $J \wt X_{0i}$ associated to $J X_{0i}$. 
 We compare the latter with the extension class on $J\wt C_{0i}$ associated to $JC_{0i}$. More precisely, in Proposition \ref{prop:moving} we use the induced isogeny $\wt \varphi_i \colon (J\wt X_{0i})^k \to J \wt C_{0i}$ to compare the two extension classes and the way they move to conclude that each irrational connected component $K \subset \wt C_{0i}$ of the smooth curve $\wt C_{0i}$ must be hyperelliptic since $\wt X_{0i}$ is hyperelliptic. 
 Consequently, a result of Lazarsfeld--Martin (see Proposition \ref{proposition:lazarsfeldmartin}) implies that the genus of $K$ equals the genus of the curve $\wt X_{0i}$. 
 We then apply a simpler version of Theorem \ref{theorem:maintheorem1:hyperelliptic} due to Naranjo--Pirola, see Theorem \ref{theorem:C-hyperell-g=3}, in which one has the additional assumptions that $k = n = 1$ and \emph{both} curves are hyperelliptic. 
 This allows us to conclude that $K \cong \wt X_{0i}$ for each irrational connected component $K \subset \wt C_{0i}$, which implies that $(J\wt X_{0i})^k \cong J\wt C_{0i}$.
 In other words, item \ref{item:(i)} holds for the normalizations of $X_{0i}$ and $C_{0i}$ for each $i$, and we aim to deduce from this the statement in \ref{item:(i')}. 

The degenerations come together with specialization maps 
 $
 H^1(X_{0i},\Z) \hookrightarrow H^1(X,\Z) = H.
 $
 The above implies that for each $i \in \set{1,2,3,4}$, there exists a matrix $\alpha_i \in \rm{M}_k(\Z)$ with non-zero determinant, such that if $W_0H^1(X_{0i},\Z) \coloneqq H^1(X_{0i},\Z) \cap W_0H^1(X_{0i},\Q)$, then 
 \begin{align*} 
 M \cap H^1(X_{0i},\Z)^{\oplus k} \equiv \alpha_i \cdot H^1(X_{0i},\Z)^{\oplus k} \quad \mod (W_0H^1(X_{0i},\Z))^{\oplus k} \quad \quad \forall i \in \set{1,2,3,4}. 
 \end{align*}
 The way in which we chose our four degenerations $X \rightsquigarrow X_{0i}$ allows us to apply two technical linear algebra statements to the above congruences, see Lemmas \ref{lemma:matricessaturatedsubmodules} and \ref{lemma:appendix-k=arbitrary-new}. 
 The result is that $M = \alpha_i \cdot H^{\oplus k}$ for each $i$. In particular, this proves the above item \ref{item:(i')} as desired. 
  
 To finish the proof of Theorem \ref{theorem:maintheorem1:hyperelliptic}, it remains to prove item \ref{item:(ii)} above. For this, we establish the following result of independent interest; for a more general version of the statement, see Theorem \ref{cor:iso-of-products-of-curves-main-main} in Section \ref{section:polarizations}.  

\begin{theorem}\label{cor:iso-of-products-of-curves-main}
Let $g\in \Z_{\geq 1}$ and let $Z\subset \mathcal M_g$ be an irreducible subvariety which contains the hyperelliptic locus.
Let $X$ be a curve that defines a very general point in $Z$.
If for some $k,n \geq 1$, there is an isomorphism of unpolarized abelian varieties $JC_1 \times \cdots \times JC_n \cong (JX)^k$ for some smooth projective connected curves $C_i$ of positive genus, then $k = n$ and $C_i \cong X$ for each $i$.  
\end{theorem}

We emphasize that Theorem \ref{cor:iso-of-products-of-curves-main} works in any genus $g\geq 1$, while the assumption $g\geq 4$ in Theorem \ref{theorem:maintheorem1:hyperelliptic} is necessary. For instance, Theorem \ref{cor:iso-of-products-of-curves-main} implies that for an elliptic curve $E$ with transcendental $j$-invariant, no power $E^k$ with $k \geq 2$ is isomorphic as unpolarized abelian varieties to the Jacobian of a smooth projective connected curve. 

The assumption in Theorem \ref{cor:iso-of-products-of-curves-main} cannot be weakened to only ask that $(JX)^k$ is isogenous to a product of Jacobians of curves; for example, the third power of a very general elliptic curve is isogenous to countably many Jacobians of smooth projective connected curves of genus three.  
The subtlety of the result lies in the fact that for $k \geq 2$, the abelian variety $(JX)^k$ carries various principal polarizations and so the Torelli theorem can a priori not be applied directly. 
Besides Theorem \ref{theorem:verygeneral-powers-elliptic-isogenies} and the aforementioned result  from \cite{kneser-quadratisch}, our proof depends on various fortunate numerical coincidences, see Section \ref{sec:specialsubvarieties}.

\subsection{Conventions} \label{section:conventions}
We work over the field of complex numbers.
Varieties are integral separated schemes of finite type over $\C$. 
In particular, varieties  (and subvarieties) are integral and hence irreducible by convention.

A complex point $b\in B$ of a complex variety $B$ is very general if there is a finitely generated subfield $K\subset \C$ and a model $B_0$ of $B$ over $K$, i.e.\ $B=B_0\times_K\C$, such that $b$ maps to the generic point of $B_0$ via the natural projection $B\to B_0$.
If $X\to B$ is a family, then a very general fibre $X_b$ is a fibre over a very general point $b$ as above such that in addition, $X\to B$ can be defined over the field $K$.
(The question whether the image of $b$ is the generic point in $B_0$ depends on the choice of $K$ and we simply ask that this holds for some $K$; one can make this independent of choices by asking that $K$ has minimal transcendence degree such that models as above exist.)
The set of all very general points of a complex variety is the complement of a countable union of proper closed subsets.
If $M$ is an integral fine moduli space parametrizing complex varieties (e.g.\ curves) with some properties, then an object of $M$ is called very general if the corresponding moduli point in $M$ is a very general point of the universal family. 
We may thus think about a very general object of $M$ as (a base change  of) the geometric generic fibre of the universal family.
This allows us to specialize very general objects of $M$ to any other given object, which is the key property that we will use in this paper. 
If $M$ is the coarse moduli space for a moduli problem that admits a finite cover by a scheme $M' \to M$ over which a universal family exists, e.g.\ the cover of $\ca M_g$ obtained by adding level $\geq 3$ structure, then we define a very general point of $M$ as the image of a very general point of $M'$.
For instance, an elliptic curve is very general if and only if its $j$-invariant is transcendental.

A curve is a reduced projective scheme of pure dimension one over $\C$, or the analytification of such a scheme. A curve $X$ of arithmetic genus $g \geq 2$ is semi-stable (resp.\ stable) if its singularities are at most nodal and each rational connected component of its normalization contains at least two (resp.\ three) points lying over nodes of $X$. A curve is of compact type if its dual graph is a tree. 
A family of curves is a proper flat morphism $p \colon \ca X \to B$ of finite type schemes, or complex analytic spaces, such that for each $b \in B$, the fibre $X_b = p^{-1}(b)$ is a curve which is (unless mentioned otherwise) connected. We say that a family of curves $p \colon \ca X \to B$ is a family of semi-stable (resp.\ stable, smooth, compact type, nodal) curves if the curve $X_b$ is semi-stable (resp.\ stable, smooth, of compact type, has at most nodal singularities) for each $b \in B$. 
If $X$ denotes a complex quasi-projective variety, then we denote by $W_kH^i(X,\Q)$ the $k$-th piece of the weight filtration of the associated mixed Hodge structure, see e.g.\ \cite{petersteenbrink}.
If $H^i(X,\Z)$ is torsion-free (e.g.\ if $i=1$), then we also write $W_kH^i(X,\Z) \coloneqq H^i(X,\Z)\cap W_kH^i(X,\Q)$.
If $A$ is an abelian variety with dual abelian variety $A^\vee$, then a polarization on $A$ is the isogeny $\lambda \colon A \to A^\vee$ associated to an ample line bundle $\ca L$ on $A$; thus $\lambda(x) = [t_x^\ast(\ca L) \otimes \ca L^{-1}] \in \Pic^0(A) = A^\vee$ for $x \in A$, where $t_x \colon A \to A$ denotes the translation by $x$ map.

\subsection{Acknowledgements} 
Starting with the influential work of Bardelli--Pirola \cite{bardellipirola-curvesofgenusg}, many authors have studied irreducible curves of low genus on (sufficiently general) $g$-dimensional abelian varieties, see e.g.\ \cite{naranjopirola1994, lange-sernesi, marcucci-genusofcurves, marcucci-naranjo-pirola, naranjopirola2018}. 
We were greatly inspired by these works. 
We would also like to thank Frans Oort and the referee for their valuable comments on this paper. 
We are grateful to Kang Zuo for answering our questions concerning \cite{luzuo-shimuracurves}.  

This project has received funding from the European Research Council (ERC) under the European Union's Horizon 2020 research and innovation programme under grant agreement N\textsuperscript{\underline{o}} 948066  (ERC-StG RationAlgic).
 
\section{Preliminaries} \label{section:setup}

In this section, we gather various preliminary results. 
 
\subsection{Isogenies between powers of abelian schemes and cohomology} We start with:

\begin{lemma} \label{lemma:isogenyVHS}
Let $k \in \Z_{\geq 1}$. Let $S$ be a connected complex manifold, $f \colon \ca A \to S$ and $g \colon \ca B \to S$ families of compact complex tori, and $\psi \colon \ca  A^k \to \ca B$ a family of isogenies over $S$, where $\mathcal A^k$ denotes the $k$-fold fibre product of $\mathcal A$ over $S$.  
Suppose that, for some $t \in S$, the induced map $\psi_t^\ast \colon H^1(B_t,\Z) \to H^1(A_t,\Z)^{\oplus k}$ satisfies $\psi_t^\ast (H^1(B_t,\Z)) = \alpha \cdot H^1(A_t,\Z)^{\oplus k} \subset H^1(A_t,\Z)^{\oplus k}$ for some $\alpha \in \rm{M}_k(\Z)$.  
Then there is an isomorphism $\ca B \cong \ca A^k$ over $S$.  
\end{lemma}

\begin{proof}
The case $k=1$ goes back to \cite{bardellipirola-curvesofgenusg}; the more general version stated above is similar.

First of all, in order to prove $\ca B \cong \ca A^k$ over $S$, it suffices to show that the variation of integral Hodge structure $ R^1g_\ast\Z$ on $S$ is isomorphic to $(R^1f_\ast\Z)^{\oplus k}$.
The matrix $\alpha$ has full rank, because $\psi_t$ is an isogeny and $\psi_t^\ast (H^1(B_t,\Z)) = \alpha \cdot H^1(A_t,\Z)^{\oplus k} $ by assumption.
Hence, fibrewise multiplication by $\alpha$ yields an embedding of integral Hodge structures $\alpha\colon  (R^1f_\ast\Z)^{\oplus k} \to (R^1f_\ast\Z)^{\oplus k}$.
Similarly, fibrewise pushforward along the family of isogenies $\psi$ yields an embedding $ \psi^\ast\colon  R^1g_\ast\Z \to (R^1f_\ast\Z)^{\oplus k}$.
To prove that  $ R^1g_\ast\Z\cong (R^1f_\ast\Z)^{\oplus k}$, it thus suffices to prove that the images of the above embeddings coincide, i.e.\ that the following equality holds:
$$
\im(\psi^\ast\colon  R^1g_\ast\Z \to (R^1f_\ast\Z)^{\oplus k})=\im(\alpha\colon  (R^1f_\ast\Z)^{\oplus k} \to (R^1f_\ast\Z)^{\oplus k} ).
$$
Since $S$ is connected, the above identity can be checked at the single point $t\in S$, where it holds by assumption.
This concludes the proof of the lemma. 
\end{proof}

As a corollary, we obtain the following useful criterion. 

\begin{lemma} \label{lemma:isomorphismcriterion}
    Let $S$ be a smooth connected complex manifold and let $f \colon \ca A \to S$ and $g \colon \ca B \to S$ be families of compact complex tori.
    Let $\psi\colon \ca A^k \to \ca B$ be a family of isogenies of complex tori over $S$. Suppose that, for some $t \in S$, we have $\End(A_{t}) = \Z$ and $B_{t} \cong A_{t}^k$ for some $k \in \Z_{\geq 1}$. Then there is an isomorphism $\ca B \cong \ca A^k$ over $S$. 
\end{lemma}
\begin{proof}
    By Lemma \ref{lemma:isogenyVHS}, it suffices to show that the induced map $\psi_t^\ast \colon H^1(B_{t},\Z) \to H^1(A_{t},\Z)^{\oplus k}$ has image $ \alpha \cdot H^1(A_{t},\Z)^{\oplus k} \subset H^1(A_{t},\Z)^{\oplus k}$ for some $\alpha \in \rm{M}_k(\Z)$. This is clear: as $\End(A_{t}) = \Z$, the composition 
    \[
A_{t}^k \stackrel{\psi_t} \longrightarrow B_{t} \cong A_{t}^k 
    \]
    is given by a matrix $\alpha \in \rm{M}_k(\Z)$, hence the same is true on cohomology. 
\end{proof}

 \subsection{Gau{\ss} maps.}  The goal of this section is to prove Proposition \ref{proposition:hyperelliptic-gauss} below. This proposition says that if $C$ is a smooth connected non-hyperelliptic curve and $m \in \Z_{\geq 1}$, then the surface
$ 
 m(C - C) \subset JC
$ 
contains no hyperelliptic curves. 

For a dominant morphism of varieties $f \colon X \to Y$ with $\dim(X) = \dim(Y)$, let the branch locus $B(f) \subset Y$ be the reduced closed subscheme which is the complement of the largest open subset of $Y$ over which $f$ is \'etale. 
Similarly, define the ramification locus $R(f) \subset X$ of $f$ as the complement of the largest open subset $U\subset X$ such that $f|_U\colon U\to X$ is \'etale. 
These definitions readily extend to the case of a dominant rational map $f \colon X \dashrightarrow Y$ of varieties of the same dimension. 
Namely, if $U \subset X$ is a non-empty open subset on which $f$ restricts to a morphism $f|_U \colon U \to Y$, then we define $R(f)$ as the closure of $R(f|_U)$ in $X$, and $B(f)$ as the closure of $B(f|_U)$ in $Y$. 

Let $C$ be a smooth connected non-hyperelliptic curve (hence of genus $g\geq 3$), and let 
\[
\phi_{K_C} \colon C \longrightarrow \P (H^0(C, K_C) )  = \P\left( T_0JC\right) 
\]
be the canonical embedding. This gives a morphism
\[
s \colon C \times C \longrightarrow  \rm{Gr}(1, \P(T_0JC)), \quad (p,q) \mapsto \left(\text{line spanned by $\phi_{K_C}(p)$ and $\phi_{K_C}(q)$}\right).
\]
\begin{lemma}  \label{lemma:delta-birational}
In the above notation, the morphism $s \colon C \times C \to s(C \times C)$ is generically finite. 
Moreover, if we let $\deg(s) $ denote its degree, then we have 
\[
\deg(s) = 
\begin{cases}
2 & \text{ if } g \geq 4, \\
12 & \text{ if } g = 3. 
\end{cases}
\]
The ramification locus $R(s)$ of $s$ is the diagonal $\Delta_C \subset C \times C$, and the branch locus $B(s)$ of $s$ is the curve $D$ in $\rm{Gr}(1, \P(T_0JC))$ parametrizing lines which are tangent to the canonical curve $\phi_{K_C}(C) \subset \P\left( T_0JC\right)$. Moreover, the induced map $C = \Delta_C = R(s) \to B(s) = D$ is birational. 
\end{lemma}
\begin{proof}
This is well-known (see e.g.\ \cite[Remark 3.1.2]{bardellipirola-curvesofgenusg}). 
\end{proof}

\begin{lemma} \label{lem:diagonal-contraction}
    Let $C$ be a smooth connected non-hyperelliptic curve. Let $\Delta_C \subset C \times C$ be the diagonal, and consider the origin $0 \in C-C \subset JC$. The difference map $a \colon C \times C \to C-C$ restricts to an isomorphism $\left(C\times C\right)\setminus \Delta_C\xrightarrow{\sim} \left(C-C \right) \setminus \set{0}$. 
\end{lemma}
\begin{proof}
Let $p,q, p', q' \in C$. Then $[p - q] = [p' - q'] \in JC$ implies that the divisors $p+q'$ and $p'+q$ on $C$ are linearly equivalent. Hence, they are equal, because $C$ is not hyperelliptic. 
\end{proof}

\begin{lemma} \label{lem:line}
Let $C$ be a smooth connected non-hyperelliptic curve of genus three, with canonical embedding $C \hookrightarrow \P^2$. 
Let $\ell \subset \P^2$ be a very general line and let $p_1, p_2, p_3, p_4$ be four distinct points on $C$ such that $\ell \cap C = \set{p_1, p_2, p_3, p_4} \subset C$. 
Then for any integer $m \geq 1$, the twelve points $m[p_i-p_j] \in JC$ $(i\neq j \in \set{1,2,3,4})$ are pairwise distinct. 
\end{lemma}
\begin{proof}
Assume that there is a very general line $\ell \subset \P^2$  and points $p,q,p',q' \in C\cap \ell$ with $p\neq q$ and $p'\neq q'$, such that $m[p-q] = m[p'-q'] \in JC$ for some  $m \geq 1$.
We need to show $p=p'$ and $q=q'$. 
For a contradiction, we assume that this is not the case.

The line $\ell$ is spanned by the two points $p$ and $q$.
Since $\ell$ is very general, so is $(p,q)\in C\times C$.
In a first step, we note that $p=q'$ and $q=p'$ is impossible, as it leads to $2m[p-q]=0\in JC$ for very general $(p,q)\in C\times C$, which is absurd.
Since $(p,q)\neq (p',q')$ by assumption, we get that $\{p,q\}\neq \{p',q'\}$.

We specialize the line $\ell$ to a line $\ell_0$ so that $p$ and $q$ collapse to a single very general point $p_0\in C$ and $\ell_0$ is the tangent line of $C$ at $p_0$. 
We denote the limit points of $p'$ and $q'$ by $p'_0$ and $q'_0$, respectively. 
Since $\{p,q\}\neq \{p',q'\}$, the equality $p'_0= q'_0$ would imply that $\ell_0$ is a bitangent, which is impossible as their number is finite (there are exactly 28 such lines).
Hence, $p'_0\neq q'_0$. 
Moreover, the identity $m[p-q] = m[p'-q'] \in JC$ specializes to the identity $0 = m [p_0' - q_0'] \in JC$. 
We can further specialize the point $p_0\in C$ at which $\ell_0$ is tangent and find that there is a one-dimensional family of points $(p_0',q_0')\in C\times C$ such that $0 = m [p_0' - q_0'] \in JC$ and $p_0' \neq q_0'$. 
Since the $m$-torsion points of $JC$ are discrete, taking the closure of the above one-dimensional family yields 
a curve in $C\times C$ which is different from the diagonal and which is contracted by the difference map $C\times C\to C-C$.
This contradicts Lemma \ref{lem:diagonal-contraction}, and hence concludes the proof. 
\end{proof}

Let $A$ be an abelian variety of dimension $g \geq 2$ and let $V \subset A$ be a closed subvariety of dimension $k$, with $1 \leq k \leq g-1$. Let $T_0A$ be the tangent space of $A$ at the origin, and consider the canonical trivialization $TA \cong T_0A \times A$ of the tangent bundle $TA$ of $A$. Define $\rm{Gr}(k,T_0A)$ as the Grassmannian of $k$-planes in $T_0A$. 
Recall that, in this setting, the \emph{Gau{\ss} map} 
\[
\xymatrix{
\mathscr G_{V,A} \colon V  \ar@{-->}[r] & \rm{Gr}(k, T_0A) = \rm{Gr}(k-1, \mathbb P(T_0A))
}
\]
is the rational map defined as follows. For a point $x$ in the smooth locus of $V$, the induced map on tangent spaces $T_xV \to T_xA$ is an embedding, whose image can be identified with a $k$-plane in $T_0A$ via $t$; 
we let $\mr G_{V,A}(x) \in \rm{Gr}(k,T_0A)$ be the induced point of the Grassmannian. 

Let $C$ be a smooth connected non-hyperelliptic curve. Define a morphism
\begin{align} \label{definition:a}
    a \colon C \times C \longrightarrow JC, \quad a(x,y) = [x-y] \in \Pic^0(C) = JC.
\end{align}
For $m \in \Z_{\geq 1}$, let $m(C - C) \subset JC$ be the image of the morphism $m \cdot a \colon C \times C \to JC$, and put $C - C = 1(C - C) \subset JC$. 
Consider the Gau{\ss} map $G_m = \mr G_{m(C-C), JC}$, which is a rational map
\begin{align} \label{equation:gauss}
\xymatrix{
G_m \colon m(C - C) \ar@{-->}[r] &    \rm{Gr}(1, \P(T_0JC)). 
}
\end{align}
By \cite[Lemma 3.1.1 \& Remark 3.1.2]{bardellipirola-curvesofgenusg}, the following diagram commutes:
\begin{align} \label{diagram:bardellipirola}
\begin{split}
\xymatrix{
C \times C  \ar@{->>}[r]^-a \ar[d]^s & C - C\ar@{-->}[d]^-{G_1}\ar@{->>}[r]^-m & m (C - C)\ar@{-->}[d]^-{G_m} \\
\rm{Gr}(1, \P(T_0JC)) \ar@{=}[r] & \rm{Gr}(2, T_0JC) \ar@{=}[r] & \rm{Gr}(2, T_0JC). 
}
\end{split}
\end{align}
In particular, we obtain the following rational maps:
\begin{align} \label{equation:rationalmaps}
G_1 \colon C - C \dashrightarrow s(C \times C), \quad G_m \colon m(C-C) \dashrightarrow s(C \times C).
\end{align}

\begin{proposition} \label{proposition:pirolabardelli-branch}
Let $C$ be a smooth connected non-hyperelliptic curve and let $m \geq 1$ be an integer. Consider the commutative diagram \eqref{diagram:bardellipirola} above, and the resulting morphisms \eqref{equation:rationalmaps}. 
The following assertions are true.
\begin{enumerate}
\item \label{item:BP-one} The morphisms $a \colon C \times C \to C - C$ and $m \colon C - C \to m(C-C)$ are birational.  
\item \label{item:BP-two} 
The rational maps $G_1 \colon C - C \dashrightarrow s(C \times C)$ and $G_m \colon m(C - C) \dashrightarrow s(C \times C)$ are generically finite.  
\end{enumerate}
\end{proposition}
\begin{proof}
By Lemma \ref{lemma:delta-birational}, the morphism $s\colon C\times C\to s(C\times C)$ is generically finite.
Hence, the commutativity of diagram \eqref{diagram:bardellipirola} shows that item \eqref{item:BP-two} follows from item \eqref{item:BP-one}.

Let us prove item \eqref{item:BP-one}. 
By Lemma \ref{lem:diagonal-contraction}, the map $a$ is birational, hence it remains to prove that the map $m \colon C-C \to m(C-C)$ is birational. Since diagram \eqref{diagram:bardellipirola} commutes and $a$ is birational, we get that $\deg(s) = \deg(G_1) = \deg(m) \cdot \deg(G_m)$. 
To prove $\deg (m)=1$ it thus suffices to show $\deg(s) \leq \deg(G_m)$. 
Assume first that $g \geq 4$. Let $p \in C$ be any point, and let $q \in C$ be a point such that $[p-q] \in JC$ is not $2m$-torsion. Then, on the one hand, $m[p-q] \neq m[q-p] \in JC$, and on the other hand, $G_m(m[p-q]) = G_m(m[q-p]) \in \rm{Gr}(2, T_0JC)$. Therefore, $2 \leq \deg(G_m)$, and since $\deg(s) = 2$ by Lemma \ref{lemma:delta-birational}, we get $\deg(s) \leq \deg(G_m)$, proving what we want.

Next, assume that $g = 3$, and consider the canonical embedding $C \hookrightarrow \P^2$. Let $\ell \subset \P^2$ be a very general line, so that $\ell \cap C = \set{p_1, p_2, p_3, p_4} \subset C$ for distinct points $p_1, p_2,p_3, p_4$ on $C$. 
The elements $m[p_i-p_j] \in m(C-C)$ for $i \neq j$ are all sent to the same element in $\rm{Gr}(2, T_0JC)$ under the rational map $G_m$. Moreover, the twelve points $m[p_i-p_j] \in JC$ $(i\neq j \in \set{1,2,3,4})$ are pairwise distinct by Lemma \ref{lem:line}. We conclude that $12 \leq \deg(G_m)$, and since $\deg(s) = 12$ by Lemma \ref{lemma:delta-birational}, we get $\deg(s) \leq \deg(G_m)$ and we are done. 
\end{proof}

\begin{remark} \label{remark:NP-error-twelve}
In \cite{naranjopirola2018}, item \eqref{item:BP-one} of Proposition \ref{proposition:pirolabardelli-branch} is proven for $ g \geq 4$. 
We gave some details of the arguments above because we will need the case $g=3$, in which case the claim in loc.\ cit.\ that $\deg(s) = 2$ is incorrect, see Lemma \ref{lemma:delta-birational}. 
\end{remark}

We are now in a position to prove the following proposition. 

\begin{proposition} \label{proposition:hyperelliptic-gauss} 
Let $C$ be a smooth connected curve of genus $g \geq 3$. Suppose that, for some $m \in \Z_{\geq 1}$, there is a non-constant morphism 
\[
f \colon X \longrightarrow m(C - C)\subset JC,
\]
where $X$ is a smooth connected hyperelliptic curve. Then $C$ is hyperelliptic. 
\end{proposition}
\begin{proof} 
Assume that $C$ is non-hyperelliptic; our goal is to arrive at a contradiction.  
Recall that $m \colon C - C \to m(C-C)$ is birational, see item \eqref{item:BP-one} of Proposition \ref{proposition:pirolabardelli-branch}.  
Assume first that $f(X)$ is contained in the branch locus $B(m)$ of $m$. The commutativity of diagram \eqref{diagram:bardellipirola} yields a rational map $X \dashrightarrow B(s)$ defined as the composition 
\[
\xymatrix{
X \ar@{-->}[r]^-f & B(m) 
\ar@{-->}[r]^-{G_m}  & B(G_m \circ m) \ar@{=}[r] & B(G_1) \subset B(G_1 \circ a) \ar@{=}[r]&  B(s).
}
\]
This rational map is non-constant as $B(m)$ must be a curve (as $f$ is non-constant) and $G_m$ is generically finite, see item \eqref{item:BP-two} in Proposition \ref{proposition:pirolabardelli-branch}.  
As $B(s)$ is birational to $C$ by Lemma \ref{lemma:delta-birational}, one obtains a non-constant morphism $X \to C$, proving that $C$ is hyperelliptic (see e.g.~\cite[Lemma 1.1]{schoen-hyperelliptic}), which yields the desired contradiction.

Therefore, the curve $f(X) \subset m(C-C)$ is not contained in $B(m)$, and we obtain a non-constant rational map $X \dashrightarrow C - C$ defined as the composition
\[
\xymatrix{
X 
\ar@{-->}[r]^-f
&
m(C-C) 
\ar@{-->}[r]^-{m^{-1}} &
C - C. 
}
\]
Note that $R(a) = \Delta_C \subset C \times C$ and $B(a) = \set{0} \subset C - C$. 
Consequently, composing the non-constant rational map $m^{-1} \circ f$ with the rational map $a^{-1} \colon C - C \dashrightarrow C \times C$, one  
obtains a non-constant rational map 
$X \dashrightarrow C$ defined as the composition 
\[
a^{-1} \circ m^{-1} \circ f \colon X \dashrightarrow m(C - C) \dashrightarrow C - C \dashrightarrow C \times C.
\]
Thus, $X$ admits a non-constant morphism $X \to C$. This is a contradiction, and we are done. \end{proof}

\begin{remarks} Let $C$ be a non-hyperelliptic curve of genus $g \geq 3$. Essential in the proof of Proposition \ref{proposition:hyperelliptic-gauss} above is to exploit the birational map $(m \circ a)^{-1} \colon m(C-C) \dashrightarrow C \times C$. This idea was inspired by \cite[page 902]{naranjopirola2018}. 
We provided some additional details of the argument for convenience of the reader.
\end{remarks}
 \subsection{Extensions of abelian varieties} \label{section:extensionnodal}
 
Let $A$ be an abelian variety with dual abelian variety $A^\vee$.  
Recall that, by the Barsotti--Weil formula, there is a canonical isomorphism 
$ 
\Ext(A, \G_m) = A^\vee$. 
If $T \cong \G_m^r$ is a torus, $G$ a connected commutative algebraic group, and $0\to  T \to G \to A \to 0$ an exact sequence of commutative algebraic groups, then by applying $\Hom(-, \G_m)$ one obtains a homomorphism 
$$
c^t \colon \Hom(T,\G_m) \longrightarrow \Ext(A, \G_m) 
= A^\vee$$
from the character group of $T$ to the dual abelian variety $A^\vee$, and this construction induces a bijection 
(compare \cite[Proposition 2]{carlson-extensions} and \cite[Chapter II, Section 2]{chai-siegelcompactifications}):
\begin{align} \label{equation:extAVdual}
\Ext(A,T) = \Hom\left( \Hom(T,\G_m), A^\vee \right). 
\end{align}

\begin{lemma} \label{lemma:extension-class-commutativity}
For $i \in \set{1,2}$, let $0 \to T_i \to G_i \to A_i \to 0$ be an exact sequence of algebraic groups, where $A_i$ is an abelian variety and $T_i\cong \G_m^{r_i}$ a torus. 
Let $f \colon G_1 \to G_2$ be a morphism of algebraic groups. 
Then $f$ restricts to a homomorphism $f|_{T_i} \colon T_1 \to T_2$ and hence induces a homomorphism $\bar f \colon A_1 \to A_2$. 
Moreover, if $c^t_i \in \Hom(\Hom(T_i,\G_m), A^\vee_i)$ is the homomorphism that corresponds to the class of $G_i$ in $\Ext(A_i,T_i)$ via \eqref{equation:extAVdual}, then the following diagram commutes:
\begin{align}\label{diagram:character-dualAV}
\begin{split}
\xymatrixcolsep{3pc}
\xymatrix{
\Hom(T_2,\G_m) \ar[d]^{c^t_2} \ar[r]^{(f|_{T_1})^\ast} & \Hom(T_1,\G_m) \ar[d]^{c^t_1} \\
A_2^\vee \ar[r]^{(\bar f)^\vee} & A_1^\vee. 
}
\end{split}
\end{align}
\end{lemma}
\begin{proof}
The first statement follows from the fact that $\Hom(T_1,A_2) = 0$. For the second statement, consider the following commutative diagram with exact rows:
    \[
    \xymatrix{
    0 \ar[r] & T_1 \ar[r] \ar[d] & G_1 \ar[r] \ar[d] & A_1 \ar[d]\ar[r] & 0 \\
    0 \ar[r] & T_2 \ar[r] & G_2 \ar[r] & A_2 \ar[r] & 0.
    }
    \]
Applying $\Hom(-,\G_m)$ to this diagram, and using the fact that for an abelian variety $A$, the isomorphism $\Ext(A,\G_m) = A^\vee$ is functorial in $A$, the commutativity of \eqref{diagram:character-dualAV} follows. 
\end{proof}

\subsection{Extension classes of nodal curves} \label{subsection:extensionclassesofnodalcurves}

Next, we recall some known results on extensions of Jacobians of smooth projective connected curves, following \cite[Sections 2.2--2.4]{alexeev-compactifiedjacobians}.  We will make use of the following definition. 
\begin{definition} \label{definition:orientation}
Let $X$ be a connected nodal curve, and let $ \Gamma(X)$ be its dual graph. An \emph{orientation} of $\Gamma(X)$ is the choice of an ordering $(P^+, P^-)$ on every pair of points on the normalization $\wt X$ of $X$ lying above the same node.
\end{definition}
Let $X$ be a connected nodal curve. Let $X_1, \dotsc, X_n$ be the irreducible components of $X$, and $ \widetilde X_i \to X_i$ their normalizations. 
We denote the dual graph of $X$ by $\Gamma = \Gamma(X)$ and fix any orientation of $\Gamma$.
Let $$ \wt X = \coprod_{i = 1}^n \wt X_i \longrightarrow X$$ be the normalization of $X$. The dual abelian variety of $J \wt X$ is identified with itself via the principal polarization, and the character group of the torus $T = \Ker(JX \to J\wt X)$ is canonically identified with $H_1(\Gamma,\Z)$.
Hence, 
\begin{align} \label{equation:nodalidentification}
\Ext(J\wt X, T) = \Hom(\Hom(T,\G_m), J\wt X) = \Hom(H_1(\Gamma,\Z), J\wt X).
\end{align}
Moreover, the homomorphism 
\begin{align} \label{equation:ct}
c^t \colon H_1(\Gamma,\Z) \longrightarrow J\wt X
\end{align}
corresponding to $[JX] \in \Ext(J \wt X,T)$ via \eqref{equation:nodalidentification} is described explicitly in the following way. 
Every edge $e$ of $\Gamma$ corresponds to a node $P$ of $X$, and the orientation defines an ordered pair of points $(P^{+}, P^{-})$ on $\wt X$.
Put
\[
c^t(e) = P^+ - P^{-} \in \Pic(\wt X), 
\]
and extend this by linearity to the free module $C_1(\Gamma,\Z)$ on the edges of $\Gamma$. Let $C_0(\Gamma,\Z)$ be the free $\Z$-module on the set of vertices of $\Gamma$. For an edge $e$, let $\text{end}(e)$ and $\text{beg}(e)$ be the end-vertex and begin-vertex of $e$, as determined by the orientation of $e$, and define
\[
\partial \colon C_1(\Gamma,\Z) \longrightarrow C_0(\Gamma,\Z) \quad \text{by} \quad \partial(e) = \text{end}(e) - \text{beg}(e) \quad \text{for an edge } e \text{ of } \Gamma.
\]
If
$h \in H_1(\Gamma,\Z) = \Ker(\partial \colon C_1(\Gamma,\Z) \to C_0(\Gamma,\Z)),$ one has $c^t(h)  \in J\wt X$, and this construction defines the homomorphism \eqref{equation:ct}, see \cite[Section 2.4]{alexeev-compactifiedjacobians}. 

\subsection{Graph homology and extension classes.} We continue with the notation of Section \ref{subsection:extensionclassesofnodalcurves}. 
Let $v_1, \dotsc, v_m$ be a set of vertices and $e_1, \dotsc, e_{m}$ a set of edges of $\Gamma$, such that for each $j \in \set{1, \dotsc, m-1}$, the edge $e_j$ connects the vertices $v_j$ and $v_{j+1}$, and $e_m$ connects $v_m$ and $v_1$. In particular, if $m = 1$, then $e_1$ is a loop connecting $v_1$ to itself. 

Let $$\wt X_{i_1}, \wt X_{i_2}, \dotsc , \wt X_{i_m}$$ be the connected components of the normalization $\wt X$ of $X$ that correspond to the vertices $v_1, \dotsc, v_m$. For $j \in \set{1, \dotsc, m-1}$, the orientation of $e_j$ defines an ordered pair of points $(P_j^+, P_j^-)$ such that $P_j^+$ lies either on $\wt X_{i_j}$ or on $\wt X_{i_{j+1}}$, and the opposite is true for $P_j^-$; define $\varepsilon_j \in \set{1,-1}$ by declaring that $\varepsilon_j = 1$ if $P_j^-$ lies on $\wt X_{i_j}$ and $\varepsilon_j = -1$ otherwise. 
We obtain a cycle 
\begin{align} \label{equation:homologycycle}
\gamma = \sum_{j = 1}^m \varepsilon_j \cdot e_j \in C_1(\Gamma,\Z),
\end{align}
and one readily observes that $\partial(\gamma) = 0$, so that $\gamma\in H_1(\Gamma,\Z) \subset C_1(\Gamma,\Z)$. 

\begin{lemma} \label{lemma:basis:homology:suitable}
Let $X$ be a connected nodal curve with dual graph $\Gamma$. 
Then the following holds.
\begin{enumerate}
    \item There exists a linearly independent subset $S = \set{\gamma_1, \dotsc, \gamma_k} \subset H_1(\Gamma,\Z)$ consisting of homology classes of the form \eqref{equation:homologycycle} such that $S$ defines a basis of $H_1(\Gamma,\Q)$. \label{item:lemma:basis:homology:suitable}
    \item \label{item:lemma:alexeev1motif}
    Let $\wt X_1, \dotsc, \wt X_n$ be the connected components of the normalization $\wt X$ of $X$. 
Let $\gamma \in H_1(\Gamma,\Z)$ be a class of the form \eqref{equation:homologycycle}, and consider the homomorphism $c^t \colon H_1(\Gamma,\Z) \to J \wt X$, see \eqref{equation:nodalidentification} and \eqref{equation:ct}. 
There are points $p_{i}, q_i \in \wt X_i$ for each $i \in \set{1, \dotsc, n}$, such that
\begin{align} \label{equation:formct}
c^t(\gamma) = (p_{1} - q_{1}, \dotsc, p_{n} - q_{n}) \in J\wt X_1 \times \cdots \times J \wt X_n = J \wt X. 
\end{align} 
Note that we do not require that all the points $p_{i}$ and $q_{i}$ are distinct. 
\end{enumerate}
\end{lemma}

\begin{proof}
To prove item \eqref{item:lemma:basis:homology:suitable}, recall first the following fact. Let $Y$ be a path-connected one-dimensional CW complex with basepoint $y_0$, 
a $0$-cell. Then every loop in $Y$ is homotopic to a loop consisting of a finite sequence of
edges traversed monotonically, see \cite[Section 1.1, Exercise 19]{hatcher-AT}. Consequently, homology classes of the form \eqref{equation:homologycycle} generate $H_1(\Gamma,\Z)$, yielding the lemma. 
(We do not ask that $S$ is an integral basis, because not any generating set of a free $\Z$-module contains a basis.) 

To prove item \eqref{item:lemma:alexeev1motif}, we note that in the notation used above equation \eqref{equation:homologycycle}, we have
\[
c^t(\gamma) = \left( 
P^+_{i_1} - P^-_{i_{m-1}}, P^+_{i_2} - P^-_{i_1}, P^+_{i_3} - P^-_{i_2}, \dotsc, P^+_{i_{m}} - P^-_{i_{m-1}}
\right) 
\in \prod_{j = 1}^m J \wt X_{i_j} \subset J \wt X.
\]
The lemma follows.   
\end{proof}

\subsection{Extension classes of nodal hyperelliptic curves.}

We turn to the hyperelliptic case.  

\begin{lemma} \label{lemma:hyper-extension}
    Let $X$ be an irreducible one-nodal hyperelliptic curve of genus $g \geq 1$, with normalization $\wt X \to X$. Consider the homomomorphism $c^t \colon H_1(\Gamma,\Z) \to J\wt X$, where $\Gamma$ denotes the dual graph of $X$, see equation \eqref{equation:ct}. If $\gamma$ is a generator for $H_1(\Gamma,\Z) \cong \Z$, there is a point $x \in \wt X$ such that $c^t(\gamma) = x - \iota(x)$, where $\iota \colon \wt X \to \wt X$ is the hyperelliptic involution. 
\end{lemma}

\begin{proof}
By the description of the map $c^t$, we have $c^t(\gamma) = p-q$, where $p,q$ denote the points on the normalization $\wt X$ of $X$ that are glued to form the nodal curve $X$.
As $X$ is hyperelliptic, the hyperelliptic involution on $\wt X$ descends to an involution on $X$, which implies that $q=\iota(p)$.
\end{proof}

\begin{lemma} \label{lemma:compactabeliansubvarieties} 
Let $X$ be a very general one-nodal hyperelliptic curve. 
There is no positive dimensional abelian subvariety of $JX$. 
 \end{lemma}
 \begin{proof}
 Let $T = \Ker(JX \to J\wt X)$ and consider the extension
  \begin{align} \label{equation:extension-isogeny-split}
 0 \longrightarrow T \longrightarrow JX \longrightarrow J \wt X \longrightarrow 0 .
 \end{align}
 Since $X$ is very general, $J \wt X$ is simple.
 Hence any positive dimensional abelian subvariety of $JX$ must be isogenous to $J\wt X$.
 If such a subvariety exists, then the extension \eqref{equation:extension-isogeny-split}  splits up to isogeny, and so it suffices to exclude the latter.
By Lemma \ref{lemma:hyper-extension}, the isomorphism $\Ext(J \wt X, T) \cong J \wt X$ (cf.~Section \ref{section:extensionnodal}) identifies $[JX]$ with $x-\iota(x)$ for some $x \in \wt X$, hence \eqref{equation:extension-isogeny-split} splits up to isogeny if and only if $x-\iota(x) \in J \wt X$ is torsion. As the one-nodal hyperelliptic curve $X$ is very general, this is not the case.
 \end{proof}

\subsection{Nodal degeneration and vanishing cycles} \label{subsection:nodaldegenerationandcohomology}

Recall the following result. 

\begin{proposition} \label{proposition:fibrewise-retraction} 
    Let $D$ be the open unit disc with origin $0 \in D$. Let $X$ be a complex analytic space and $f \colon X \to D$ a proper map. Put $X_0 = f^{-1}(0)$. If $X$ is smooth, and $f$ is smooth over $D - \set{0}$, then the inclusion $X_0 \hookrightarrow X$ is a homotopy equivalence. If $X$ is any complex analytic space and $f$ any proper map $X \to D$, then the same is true up to shrinking $D$ around $0$. 
\end{proposition}
\begin{proof}
See for instance \cite[Proposition C.11 and Remark C.12.ii]{petersteenbrink}.
\end{proof}

\begin{lemma}\label{lemma:dejong}
    Let $X$ be an analytic space, $(D, 0)$ the pointed unit disc, and $f \colon X \to D$ a family of nodal curves over $D$ which is smooth over $D^\ast = D - \set{0}$. There exists a complex manifold $\widetilde X$ and a projective morphism $h \colon \widetilde X \to X$ which is an isomorphism over the regular locus of $X$, such that $f \circ h$ defines a family of nodal curves $ \widetilde X \to D$ which is smooth over $D^\ast$. 
\end{lemma}
\begin{proof}
See \cite[Lemma 3.2]{dejong-alterations}; the main point being that we can find  a resolution $\widetilde X$ of $X$ such that the fibres of $\widetilde X\to D$ are reduced.
\end{proof}
Consider a proper holomorphic map $$f\colon\ca X\longrightarrow D $$ from an $n$-dimensional complex manifold $\ca X$ to a disc $D$. Assume that $f$ is a submersion over the punctured disc $D^\ast$, and that over $0 \in D$, there are $k$ critical points $x_{1}, \dotsc, x_{k} \in X_0 =  f^{-1}(0)$ for some $k \in \Z_{\geq 1}$, and that these are non-degenerate. Assume that $\ca X$ is regular (something we can always achieve by modifying $\ca X$, see Lemma \ref{lemma:dejong}). 

\begin{lemma} \label{lemma:voisin-i}
Continue with the above notation, and let $t \in D^\ast$. There are $k$ disjoint spheres $S_1^{n-1}, \dotsc, S_k^{n-1} \subset X_t$ and
    a deformation retraction of $\ca X$ onto the union of $X_t$ and $k$ disjoint $n$-dimensional balls $B^n_1, \dotsc, B^n_k$, where the ball $B^{n_i}$ is glued to $X_t$ along the sphere $S^{n-1}_i \subset X_t$. 
\end{lemma}
\begin{proof}
    This is a straightforward generalization of \cite[Theorem 2.16]{voisin-II}. 
\end{proof}

\begin{corollary} \label{corollary:voisin-ii}
Continue with the above notation. Let $t \in D^\ast$ and let $i:X_t \hookrightarrow \ca X$ be the inclusion. Then $i_\ast \colon H_m(X_t,\Z) \to H_m(X,\Z)$ is an isomorphism for $m < n-1$. For $m = n-1$, the map $i_\ast$ is surjective, with kernel generated by the cohomology classes 
of the spheres $S^{n-1}_1, \dotsc, S^{n-1}_k \subset X_t$. 
\end{corollary}

\begin{proof} 
    This is a straightforward generalization of \cite[Corollary 2.17]{voisin-II}. 
\end{proof}

\begin{lemma} \label{lemma:voisin:PL}
Continue with the above notation and let $t \in D^\ast$. 
Let $\delta_1, \dotsc, \delta_n \in H^{n-1}(X_t,\Z)$ be the Poincar\'e duals of the homology classes of the vanishing spheres $S_1^{n-1}, \dotsc, S_k^{n-1} \subset X_t$, see Lemma \ref{lemma:voisin-i}. Let $(-,-) \colon H^{n-1}(X_t,\Z) \times H^{n-1}(X_t,\Z) \to \Z$ be the cup-product pairing. For some $\epsilon_n \in \set{\pm 1}$, depending only on $n$, the natural generator $T \in \Aut(H^1(X_t,\Z))$ of the monodromy group satisfies
    \[
    T(\alpha)  = \alpha + \epsilon_n \cdot \sum_{i = 1}^k \left(\alpha, \delta_i\right) \delta_i \quad \quad  \forall \alpha \in H^{n-1}(X_t,\Z). 
    \]
\end{lemma}
\begin{proof}
        This is a straightforward generalization of \cite[Theorem 3.16]{voisin-II}.
\end{proof}
Next, we verify that, for a nodal degeneration of curves over a disc with smooth general fibre, the monodromy invariant part of the first integral cohomology group of the general fibre does not change after any finite base change. 

\begin{lemma}
\label{lemma:family:3} 
Let $D \ni 0$ be the pointed unit disc. 
Let
$
\ca X\to D
$
be a family of nodal curves, smooth over $D^\ast = D \setminus \set{0}$.  
Let $\tau\colon D' = D\to D$ be the map $z \mapsto z^m$, let $0 ' \in D'$ be the preimage of $0 \in D$, 
and consider the base change $\ca X'\coloneqq \ca X \times_{D} D'$. 
Let $t \neq 0 \in D$ and fix a preimage $t'\in D'$ of $t\in D$. Let $T \in \Aut(H^1(X_t,\Z))$ and $T' \in \Aut(H^1(X'_{t'},\Z)$ be the monodromy operators induced by the restrictions of the families to $D^\ast$ and ${D'}^\ast$. Then the invariant subspaces of $T$ and $T'$ coincide, that is,
$$
H^1(X_t,\Z)^{T}=H^1(X'_{t'},\Z)^{T'} .
$$
\end{lemma}
\begin{proof}
To prove the lemma, we may assume that $\ca X$ is regular, see Lemma \ref{lemma:dejong}. Suppose that the central fibre $X_0$ has $k$ nodes, and let $\delta_1, \dotsc, \delta_k \in H^1(X_t,\Z)$ be the cohomology classes attached to the vanishing spheres $S_1, \dotsc, S_k \subset X_t$, see Lemma \ref{lemma:voisin-i}.  
By Lemma \ref{lemma:voisin:PL}, we have the following generalization of the Picard--Lefschetz formula: for each $\alpha \in H^1(X_t,\Z)$, one has 
\[
T(\alpha) = \alpha + \sum_{i = 1}^k (\alpha \cdot \delta_i)\delta_i \in H^1(X_t,\Z). 
\]
On the one hand, the monodromy operator $T'$ on $H^1(X'_{t'},\Z)=H^1(X_{t},\Z)$ satisfies $T'=T^{m}$. On the other hand, we have 
\[
T^m(\alpha) = \alpha + m \cdot \sum_{i = 1}^k (\alpha \cdot \delta_i)\delta_i \in H^1(X_t,\Z), \quad \alpha \in H^1(X_t,\Z). 
\]
This implies that, for $\alpha \in H^1(X_t,\Z)$, one has:
\[
T'(\alpha) = \alpha \quad \iff \quad m \cdot \sum_{i = 1}^k (\alpha \cdot \delta_i)\delta_i = 0 \quad \iff \quad \sum_{i = 1}^k (\alpha \cdot \delta_i)\delta_i = 0 \quad \iff \quad  T(\alpha) = \alpha. 
\]
This proves the lemma.  
\end{proof}
 
\begin{remark} \label{rem:item:family:4}
In the course of the proof of Theorem \ref{theorem:maintheorem1:hyperelliptic} we will be forced to  perform various base changes. This is a priori a subtle issue for the following reasons. 
    We plan to degenerate to different nodal fibres, which is equivalent to the degeneration to one fixed nodal fibre followed by the application of a monodromy operator.
    However, the monodromy action on cohomology with finite coefficients may become trivial after a base change, while the basic criterion in Lemma \ref{lemma:isogenyVHS} that we aim to exploit is in fact equivalent to the analogous assertion for (sufficiently divisible) finite coefficients. 
    For this reason, the fact that even after an arbitrary base change one can deduce additional information from degeneration to various nodal fibres (a fact which was already exploited in \cite{bardellipirola-curvesofgenusg,naranjopirola2018}), seems somewhat surprising.
    The key reason which makes these arguments work is given in Lemma \ref{lemma:family:3}  above, which says that the monodromy invariant subspace is not affected by any finite base change. 
\end{remark}

We conclude the section with the following lemma. 

\begin{lemma} \label{lemma:basechange-surjective-monodromy}
Let $\ca C$ be complex analytic space, $D \subset \C$ an open disc around $0 \in \C$, and $q \colon \ca C \to D$ a family of nodal curves over $D$. Suppose that for each $s \in D^\ast = D - \set{0}$, the curve $C_s = q^{-1}(s)$ is of compact type. Then up to shrinking $D$ around $0$, the following holds: for $t \in D^\ast$ and $T \in \Aut(H^1(C_t,\Z))$ a generator of the monodromy group, the natural map 
\[
H^1(\ca C,\Z) \longrightarrow H^1(C_t,\Z)^T
\]
is surjective. 
\end{lemma}

\begin{proof}
Let $\mathcal C_i$ ($i\in I$) denote the irreducible components of $\ca C$. 
By flatness of $q$, the induced morphism $q_i \colon \ca C_i \to D$ is surjective.
We claim that the general fibre  of $q_i$ is irreducible.
Clearly, it suffices to prove this after shrinking the disc $D$.
The main point is then that $\mathcal C\to D$ is a family of nodal curves, hence the fibre above $0$ is reduced.
This in turn implies that for each irreducible component $C_{0,j}$ of the special fibre of $q$, there is up to shrinking a section of $q$ that passes through a general point of $C_{0,j}$.
This shows that the component $\mathcal C_i$ of $\mathcal C$ that contains $C_{0,j}$ admits a section and hence has irreducible general fibre.
Running through all components $C_{0,j}$ of the special fibre of $q$, we get this way that each $q_i \colon \ca C_i \to D$ has irreducible general fibre.

Since the general fibre of $q\colon\mathcal C\to D$ is of compact type, the general fibre of $q_i$ is smooth and the index set $I$ forms the vertices of a tree that indicates which components of a general fibre of $q$ meet. If $i,j\in I$ are joined by an edge $e_{i,j}$, or equivalently, if the general fibres of $\mathcal C_i$ and $\mathcal C_j$ are glued at a point, then we get a section $e_{i,j}:D^*\to\mathcal C$, which has to extend across the puncture by properness of $q$.
Since the arithmetic genus of the fibres of each $q_i\colon \mathcal C_i\to D$ is constant, and because the same holds for $q\colon \mathcal C\to D$, we see that $\mathcal C$ is given by the quotient 
$$
\mathcal C=\p{\coprod_{i\in I} \mathcal C_i}/\sim
$$
where we glue for each edge $e_{i,j}$ between some indices $i,j\in I$ according to the section constructed above.
Since the fibres of $q$ are nodal, we see moreover that the points that are glued via $\sim$ on the special fibre lie in the smooth locus of $\bigsqcup_{i\in I} \mathcal C_i\to D$.

A simple Mayer--Vietoris argument now reduces us to show that
\[
H^1(\ca C_i,\Z) \longrightarrow H^1(C_{it},\Z)^T
\]
is surjective for each $i\in I$.
In other words, we have reduced the result to the case where $q$ is smooth over the punctured disc $D^*$.
In this case we apply Lemma \ref{lemma:dejong} and get a modification $\tau:\widetilde {\mathcal C}\to \mathcal C$ given by successive blow-ups of the singular points in the central fibre, such that the fibres of $\widetilde {\mathcal C}\to D$ are reduced and hence nodal curves.
By \cite[Theorem 7.8]{kollar-plurigenera}, $
\tau_\ast \colon \pi_1(\wt{\ca C}) \to \pi_1(\ca C)
$
is an isomorphism.
Passing to the abelianization and applying $\Hom(-,\Z)$, we find that  $\tau^\ast\colon H^1(\mathcal C,\Z)\to H^1(\widetilde{\mathcal C},\Z) $ is an isomorphism.
This reduces us to the case where $\mathcal C$ is regular and $q$ is smooth over $D^*$. 
By Corollary \ref{corollary:voisin-ii},  $H_1(C_t,\Z) \to H_1(\ca C,\Z)$ is surjective.
It follows that $H^1(\ca C,\Z) \to H^1(C_t,\Z)$ is injective with torsion-free cokernel.
By the local invariant cycle theorem (cf.\ \cite{morrison-clemensschmid}), the map $H^1(\ca C,\Z) \to H^1(C_t,\Z)^T$ becomes surjective after tensoring with $\Q$; as its cokernel is torsion-free, it is surjective.
This concludes the proof of the lemma.
\end{proof}

\subsection{Degenerations of hyperelliptic curves} In this section, we construct a family of stable hyperelliptic curves satisfying suitable properties. The base will be higher dimensional with several divisors each of which giving rise to a family of one-nodal hyperelliptic curves, allowing us to degenerate a very general hyperelliptic curve in different directions.
We will ultimately need these different degenerations in order to prove Theorem \ref{theorem:maintheorem1:hyperelliptic}, see Section \ref{section:maintheorem}. 

The following lemma is certainly well-known; we include some details for convenience of the reader. 

\begin{lemma} \label{lemma:hyperelliptic-degeneration} 
Let $n, g \geq 2$ be positive integers with $n \leq g$. Consider the affine space $\A^{2g-1}$ with coordinates $z_1, \dotsc, z_{2g-1}$.  
There exists a non-empty Zariski open subset $U \subset \A^{2g-1}$, irreducible divisors $\Delta_i \subset U$ for $i = 1, \dotsc, n$, and a family of genus $g$ stable hyperelliptic curves
\begin{align} \label{family:stable-hyper}
p \colon \ca X\longrightarrow U
\end{align}
such that the following holds.
\begin{enumerate}
    \item \label{item:family:1} Let $\Delta \coloneqq \cup_i \Delta_i$. The family \eqref{family:stable-hyper} is smooth over $U - \Delta$, the fibre $X_t = p^{-1}(t)$ for very general $t \in U$ is a very general hyperelliptic curve of genus $g$, and for each $i \in \set{1,\dotsc, n}$, the fibre $X_{0_i}$ above a very general point $0_i \in \Delta_i$ is a very general one-nodal hyperelliptic curve of arithmetic genus $g$.  
      \item For $0_i \in \Delta_i$ $(i\in \set{1, \dotsc, n})$, there exists a one-dimensional disc $D_i \subset U$ intersecting $\Delta$ transversally in $0_i \in \Delta_i$ such that the restriction $p|_{D_i} \colon \ca X|_{D_i} \to D_i$ is a Lefschetz degeneration with nodal central fibre above $0_i \in D_i$. \label{item:family:2:Lefschetz}
      \item \label{item:family:3:symplectic}
      Let $t \in U- \Delta$ and $t_i \in D_i - \set{0_i}$ be base points. 
      Let $\delta_i\in H^1(X_{t_i},\Z)$ be the vanishing cycle associated to $p|_{D_i} \colon \ca X|_{D_i} \to D_i$, and view $\delta_i$ as an element of $H^1(X_t,\Z)$ via parallel transport along a path $\rho_i$ from $t$ to $t_i$. 
      Then $\delta_1,\dots, \delta_n$ can be completed to a symplectic basis 
\begin{align} \label{equation:symplecticbasis}
H^1(X_t,\Z) = \langle \delta_1, \dotsc, \delta_g; \gamma_1, \dotsc, \gamma_g \rangle .
\end{align}  
\end{enumerate}
\end{lemma}
\begin{proof}
Let $a_{n-1},\dots ,a_{2g-1}\in \C$ be general complex numbers, and consider the following equation:
$$
y^2=((x-a_{n-1})^2-z_1)\cdot ((x-a_n)^2-z_2)\cdot\prod_{i=1}^{n-2}((x-z_{n+i}-a_{n+i})^2-z_{i+2})  \cdot \prod_{j=2n-1}^{2g-1}(x-z_j-a_j).
$$
(The slight asymmetry in the quadratic terms reflects automorphisms of $\CP^1$ and stems from the fact that we want to have a family that depends on $2g-1$ parameters $z_1,\dots ,z_{2g-1}$.) 
This defines a family of affine hyperelliptic curves over $\A^{2g-1}$, branched at the $2g+2$ points
$$
x=\infty,\ \ x=\pm \sqrt{z_1} + a_{n-1},\ \ x=\pm \sqrt{z_2}+a_n,\ \ x=\pm\sqrt{z_{i+2}}+z_{n+i}+a_{n+i},\  \ x=z_j+a_j
$$
for $i=1,\dots ,n-2$ and $j=2n-1,\dots ,2g-1$.
This extends to a projective family of hyperelliptic curves and we denote by $U\subset \A^{2g-1}$ the open subset where the corresponding hyperelliptic curve attains at most one node.
The corresponding projective family of hyperelliptic curves is denoted by $p \colon \ca X \to U$ and we note that $\ca X$ is regular. By construction, for $u = (z_1, \dotsc, z_{2g-1}) \in U$, the fibre $X_u = p^{-1}(u)$ is either smooth or attains exactly one node; the latter happens if and only if $z_i = 0$ for some $i \in \set{1, \dotsc, n}$.   
Let $\Delta_i \coloneqq U\cap \{z_i=0\}$ for $i\in \set{1,\dots ,n}$.
The above description of the ramification points of the hyperelliptic covering $X_u\to \CP^1$ for $u\in U$ shows that the moduli map $\Delta_i\to \overline{\mathcal M}_g$ is generically  finite onto its image (because the hyperelliptic locus in $\mathcal M_g$ has dimension $2g-1$). 
Altogether this proves item \eqref{item:family:1} in the lemma.

For $i\in \set{1, \dotsc, n}$, let $D_i \subset U$ be a disc that intersects $\Delta$ transversally in a general point $0_i \in \Delta_i$.
Up to shrinking the disc $D_i$, we can assume that the total space $\ca X|_{D_i}$ of the restriction $\ca X|_{D_i} \to D_i$ is regular: the only possible singularity is at the node of the central fibre, where analytically locally an equation of $\mathcal X|_{D_i}$ is given by $t=x^2-y^2$, which yields a regular surface. 
Thus, $\ca X|_{D_i} \to D_i$ is a Lefschetz degeneration, proving item \eqref{item:family:2:Lefschetz}.

The fibre $X_t$ is a double covering $X_t\to \CP^1$ branched along $2g+2$ points $p_0,p_1,\dots ,p_{2g+1}$.
By construction, $\mathcal X\to U$ is given by an equation of the form $y^2=\prod_{i=1}^{n}f_i(x)\cdot \prod_{j=2n-1}^{2g-1} g_j(x)$, where $g_j(x)=x-z_j-a_j$ and $f_i(x)$ is quadratic in $x$. 
Up to reordering, we can assume that $\{p_{2i-1},p_{2i}\}$ corresponds to the roots of $f_i(x)$.
In particular, $p_{2i-1}$ and $p_{2i}$ collide along the Lefschetz degeneration over the disc $D_i$, for $i=1,\dots ,n$.

We pick a path of shortest distance between $p_{2i-1}$ and $p_{2i}$ on $\CP^1$ for $i=1,\dots ,g$ and note that the preimage of this path in the hyperelliptic curve $X_t$ gives  rise to a homology class in $H_1(X_t,\Z)$ (well-defined up to sign) whose Poincar\'e dual $\delta_i\in H^1(X_t,\Z)$ is for $i=1,\dots ,n$ the vanishing cycle that corresponds to colliding $p_{2i-1}$ and $p_{2i}$. 
The classes $\delta_1,\dots ,\delta_g$ are orthogonal to each other and can be completed to a symplectic bases, proving item \eqref{item:family:3:symplectic}. 
\end{proof}

Let $p \colon \ca X \to U$ be a family of hyperelliptic curves of genus $g$ as in Lemma \ref{lemma:hyperelliptic-degeneration}. For $i \in \set{1, \dotsc, n}$, consider the embedding
$
H^1(X_{0_i},\Z) \hookrightarrow H^1(X_t,\Z)
$ defined as the composition of the inverse of the map $H^1(\ca X|_{D_i},\Z) \to H^1(X_{0_i},\Z)$ (which is an isomorphism by Proposition \ref{proposition:fibrewise-retraction}), the restriction $H^1(\ca X|_{D_i},\Z) \to H^1(X_{t_i},\Z)$ and the parallel transport $ H^1(X_{t_i},\Z)\to  H^1(X_{t},\Z) $ along $\rho_i$. 
Let $$W_0H^1(X_{0_i},\Z) = W_0H^1(X_{0_i},\Q) \cap H^1(X_{0_i},\Z)$$ be the integral part of the zeroth piece of the weight filtration. 
For $i \in \set{1, \dotsc, n}$, let  
 $T_i \in \Aut(H^1(X_t,\Z))$ be the monodromy operator associated to the path $\rho_i$ and the pointed disc $(D_i, 0_i)$.

\begin{lemma} \label{lemma:conditionssatisfied}
Consider the above notation. 
 With respect to the symplectic basis \eqref{equation:symplecticbasis}, we have 
\begin{align} \label{equation:invariant-image-lemma-statement-I}
\begin{split}
\Ima\left(
H^1(X_{0_i},\Z) \hookrightarrow H^1(X_t,\Z)
\right) &= H^1(X_t,\Z)^{T_i} \\
&= \langle \delta_1, \dotsc, \delta_g; \gamma_1, \dotsc, \gamma_{i-1}, \widehat\gamma_i, \gamma_{i+1}, \dotsc, \gamma_g \rangle,
\end{split}
\\
\label{equation:invariant-image-lemma-statement-II}
\Ima\left(W_0H^1(X_{0_i},\Z) \hookrightarrow H^1(X_t,\Z)
\right)&= \Ima\left(T_i - \id\right) = \left\langle \delta_i   \right \rangle = \Z \cdot \delta_i. 
\end{align}
Here, the module on the right in \eqref{equation:invariant-image-lemma-statement-I} denotes the submodule of $H^1(X_t,\Z)$ obtained from $H^1(X_t,\Z)$ by removing $\gamma_i$ from the symplectic basis \eqref{equation:symplecticbasis}.
\end{lemma} 

\begin{proof}
 For each $i \in \set{1, \dotsc, n}$, define \begin{align*}H_i' \coloneqq \Ima\left(
H^1(X_{0_i},\Z) \hookrightarrow H^1(X_t,\Z)
\right) \quad \text{and}\quad V_i \coloneqq \Ima\left(W_0H^1(X_{0_i},\Z) \hookrightarrow H^1(X_t,\Z)
\right).\end{align*}
 Consider the monodromy operator 
$
T_i \colon H^1(X_t,\Z) \to H^1(X_t,\Z)$. In view of items \eqref{item:family:2:Lefschetz} and \ref{item:family:3:symplectic} in Lemma \ref{lemma:hyperelliptic-degeneration}, $T_i$ is given by the formula 
 \begin{align*}
 T_i(\alpha)=\alpha + (\alpha\cdot \delta_i)\delta_i, \quad \alpha \in H^1(X_t,\Z),
 \end{align*}
 where the $\delta_i$ are the vanishing cycles $\delta_i \in H^1(X_t,\Z)$ attached to the Lefschetz degenerations $\ca X|_{D_i} \to D_i$ and the path $\rho_i$. 
 It is well-known that 
\begin{align} \label{align:rational-cohomology}
V_i \otimes \Q = \Ima(T_i - \id) \otimes \Q \subset\Ker(T_i - \id) \otimes \Q =  H_i' \otimes \Q \subset H^1(X_t,\Q). 
\end{align}
In fact, \eqref{align:rational-cohomology} holds integrally because the respective spaces are saturated in $H^1(X_t,\Z)$.
In particular, \eqref{equation:invariant-image-lemma-statement-I} follows. 
 Note that $(T_i - \id)(\alpha) = (\alpha \cdot \delta_i)\delta_i$ for each $\alpha \in H^1(X_t,\Z)$. 
 As we have $\gamma_i \cdot \delta_i = -1$, this gives
$
\left(T_i - \id \right)(\gamma_i) = (\gamma_i \cdot \delta_i) \delta_i = -\delta_i.
$
Therefore, 
$
V_i = \Ima\left(T_i - \id\right) = \left\langle \delta_i   \right \rangle = \Z \cdot \delta_i \subset H_i'$, proving \eqref{equation:invariant-image-lemma-statement-II}, and thereby the lemma.
\end{proof}

\section{Moving the extension class} \label{section:movingextension}

Let $f \colon X \to Y$ be a flat morphism of complex analytic spaces with reduced fibres. Following \cite[Definition 5.3]{chiang-lipman-simultaneous}, a \emph{simultaneous normalization} of $f$ is a finite morphism of analytic spaces $\nu \colon Z \to X$ such that $\bar f \coloneqq f \circ \nu$ is a flat morphism $Z \to Y$ whose non-empty fibres are normal, and such that for each $y \in f(X)$, the induced morphism of fibres $\nu_y \colon \bar f^{-1}(y) \to f^{-1}(y)$ is a normalization map.  
We call $f \colon X \to Y$ \emph{equinormalizable} if a simultaneous normalization of $f$ exists. If $Y$ is normal, $f \colon X \to Y$ is flat with reduced fibres, and the connected components of $X$ are equidimensional, then any simultaneous normalization $\nu \colon Z \to X$ of $f$ is a normalization of $X$, see \cite[Proposition 5.4]{chiang-lipman-simultaneous}; in particular, $\nu \colon Z \to X$ is then unique up to isomorphism. 
 
\begin{proposition} \label{prop:simultaneous-normalization}
Let $X$ and $Y$ be complex analytic spaces with $Y$ normal. 
Let $f \colon X \to Y$ be a (proper) family of nodal curves with irreducible fibres (cf.\ Section \ref{section:conventions}). Assume the number of nodes of the curve $X_y = f^{-1}(y)$ is constant for $y \in Y$. Then $f$ is equinormalizable. 
\end{proposition} 

\begin{proof}
    This can for instance be deduced from \cite[Theorem 5.6, Corollary 5.4.2, and Definition 5.1]{chiang-lipman-simultaneous}; we include some details for convenience of the reader.
    Since $f$ is a proper family of nodal curves, $\Sing(f)\to Y$ is finite and unramified, see e.g.\ \cite[\S 2.21]{dejong-alterations}.
    Since the number of nodes is constant in the family, it follows that each component of $\Sing(f)$ dominates $Y$.
    Using this we see that for each $x\in \Sing(f)$ with image $y=f(x)$, there is a suitable neighbourhood $U\subset X$ such that $\Sing(f)\cap U\cap f^{-1}(y)=\{x\}$. 
    Up to shrinking $U$, we can assume that the fibres of $U\to Y$ have at most one node.
    Up to replacing $U$ by the intersection with the preimage of a suitable neighbourhood of $y$ in $Y$, we can moreover assume that $U\cap \Sing(f)\to f(U)$ is finite and surjective.
    Moreover, since $f$ is flat, so is $U\to f(U)$.
    We may then apply Theorem 5.6 in loc.\ cit.\ to conclude that $f$
    is equinormalizable at $x$ (in the sense of the paragraph above Corollary 5.4.2 in loc.\ cit.). 
    Since $x\in X$ was arbitrary, it follows from Corollary 5.4.2(ii) in loc.\ cit.\ that $f$ is equinormalizable. 
\end{proof}
Let $p\colon\mathcal X\to H$ be a family of one-nodal hyperelliptic curves of arithmetic genus $g \geq 2$ over a connected normal complex analytic space $H$. 
Assume that:
\begin{enumerate}
    \item For some $k\geq 1$ and some family of semi-stable curves $q:\mathcal C\to H$, 
there is an isogeny
$$
\psi \colon (J\mathcal X)^k \coloneqq 
J\mathcal X\times_H\dots \times_HJ\mathcal X\longrightarrow J\mathcal C
$$
of semi-abelian varieties over $H$, where $J\mathcal X$ and $J\mathcal C$ are the relative Jacobians over $H$.
\item 
If $\nu\colon \widetilde {\mathcal X}\to \mathcal X$ is a simultaneous normalization (cf.\ Proposition \ref{prop:simultaneous-normalization}), then
\begin{enumerate}[label=(\roman*)]
    \item there is an isomorphism
\begin{align} \label{item:prop:moving:new}
\widetilde{\mathcal X}\cong  \widetilde X_0\times H
\end{align}
over $H$, where $\widetilde X_0\subset \widetilde{\mathcal X}$ is the fibre above a general base point $0\in H$.
\item \label{item:prop:moving:3} there is a non-constant morphism $H\to \widetilde X_0$, $u\mapsto x_u$, such that for any $u\in H$, the fibre $X_u$ is obtained from its normalization $\wt X_u$ by gluing the images of the points $x_u, \iota(x_u) \in \wt X_0$ under the isomorphism $\wt X_0 \cong \wt X_u$ induced by \eqref{item:prop:moving:new}. 
Here, $\iota \colon \wt X_0 \to \wt X_0$ denotes the hyperelliptic involution on $\widetilde X_0$.
\end{enumerate}
\end{enumerate}  

\begin{proposition} \label{prop:moving}
In the above notation, let $\wt C_{0,1}, \dotsc, \wt C_{0,n}$ be the non-rational connected components of the normalization $\wt C_0$ of the curve $C_0 = q^{-1}(0)$. 
Then there is an integer $N\geq 1$ and a matrix $\alpha = (a_{i,j})_{ij} \in \rm{M}_k(\Z)$ with non-zero determinant, such that for each $j \in \set{1, \dotsc, k}$,  
the image of the composition 
$$ 
g_j\colon \wt X_0\stackrel{f_j}\longrightarrow (J\wt X_0)^k\stackrel{\alpha} \longrightarrow (J\wt X_0)^k\stackrel{\wt \psi_0} \longrightarrow J\wt C_0=J \wt C_{0,1} \times \cdots \times J \wt C_{0,n}
$$ 
is contained in 
$$
N \cdot \left(\wt C_{0,1} - \wt C_{0,1} \right) \times \cdots \times N \cdot \left(\wt C_{0,n} - \wt C_{0,n} \right)\subset J \wt C_{0,1} \times \cdots \times J \wt C_{0,n},
$$
where $f_j(x) = (0, \dotsc, x - \iota(x), \dotsc, 0)$ with $x - \iota(x)$ placed on the $j$-th coordinate and where $\wt \psi_0$ is the isogeny induced by $\psi$. 
Moreover, for each $i \in \set{1, \dotsc, n}$, there exists $j \in \set{1, \dotsc, k}$ such that the resulting morphism $$\pr_i \circ g_j \colon \wt X_0 \longrightarrow N \cdot (\wt C_{0,i} - \wt C_{0,i})$$ is non-constant. 
\end{proposition}

\begin{proof}
The image of $g_j$ is analytic and the image of $f_j$ is one-dimensional.
Hence, in order to prove the proposition, we are allowed to perform a base change along an arbitrary morphism $\tau\colon H'\to H$ of complex analytic spaces as long as $0\in \im(\tau)$ and the point $x_u\in \widetilde X_0$  from item \ref{item:prop:moving:3} moves if $u$ runs along $H\cap \im(\tau)$. 
This easily reduces us to the situation where $H$ is a  one-dimensional  disc and the number of irreducible components of the curve $C_u = q^{-1}(u)$ is constant for $u \in H$. 
By the existence of the isogeny $\psi \colon (J\ca X)^k \to J\ca C$, this implies that the number of nodes of $C_u$ is constant for $u \in H$. 
Indeed, if $\nu$ is the number of nodes of $C_u$ and $c$ the number of irreducible components of $C_u$, then one has the formula $k = \nu - c + 1$.  

Since the fibres of $p\colon \mathcal X\to H$ are one-nodal and $\psi\colon (J\mathcal X)^k\to J\mathcal C$ is an isogeny, we have that $k = \text{rank}_\Z(H_1(\Gamma(C_0),\Z))$, where $\Gamma(C_0)$ denotes the dual graph on $C_0$.
By item \eqref{item:lemma:basis:homology:suitable} in Lemma 
\ref{lemma:basis:homology:suitable}, there is a linearly independent subset 
$
\set{\gamma_1(0), \dotsc, \gamma_k(0)} \subset H_1(\Gamma(C_0),\Z)
$
of homology classes of the form \eqref{equation:homologycycle}. 
Since $\mathcal C \to H$ is equisingular, the groups $H_1(\Gamma(C_u),\Z)$ form for $u\in H$ a local system on $H$.
Since $H$ is a disc, the corresponding local system is trivial and so each $\gamma_j(0)$ extends to a section of classes 
$
\gamma_j(u)\in H_1(\Gamma(C_u),\Z)$.
Hence, for each $u\in H$ we get a linearly independent subset 
\begin{equation}\label{eq:gamma_j(u)}
\set{\gamma_1(u), \dotsc, \gamma_k(u)} \subset  H_1(\Gamma(C_u),\Z)
\end{equation} 
of homology classes of the form \eqref{equation:homologycycle}. 

\begin{claim} \label{claim:global-extension}
There exists a family of abelian varieties 
$
\ca A \to H
$
with fibre $\ca A_u \cong J \wt C_u$ above a point $u \in H$, so that the semi-abelian scheme $J \ca C$ is globally over $H$ the extension of $\ca A$ by the torus $\bb G_m^{ k} \times H$ over $H$, and the isogeny $\psi \colon (J\ca X)^k \to J\ca C$ induces an isogeny 
\[
\wt \psi \colon (J \wt X)^k \longrightarrow \ca A
\]
of families of abelian varieties over $H$. 
\end{claim}

\begin{proof}[Proof of Claim \ref{claim:global-extension}]
The claim follows easily once the existence of $\ca A \to H$ is established. 
In the algebraic setting, this is \cite[Chapter I, Corollary 2.11]{faltings-chai}. As we are working in the analytic setting, we provide an argument for convenience of the reader. 

Let $(\bb V, W, \ca F)$ be the polarized integral variation of mixed Hodge structure on $H$ (cf.\ \cite[Definitions 14.44 \& 14.45]{petersteenbrink}) defined by the family of curves $q \colon \ca C \to H$. In particular, the underlying local system of $\Z$-modules $\bb V$ has stalk $\bb V_{u} = H^1(C_u,\Z)$ for a point $u \in H$. The quotient $\bb V/W_0(\bb V)$ is a local system on $H$ with stalk $$H^1(C_u,\Z)/W_0H^1(C_u,\Z) = H^1(\wt C_u,\Z)$$ for $u \in H$, and the filtration $\ca F$ induces a filtration $\overline{\ca F}$ on the holomorphic vector bundle $\left(\bb V/W_0(\bb V) \right) \otimes_{\Z} \OO_H$ that extends to a principally polarized integral variation of Hodge structure of weight one over $H$. 
This concludes the proof of the claim.
\end{proof}

Since $(J \wt X)^k \cong (J \wt X_0)^k \times H$ as families of abelian varieties over $H$, in view of the isogeny $\wt \psi$ above, the polarized abelian scheme $\ca A \to H$ from the above claim is isotrivial, and hence constant, since $H$ is simply connected. Thus, we get a canonical isomorphism
$$
\ca A \xlongrightarrow{\sim} J \wt C_0 \times H
$$
of principally polarized abelian schemes over $H$, yielding a canonical isomorphism of principally polarized abelian varieties
\begin{align} \label{align:identification-of-jacobians-over-H}
J \wt C_u = \ca A_u \xlongrightarrow{\sim} J \wt C_0
\end{align}
for each $u \in H$. 

Consider the homomorphism 
\[
c_t \colon H_1(\Gamma(C_u),\Z) \longrightarrow J \wt C_u \cong J \wt C_0 = \prod_{i = 1}^n J \wt C_{0,i},
\]
see Section \ref{subsection:extensionclassesofnodalcurves} and in particular equations \eqref{equation:nodalidentification} and \eqref{equation:ct}. 
Here, the isomorphism $J \wt C_u \cong J \wt C_0$ is the one defined in \eqref{align:identification-of-jacobians-over-H} above.
By item \eqref{item:lemma:alexeev1motif} in Lemma \ref{lemma:basis:homology:suitable}, for each $u\in H$ there are points $p_{i,j}(u), q_{i,j}(u) \in \wt C_{0,i}$ for $i \in \set{1, \dotsc, n}$ and $j \in \set{1, \dotsc, k}$, such that
\begin{align} \label{equation:extensionclasC0}
c^t(\gamma_j(u)) = \left(p_{1,j}(u) - q_{1,j}(u), \dotsc, p_{n,j}(u) - q_{n,j}(u)\right) \in & J \wt C_{0,1} \times \cdots \times J \wt C_{0,n}.
\end{align}  
Up to a suitable base change we may assume that the points $p_{i,j}(u),q_{i,j}(u)$ depend holomorphically on $u$ and give rise to sections of $\mathcal C\to H$. 

By Lemma \ref{lemma:extension-class-commutativity} and equation \eqref{equation:nodalidentification}, for each $u \in H$, the isogeny of semi-abelian varieties $\psi_u \colon (JX_u)^k\to JC_u $ induces a canonical morphism $(\psi_u)_\ast \colon H_1(\Gamma(X_u),\Z)^{\oplus k}  \to H_1(\Gamma(C_u),\Z) $ such that the following diagram commutes:
\begin{align}
\label{diagram:commutative-extension-diagram}
\xymatrixcolsep{4pc}
\xymatrixrowsep{3pc}
\begin{split}
\xymatrix{
H_1(\Gamma(X_u),\Z)^{\oplus k} \ar[r]^{(\psi_u)_\ast} \ar[d]^{c^t_{(JX_u)^k}} & H_1(\Gamma(C_u),\Z) \ar[d]^{c^t_{JC_u}} \\
(J \wt X_0)^k \ar[r]^{\wt \psi_0} & J \wt C_0.
}
\end{split}
\end{align}
By assumptions, for $u \in H$, we have $\wt X_0 \cong \wt X_u$ and the curve $X_u$ is obtained from its normalization $\wt X_u$ by gluing the images in $\wt X_u$ of the points $x_u$ and $\iota(x_u)$ on $\wt X_0$.
 In particular, the dual graph $\Gamma(X_u)$ consists of a single loop.
 We fix an orientation of this dual graph (see Definition \ref{definition:orientation}) and obtain a canonical identification $H_1(\Gamma(X_u),\Z)=\Z$.
 Let
$$
\mu(u)_r = (0, \dotsc,0, 1,0, \dotsc, 0) \in H_1(\Gamma(X_u),\Z)^{\oplus k},
$$
where $1$ is placed on the $r$-th coordinate. 
Then
\[
c^t_{(JX_u)^k}
\left(
\mu(u)_r
\right) 
= 
\left(
0, \dotsc, x_u - \iota(x_u), \dotsc, 0
\right)  \in (J \wt X_u )^k = (J \wt X_0)^k,
\] 
where $x_u - \iota(x_u)$ is placed on the $r$-th coordinate. 
Recall the set of classes $\gamma_j(u)\in H_1(\Gamma(C_u),\Z)$ from \eqref{eq:gamma_j(u)}.
As the elements $\mu(u)_1, \dotsc, \mu(u)_k$ form a basis of $H_1(\Gamma(X_u),\Z)^{\oplus k}$, and the cokernel of the embedding
$$
(\psi_u)_\ast \colon 
H_1(\Gamma(X_u),\Z)^{\oplus k} \longrightarrow H_1(\Gamma(C_u),\Z)
$$
is finite (because $\psi_u$ is an isogeny), there exists an integer $N \in \Z_{\geq 1}$ and integers
$a_{1,j}, \dotsc, a_{k,j} \in \Z$ for each $j \in \set{1, \dotsc, k}$, such that
\begin{align} \label{equation:composition-pointonjac}
N \cdot \gamma_{j}(u) = \sum_{r = 1}^k a_{r,j} \cdot (\psi_u)_\ast (\mu(u)_r) \in H_1(\Gamma(C_u),\Z)
\end{align}
for each $j \in \set{1, \dotsc, k}$. 
By the commutativity of the diagram \eqref{diagram:commutative-extension-diagram}, together with \eqref{equation:extensionclasC0}, we obtain the following equalities for each $u \in H$ and each $j \in \set{1, \dotsc, k}$: 
\begin{align*}
\wt \psi_0
(
a_{1,j} \cdot &(x_u-\iota(x_u)), \dotsc, 
a_{k,j} \cdot (x_u-\iota(x_u))
)  \\
&
= 
\sum_{r = 1}^k a_{r,j} \cdot 
\wt \psi_0 
\left(
0, \dotsc, x_u - \iota(x_u), \dotsc, 0
\right)= 
\sum_{r = 1}^k a_{r,j} \cdot 
\wt \psi_0 \left(
c^t_{(JX_u)^k}(\mu(u)_r)
\right) \\
&= 
\sum_{r = 1}^k a_{r,j} \cdot c^t_{JC_u}\left(
(\psi_u)_\ast(\mu(u)_r)
\right) 
=
c^t_{JC_u}
\left(
\sum_{r = 1}^k a_{r,j} \cdot (\psi_u)_\ast (\mu(u)_r)
\right)  \\
& = 
c^t_{JC_u}\left(
N \cdot \gamma_{j} 
\right) 
=
N \cdot 
c^t_{JC_u}(\gamma_{j}) \\
&
=
N \cdot  \left(p_{1,j}(u)- q_{1,j}(u), \dotsc, p_{n,j}(u) - q_{n,j}(u)\right)
\quad \in \;\; J \wt C_{0,1} \times \cdots \times J \wt C_{0,n}.
\end{align*}

We may now consider the integral $k\times k$ matrix $\alpha \coloneqq (a_{i,j})_{i,j}$ and we note that this matrix has non-zero determinant because 
$\gamma_1(u),\dots ,\gamma_k(u) \in H_1(\Gamma(C_u),\Z)$ from \eqref{eq:gamma_j(u)} form a rational basis.
Then 
\[
\alpha \cdot f_j(x_u) = 
\left(
a_{1,j} \cdot (x_u-\iota(x_u)), \dotsc, 
a_{k,j} \cdot (x_u-\iota(x_u))
\right),
\] 
and so 
we conclude that
\[
\wt \psi_0 \left(
\alpha \cdot f_j(x_u)
\right)
= 
N \cdot  \left(p_{1,j}(u)- q_{1,j}(u), \dotsc, p_{n,j}(u) - q_{n,j}(u)\right)
\quad \in \;\; J \wt C_{0,1} \times \cdots \times J \wt C_{0,n}.
\]  
As $\wt \psi_0$ is an isogeny, $\wt \psi_0 \left(
\alpha \cdot f_j(x_u)
\right)$ moves with $u\in H$ because the map $H\to \widetilde X_0$, $u\mapsto x_u$ is non-constant by assumption. 
This shows that the restriction of $g_j \colon \widetilde X_0\to J\widetilde C_0$ to some analytically open non-empty subset of $\widetilde X_0$ has the property claimed in the proposition, which suffices to conclude. 

It remains to prove the last assertion. Note that the curves $f_j(\wt X_0)$ for $j = 1, \dotsc, k$ generate $(J \wt X_0)^k$. As $\wt \psi_0$ and $\alpha$ are isogenies, the curves $g_j(\wt X_0)$ for $j = 1, \dotsc, k$ generate $J \wt C_0 = \prod_{i = 1}^n J \wt C_{0,i}$. Thus, for $i \in \set{1, \dotsc, n}$, the projection $\pr_i \colon g_j(\wt X_0) \to J \wt C_{0,i}$, and hence the projection $\pr_i \colon g_j(\wt X_0) \to N \cdot (\wt C_{0,i} - \wt C_{0,i})$, is non-constant for some $j \in \set{1, \dotsc, k}$. 
This concludes the proof of the proposition.
\end{proof}

\section{Hyperelliptic curves on hyperelliptic Jacobians} \label{section:genusgcurveshyperellipticjacobians}
 
The following result is due to Naranjo and Pirola, see \cite[Theorem 1.1]{naranjopirola2018}. 

\begin{theorem}[Naranjo--Pirola] \label{theorem:C-hyperell-g=3}
Let $X$ be a very general hyperelliptic curve of genus $g \geq 3$. Suppose that the Jacobian $JC$ of some hyperelliptic curve $C$ is isogenous to $JX$. Then $C \cong X$.  
\end{theorem}

 \begin{remark} \label{rem:naranjo-pirola:2}  
 The statement of \cite[Theorem 1.1]{naranjopirola2018} is more general than the above Theorem \ref{theorem:C-hyperell-g=3}, but the proof of \cite[Theorem 1.1]{naranjopirola2018} contains a gap, see \cite[line -5 in the proof of Theorem 1.1]{naranjopirola2018}. 
Here, it is claimed that the nodal curves $C_0$ and $D_0$ in loc.\ cit.\ are isomorphic, while the given arguments only suffice to conclude that their normalizations $\tilde C_0$ and $\tilde D_0$ are isomorphic. 
(Note that the authors do indeed prove that $\tilde C_0 \cong \tilde D_0$. 
Indeed, they show that $\tilde D_0$ is hyperelliptic, hence one can apply Theorem \ref{theorem:C-hyperell-g=3} to the isogeny  $\tilde f_0 \colon J\tilde D_0 \to J \tilde C_0$.)
Proving Theorem \ref{theorem:maintheorem1:hyperelliptic} will in particular fix the gap in loc.\ cit.  
Naranjo and Pirola have informed us that it is possible to find an alternative fix via the study of infinitesimal variations of Hodge structures due to Griffiths and Voisin. 
\end{remark}

\begin{remark} \label{remark:naranjo-pirola-advances}
Theorem \ref{theorem:C-hyperell-g=3} is the part of \cite[Theorem 1.1]{naranjopirola2018} that is not affected by the aforementioned gap. 
To explain this, in the notation of \cite[Theorem 1.1]{naranjopirola2018}, assume that for a very general hyperelliptic curve $C$ there exists an isogeny $f \colon JD \to JC$ where $D$ a smooth hyperelliptic curve of genus $g \geq 3$. 
The moduli count on \cite[p.\ 901]{naranjopirola2018} is used to show that the map $\phi \colon \tilde D_0 \to J \tilde D_0$ defined as $y \mapsto m(y - \iota_{\tilde D_0}(y))$ (see page 901 in loc.\ cit.)\ is birational onto its image. 
Although this moduli count is incorrect, $\phi$ is birational onto its image for the following reason: if $2m \tilde D_0 \subset J \tilde D_0$ denotes the image of $\phi$, then the geometric genera of $\tilde D_0$ and $2m \tilde D_0$ are the same because $2m \tilde D_0$ is dominated by $\tilde D_0$ and generates $J \tilde D_0$. 
The isogeny $\tilde f_0^\ast \colon J\tilde C_0\to J\tilde D_0$ in loc.\ cit.\ sends, by comparison of the respective extension classes, the curve $2 \tilde C_0$ into the curve $2m\tilde D_0$.
By the above, this provides us with a dominant rational map $\tilde C_0\dashrightarrow  \tilde D_0$, which must be an isomorphism because both curves have the same genus at least two. 
The remaining arguments in loc.\ cit.\ (together with Proposition \ref{proposition:marcucci:improved} in Appendix \ref{appendix}) suffice to prove that $JD \cong JC$, hence $D \cong C$ by Torelli and the genericity assumptions.   
\end{remark}

The goal of this section is to deduce from Theorem \ref{theorem:C-hyperell-g=3} the following generalization. 

\begin{theorem} \label{proposition:matrixalphai}
Let $X$ be a very general hyperelliptic curve of genus $\geq 3$. Let $C_1, \dotsc, C_n$ be hyperelliptic curves of genus $\geq 1$ such that there exists an isogeny 
$
JC_1 \times \cdots \times JC_n \to (JX)^k
$
for some $k \geq 1$. 
Then $n = k$, there is an isomorphism $C_i \cong X$ for each $i \in \set{1,\dotsc, n}$, and the induced isogeny 
$
(JX)^k \cong J  C_{1} \times \cdots \times J  C_{n} \to (J X)^k
$
is given by a matrix $\alpha \in \rm{M}_k(\Z)$. 
\end{theorem}

The proof of Theorem \ref{proposition:matrixalphai} relies on Theorem \ref{theorem:C-hyperell-g=3} and the following two results, the first of which is probably well-known, and the second of which is due to Lazarsfeld and Martin \cite{lazarsfeld2023measures}. 

\begin{lemma} \label{lem:factor of A}
    Let $A$ be an abelian variety with $\End(A) = \Z$.
    The natural maps $\rm M_k(\Z)\to \End(A^k)$ and $\GL_k(\Z)\to \Aut(A^k)$ are isomorphisms. Moreover, if there exist abelian varieties $B_1,\dots ,B_n$ and an isomorphism $\varphi\colon\prod_{i = 1}^n B_i \stackrel{\sim}\to A^k$, then there exists $\alpha\in \Aut(A^k) = \GL_k(\Z)$ such that the composition $\alpha \circ \varphi\colon\prod_{i = 1}^n B_i\to A^k $ respects the product structures on both sides.
    In particular, in that case, there is a partition $k=k_1+\dots +k_n$ such that $B_i\cong A^{k_i}$ for each $i$.
\end{lemma}
\begin{proof}
The space of endomorphisms $\End(A^k)$ is naturally given by $k\times k$ matrices whose entries are endomorphisms of $A$.
Since $\End(A)=\Z$, we find that $\End(A^k)\cong \rm{M}_k(\Z)$, which proves the first claim in the lemma.
This also implies $\Aut(A^k)\cong \GL_k(\Z)$.

The $i$-th factor $B_i$ yields a projector $p_i\in \End(A ^k)$.
By what we have said above, $p_i$ can be identified with a $k\times k$ matrix with $p_i\cdot p_i=p_i$ and $p_i\circ p_j=p_j\circ p_i$ for all $i,j$.
By simultaneous diagonalization of permuting projectors, we find  a change of coordinates, i.e.\ an automorphism $\alpha\in \Aut(A^k)=\GL_k(\Z)$, such that $\alpha\circ \varphi$ has the property claimed in the lemma. 
\end{proof}

\begin{proposition}[Lazarsfeld--Martin] \label{proposition:lazarsfeldmartin}
Let $X$ be a very general hyperelliptic curve of genus $ g \geq 3$ and let $Z \subset JX \times JX$ be an irreducible curve whose normalization is hyperelliptic. Then $Z$ generates a proper subtorus of $JX \times JX$. 
\end{proposition}
\begin{proof}
See \cite[Proposition 3.1]{lazarsfeld2023measures}. In the statement of that proposition, there is the additional assumption that $Z$ lifts to a curve on $X \times X$, but this assumption is not used in the proof. 
Indeed, the proof in loc.\ cit.\ immediately starts with the hyperelliptic curve $Z$ inside $JX\times JX$ and spreads this out to a family of hyperelliptic curves $Z_s\subset JX_s\times JX_s$.
Then the hyperelliptic Jacobian $JX_s$ is specialized to  $JX_s=B_s\times E$, where $E$ is a fixed elliptic curve and $B_s$ is the Jacobian of a very general hyperelliptic curve of genus $g-1$ which varies with $s$.
Under the assumption that $Z$ generates $JX\times JX$, the same will be true for $Z_s$ for all $s$.
In particular, under this assumption the image of $Z_s$ via the projection $JX_s\times JX_s\to E\times E$ is a curve in $E\times E$.
It is then shown (see \cite[Claim in Section 3]{lazarsfeld2023measures}) that this curve varies with $s$.
Since the normalization of $Z_s$ is hyperelliptic, its image in $E\times E$ yields a rational curve on the Kummer surface associated to $E\times E$.
The latter are rigid, because Kummer surfaces are not ruled.
This contradiction concludes the argument.
\end{proof}

\begin{proof}[Proof of Theorem \ref{proposition:matrixalphai}] 
As $J X$ is simple, we have $g(C_i) \geq g(X)$ for each $i$, and $n \leq k$. If $g(C_i) > g(X)$ for some $i$, then $k \geq 2$ and there exists a surjection
\[
JC_i \twoheadrightarrow JX \times J X. 
\] 
As $C_i$ is hyperelliptic, and $X$ is very general hyperelliptic, this contradicts Proposition \ref{proposition:lazarsfeldmartin}. We conclude that $g(C_i) = g(X)$ for each $i$, and that $n = k$. In particular, for each $i$, there exists an isogeny
$
J C_i \to J X.
$
By Theorem \ref{theorem:C-hyperell-g=3}, we have $C_i \cong X$ for each $i$. As the hyperelliptic curve $X$ is very general, the composition $(J X)^k \cong \prod_{i = 1}^n J C_i \to (J X)^k$ is given by a matrix in $ \rm{M}_k(\Z)$, see Lemma \ref{lem:factor of A}. 
\end{proof} 

\section{Polarizations on powers of abelian varieties and bilinear forms} \label{subsection:polarizations-on-powers}

This section has two goals. 
Consider a principally polarized abelian variety $A$ with endomorphism ring $\Z$. 
Firstly, we classify isomorphism classes of principal polarizations on any power of $A$, see Section \ref{subsection:polarizations-on-powers}. Secondly, we investigate principal polarizations on any abelian variety $B$ isogenous to a power of $A$, see Section \ref{appendix:polarizations}.  

 \subsection{Polarizations on powers of a very general abelian variety} 
 \label{subsubsection:polarizations-on-powers}

Let $A$ be an abelian variety with dual abelian variety $A^\vee$. For a line bundle $\ca L$ on $A$, the map $\varphi_{\ca L} \colon A \to \Pic^0(A) = A^\vee$ defined as $x \mapsto t_x^\ast(\ca L) \otimes \ca L^{-1}$ is a homomorphism of abelian varieties, and the association $\ca L \mapsto \varphi_{\ca L}$ induces an injective map 
$$
\Phi\colon \NS(A) \longhookrightarrow \Hom(A, A^\vee).
$$  
The image of $\Phi$ is contained in the subset $\Hom^{sym}(A, A^\vee) \subset \Hom(A, A^\vee)$ of maps $\phi \colon A \to A^\vee$ that satisfy $\phi^\vee = \phi $ (viewed as maps $ A^{\vee\vee} = A \to A^\vee$). 
A line bundle $\ca L$ on $A$ is ample if and only if $\varphi_{\ca L} \colon A \to A^\vee$ is an isogeny, in which case the class $[\ca L] \in \rm{NS}(A)$ (resp.\ the homomorphism $\varphi_{\ca L}$) is called a polarization.
The polarization $[\ca L]$ is principle if $h^0(A, \ca L) = 1$, or equivalently, if  $\varphi_{\ca L}$ is an isomorphism. See e.g.\ \cite{Milne1986} or \cite{birkenhakelange-complexAVs} for more details. 

\begin{lemma} \label{lemma:complextorusduality}
    Let $A = V/\Lambda$ be a complex torus. Let $k\in \Z_{\geq 1}$ and consider the natural embedding
    \[
    \iota_A \colon \rm{M}_k(\Z) \longhookrightarrow \End(A^k). 
    \]
    Let $\alpha \in \rm{M}_k(\Z)$ with attached endomorphism $\iota_A(\alpha) \in \End(A^k)$. 
    Let $\iota_A(\alpha)^\vee \in \End((A^k)^\vee)$ be the endomorphism dual to $\iota_A(\alpha)$.
    Then, with respect to the canonical isomorphism $(A^k)^\vee = (A^\vee)^k$, we have $\iota_A(\alpha)^\vee = \iota_{A^\vee}(\alpha^t)$ where $\alpha^t \in \rm{M}_k(\Z)$ is the transpose of the matrix $\alpha$.  
\end{lemma}

\begin{proof}
The lemma follows from the following well-known linear algebra statement: if $\Lambda$ is a free $\Z$-module of positive and finite rank, and if $\alpha  \in \rm{M}_k(\Z)$, then the endomorphism $f_\Lambda(\alpha) \colon \Lambda^k \to \Lambda^k$ that $\alpha$ induces satisfies
$
f_\Lambda(\alpha)^\vee = f_{\Lambda^\vee}(\alpha^t)$ as morphisms $(\Lambda^\vee)^k \longrightarrow (\Lambda^\vee)^k$. Here, $f_{\Lambda^\vee}(\alpha^t)$ is the endomorphism of $(\Lambda^\vee)^k$ attached to the transpose $\alpha^t \in \rm{M}_k(\Z)$ of $\alpha$. 
\end{proof}

Let $k \in \Z_{\geq 1}$ and let $A$ be an abelian variety, 
principally polarized by $\lambda \colon A \xrightarrow{\sim} A^\vee$. Let $\lambda^k \colon A^k \to (A^k)^\vee$ be the product polarization on $A^k$. 
From now on we drop the notation $\iota_A$ introduced in Lemma \ref{lemma:complextorusduality} by letting $\alpha \in \End(A^k)$ denote the endomorphism attached to a matrix $\alpha \in \rm{M}_k(\Z)$.
We define an injective map
\begin{align} \label{align:lambda-alpha}
\begin{split}
\rm{M}_k(\Z) &\longhookrightarrow \Hom(A^k, (A^k)^\vee),\\
\alpha &\mapsto \lambda_\alpha \coloneqq \lambda^k \circ \alpha. 
\end{split}
\end{align} 
Observe that any $\gamma\in \GL_k(A) $ acts naturally on the set of morphisms $\mu \colon A^k \to (A^k)^\vee$ via 
\begin{align} \label{align:action-aut-polarizations}
\mu\mapsto \gamma \circ \mu \circ \gamma^t,
\end{align}
where we view $\gamma$ as an automorphism of $A^k$, $\mu$ as a morphism $A^k \to (A^\vee)^k$, and the transpose $\gamma^t \in \GL_k(\Z)$ of $\gamma$ as an automorphism of $(A^\vee)^k$. If $\mu$ is a polarization, then $\gamma \circ \mu \circ \gamma^t$ is again a polarization. Indeed, it is clear that $(\gamma^{t})^\vee \circ \mu \circ \gamma^t$ is a polarization, where the automorphism $(\gamma^t)^\vee \colon (A^\vee)^k\xrightarrow{\sim} (A^\vee)^k$ is the automorphism of $(A^\vee)^k$ induced by $\gamma^t$ via duality; moreover, by Lemma \ref{lemma:complextorusduality} we have $(\gamma^t)^\vee = (\gamma^t)^t = \gamma$ as automorphisms $(A^\vee)^k\xrightarrow{\sim} (A^\vee)^k$. 

If $\End(A) =\Z$, then $\End(A^k) = \rm{M}_k(\Z)$ and $\Aut(A^k) = \GL_k(\Z)$ by Lemma \ref{lem:factor of A}, and the map $\alpha \mapsto \lambda_\alpha$ yields an isomorphism $\rm{M}_k(\Z) \xrightarrow{\sim} \Hom(A^k, (A^k)^\vee)$.  

\begin{lemma} \label{lemma:polarizationcorrespondence}
Let $(A, \lambda)$ be a principally polarized abelian variety. Then the following holds.
\begin{enumerate}
\item \label{item:polarizationcorrespondence-0} Let $\alpha \in \rm{M}_k(\Z)$ such that $\alpha$ has non-zero determinant. Then the map $\lambda_\alpha \colon A^k \to (A^k)^\vee$ associated to $\alpha$ is a polarization on $A^k$ if and only if $\alpha$ is symmetric and positive definite. In particular, the map \eqref{align:lambda-alpha} restricts to an injective map
\begin{align}
\label{equation:bijection}
\begin{split}
\rm{P}_k(\Z) &\longhookrightarrow
\set{\tn{\emph{polarizations on} $A^k$}}, \\ \alpha &\mapsto \lambda_\alpha,
\end{split}
\end{align}
where $\rm{P}_k(\Z) \subset \rm{M}_k(\Z)$ denotes the subset of positive definite symmetric matrices. 
\item \label{item:polarizationcorrespondence-2}
The map \eqref{align:lambda-alpha} is equivariant with respect to the $\GL_k(\Z)$-action on both sides, where $\GL_k(\Z)$ acts on $\rm{P}_k(\Z)$ by $\gamma \cdot \alpha = \gamma \alpha \gamma^t$ for $\alpha \in \rm{P}_k(\Z)$ and $\gamma \in \GL_k(\Z)$, and where $\GL_k(\Z)$ acts on $A^k$ via \eqref{align:action-aut-polarizations} and the natural embedding $\GL_k(\Z) \subset \Aut(A^k)$. 
\item \label{item:polarizationcorrespondence-1} If $\End(A) = \Z$, then \eqref{equation:bijection} defines a bijection \begin{align}
\label{equation:bijection:ii}
\begin{split}
\rm{P}_k(\Z) &\xlongrightarrow{\sim}
\set{\tn{\emph{polarizations on} $A^k$}}, \\ \alpha &\mapsto \lambda_\alpha.
\end{split}
\end{align} 
\item \label{item:polarizationcorrespondence-4}
Assume $\End(A) = \Z$. Let $\alpha, \beta \in \rm{P}_k(\Z)$ and $\gamma \in \GL_k(\Z) = \Aut(A^k)$. Then $\alpha = \gamma \beta \gamma^t$ if and only if $\gamma^t$ defines an isomorphism of polarized abelian varieties $(A^k, \lambda_\alpha) \xrightarrow{\sim}(A^k, \lambda_\beta)$.  
\item \label{item:polarizationcorrespondence-3} If $\End(A) = \Z$, then the map \eqref{align:lambda-alpha} induces a bijection between the set of (isomorphism classes of) unimodular positive definite symmetric bilinear forms on $\Z^k$ and the set of (isomorphism classes of) principal polarizations on $A^k$. 
\end{enumerate}
\end{lemma}

\begin{proof}
Let us prove item \eqref{item:polarizationcorrespondence-0}. Let $\alpha \in \rm{M}_k(\Z)$ such that $\alpha$ has non-zero determinant. We need to show that $\lambda_\alpha \colon A^k \to (A^k)^\vee$ is a polarization on $A^k$ if and only if $\alpha$ is symmetric and positive definite. 
For this, in view of \cite[Theorem 5.2.4]{birkenhakelange-complexAVs}, it suffices to prove that $\alpha$ is symmetric and positive definite if and only if the endomorphism $\alpha \in \End(A^k)$ is symmetric and totally positive. Here, \emph{symmetric} is understood to be with respect to the Rosati involution $\dagger \colon \End(A^k) \to \End(A^k)$ defined by the product principal polarization $\lambda^k$ on  $A^k$ attached to $\lambda$ (thus, $f \in \End(A^k)$ is symmetric if $f^\dagger = f$) and an endomorphism $\varphi \colon X \to X$ of an abelian variety $X$ is said to be \emph{totally positive} if the zeros of the characteristic polynomial of the analytic representation $\varphi^{an} \colon \tn{Lie}(X) \to \tn{Lie}(X)$ of $\varphi$ are all positive. 

Consider the canonical embedding 
\begin{align*}
\rm{M}_k(\Z) 
\longhookrightarrow
\End(A^k).
\end{align*} 
The above shows that to prove item \eqref{item:polarizationcorrespondence-0}, it suffices to prove the following assertions. 
\begin{enumerate}[label=(\alph*)]
    \item \label{item:symmetric} A matrix $\alpha \in \rm{M}_k(\Z)$ is symmetric if and only if the associated endomorphism $\alpha \in \End(A^k)$ is symmetric. 
    \item\label{item:totallypositive} 
    Let $\alpha \in \rm{M}_k(\Z)$ be a symmetric matrix with non-zero determinant. Then the induced $\R$-linear transformation $\alpha\colon \R^k \xrightarrow{\sim} \R^k$ has positive real eigenvalues if and only if the associated endomorphism  $\alpha \in \End(A^k)$ is totally positive. 
\end{enumerate}
Item \ref{item:symmetric} follows from Lemma \ref{lemma:complextorusduality}. 
Let us prove item \ref{item:totallypositive}.  
 Let $\alpha^{an} \colon \rm{Lie}(A)^k \to \rm{Lie}(A)^k$ denote the analytic representation of the endomorphism $\alpha \in \End(A^k)$.  
We must show that 
the eigenvalues of the induced $\R$-linear transformation $\alpha\colon \R^k \xrightarrow{\sim} \R^k$ are positive if and only if the eigenvalues of the complex linear map $\alpha^{an} \colon \rm{Lie}(A)^k \to \rm{Lie}(A)^k$ are positive. This follows readily from the fact that  
if $V$ is a complex vector space of finite positive dimension and $k \geq 1$ an integer, then any symmetric matrix $\alpha \in \GL_k(\R)$ has positive real eigenvalues when viewed as an $\R$-linear transformation $\alpha \colon \R^k \xrightarrow{\sim} \R^k$ if and only if the induced $\C$-linear map $V^k = \R^k \otimes_\R V \xrightarrow{\sim} \R^k \otimes_\R V= V^k$ has positive real eigenvalues.

Next, we prove item \eqref{item:polarizationcorrespondence-2}. 
Let $\alpha \in \rm{P}_k(\Z)$ with associated polarization $\lambda_\alpha \colon A^k \to (A^k)^\vee$. 
Let $(a_1, \dotsc, a_k) \in A^k$. 
Then $\lambda^k \circ \gamma = \gamma \circ \lambda^k$ as maps $A^k \to (A^\vee)^k$ since $\lambda^k$ is the image of the identity matrix under the map \eqref{align:lambda-alpha}. 
Therefore,
\[
\lambda_{\gamma \alpha \gamma^t} = \lambda^k \circ \gamma \circ \alpha \circ \gamma^t = \gamma \circ \lambda^k \circ \alpha \circ \gamma^t = \gamma \circ \lambda_\alpha \circ \gamma^t, 
\] 
hence item \eqref{item:polarizationcorrespondence-2} follows. 

We now prove item \eqref{item:polarizationcorrespondence-1}. Thus, we assume $\End(A) = \Z$. Note that this assumption implies that the map \eqref{align:lambda-alpha} is a bijection. Hence any polarization $\mu \colon A^k \to (A^k)^\vee$ is of the form $\mu = \lambda_\alpha$ for a unique $\alpha \in \GL_k(\Z)$, and by item \eqref{item:polarizationcorrespondence-0} such a matrix $\alpha$ is symmetric and positive definite, proving what we want. 

To prove item \eqref{item:polarizationcorrespondence-4}, assume $\End(A) = \Z$. Let $\alpha, \beta \in \rm{P}_k(\Z)$ and $\gamma \in \GL_k(\Z) = \Aut(A^k)$.
Then $\gamma^t$ defines an isomorphism of polarized abelian varieties $(A^k, \lambda_\alpha) \cong (A^k, \lambda_\beta)$ if and only if $\lambda_\alpha = \gamma \circ \lambda_\beta \circ \gamma^t$. Since $\gamma \circ \lambda_\beta \circ \gamma^t = \lambda_{\gamma \beta \gamma^t}$ by item \eqref{item:polarizationcorrespondence-2}, this happens if and only if $\lambda_\alpha = \lambda_{\gamma \beta \gamma^t}$, which in turn happens if and only if $\alpha = \gamma \beta \gamma^t$ by the injectivity of \eqref{equation:bijection}. Item \eqref{item:polarizationcorrespondence-4} follows. 

Finally, for $\alpha \in \rm{P}_k(\Z)$, the polarization $\lambda_\alpha$ is a principal if and only if $\alpha$ is unimodular. 
Item \eqref{item:polarizationcorrespondence-3} of the lemma follows then from item \eqref{item:polarizationcorrespondence-4}, and we are done. 
\end{proof}

\begin{lemma}
\label{lemma:polarizationcorrespondence-4} 
Let $(A, \lambda)$ a principally polarized abelian variety with $\End(A) = \Z$. Let $k$ be a positive integer, let $\alpha \in \rm{P}_k(\Z)$, and consider the associated polarization $\lambda_\alpha$ on $A^k$. Consider the positive definite integral quadratic space $(\Z^k, \alpha)$ associated to the matrix $\alpha$. 
 \begin{enumerate} 
 \item \label{item:correspondencedecompositions} The association $\alpha \mapsto \lambda_\alpha$ in \eqref{equation:bijection:ii} induces a bijection between the set of isomorphism classes of decompositions $(A^k, \lambda_\alpha) \cong \prod_i (B_i, \lambda_i)$ for some polarized abelian varieties $(B_i, \lambda_i)$ and the set of isomorphism classes of decompositions $(\Z^k, \alpha) \cong \oplus_i (\Z^{k_i}, \alpha_i)$ of $(\Z^k, \alpha)$ into an orthogonal direct sum of positive definite integral quadratic spaces.  
 \item \label{item:polarizedabeliansubvariety}
 The polarized abelian variety $(A^k, \lambda_\alpha)$ is indecomposable as a polarized abelian variety if and only if $(\Z^k, \alpha)$ is an indecomposable positive definite integral quadratic space.  
 \end{enumerate}
\end{lemma}

\begin{proof}
If $(\Z^k, \alpha) \cong \oplus_i (\Z^{k_i}, \alpha_i)$ as integral quadratic spaces, then $(A^k, \lambda_\alpha) \cong \prod_i (A^{k_i}, \lambda_{\alpha_i})$ by Lemma \ref{lemma:polarizationcorrespondence}. Conversely, consider an isomorphism $ \prod_i (B_i, \lambda_i)  \stackrel{\sim}\to (A^k, \lambda_\alpha) $ of polarized abelian varieties. By Lemma \ref{lem:factor of A}, for each $i$ there exists a non-negative integer $k_i \leq k$ and an isomorphism of abelian varieties $B_i \cong A^{k_i}$. Thus, there is a polarization $\lambda_i'$ on the abelian variety $A^{k_i}$ for each $i$ such that $\prod_i (B_i, \lambda_i) \cong \prod_i (A^{k_i}, \lambda_i')$. Each $\lambda_i'$ is again of the form $\lambda_i' = \lambda_{\alpha_i}$ for $\alpha_i \in \rm{P}_{k_i}(\Z)$, and the resulting isomorphism of polarized abelian varieties $(A^k, \lambda_\alpha) \cong \prod_i (A^{k_i}, \lambda_{\alpha_i})$ is induced by an isomorphism of positive definite integral quadratic spaces $(\Z^k, \alpha) \cong \oplus_i (\Z^{k_i}, \alpha_i)$, see Lemma \ref{lemma:polarizationcorrespondence}. 
Item \eqref{item:correspondencedecompositions} follows, and item \eqref{item:polarizedabeliansubvariety} is a direct consequence of item  \eqref{item:correspondencedecompositions}.
\end{proof}

\subsection{Polarizations on abelian varieties isogenous to a power of an abelian variety} \label{appendix:polarizations}

\begin{lemma} \label{lemma:push-forward-polarization}
    Let $A$ be an abelian variety and let $\lambda \colon A \to A^\vee$ be a polarization. There is a canonical isomorphism $H_1(A^\vee,\Z) = H_1(A,\Z)^\vee$. Moreover, if $E \colon H_1(A,\Z) \times H_1(A,\Z) \to \Z$ is the alternating form corresponding to $\lambda$, then the push-forward
    \[
    \lambda_\ast \colon H_1(A,\Z) \longrightarrow H_1(A^\vee,\Z) = H_1(A,\Z)^\vee
    \]
    satisfies
    $
    \lambda_\ast(x)(y) = E(x,y)$ for all $x,y \in H_1(A,\Z). 
    $
\end{lemma}
\begin{proof}
This is well-known and follows for instance from \cite[Lemma 2.4.5]{birkenhakelange-complexAVs}.
\end{proof}

    Let $(A, \lambda_A)$ and $(B, \lambda_B)$ be principally polarized abelian varieties.
    Let $A^k \to B$ be an isogeny for some integer $k \geq 1$, and suppose that the principal polarization $\lambda_B$ of $B$ pulls back to the polarization $\beta \cdot \lambda_{A^k}$ on $A^k$ defined as 
    \[
   \beta \cdot \lambda_{A^k} \colon A^k \longrightarrow (A^\vee)^k, \quad x \mapsto \beta \cdot \lambda_{A^k}(x),
    \]
    where $\beta$ is a positive definite symmetric integral $k \times k$ matrix and $\lambda_{A^k}$ is the natural product polarization on $A^k$ induced by $\lambda_A$. 
    
\begin{lemma} \label{lemma:betapullback}
In the above notation, let $M\coloneqq H_1(B,\Z)$ and $H\coloneqq H_1(A,\Z)$, and let $E_M$ and $E_{H}^{\oplus k}$ be the symplectic forms on $M$ and $H^{\oplus k}$ associated to the respective principal polarizations. 
    With respect to the natural embedding $M \subset H^{\oplus k}$ induced by the isogeny $A^k \to B$, we have
    \[
    E_M(x,y) = E_{H}^{\oplus k}(x, \beta^{-1} \cdot y) = E_{H}^{\oplus k}(\beta^{-1} \cdot x, y) \quad \quad  \text{ for each } \quad x,y \in M. 
    \]
\end{lemma}
\begin{proof}
Via the principal polarizations, the given isogeny $\phi \colon A^k \to B$ induces an isogeny $\psi \colon B \to A^k$. We claim that $\psi \circ \phi = \beta$ as isogenies $A^k \to A^k$. To see this, recall that the pull-back $\phi^\ast(\lambda_B)$ of the principal polarization $\lambda_B$ of $B$ is the isogeny $A^k \to (A^k)^\vee$ given by the composition
    \begin{align*} 
    \xymatrix{
    A^k \ar[r]^\phi & B \ar[r]^{\lambda_B} &B^\vee \ar[r]^{\phi^\vee} & (A^k)^\vee. 
    }
    \end{align*}
    As this is $\beta$ times the natural principal polarization on $A^k$, the claim follows. 

    On the level of lattices, the maps $\phi$ and $\psi$ induce embeddings
    \[
    \xymatrix{
    H^{\oplus k} \ar@{^{(}->}[r]^f & M \ar@{^{(}->}[r]^g & H^{\oplus k}.
    }
    \]
    The claim above implies that $g \circ f = \beta$ as linear maps $H^{\oplus k} \to H^{\oplus k}$. 
    By assumption, we have $\phi^\ast(\lambda_B) = \beta \cdot \lambda_{A^k}$. 
    Therefore, by Lemma \ref{lemma:push-forward-polarization}, we have
    \[
E_M(f(x), f(y)) =     \left(f^\ast E_M\right)(x,y) = E_{H}^{\oplus k}(\beta \cdot x, y) \quad \quad \quad \forall x,y \in H^{\oplus k}. 
    \]
If we view $M$ as a sublattice of $H^{\oplus k}$ via $g$, then the above equality implies that 
\[
E_M(\beta \cdot x, \beta \cdot y) = E_M(f(x), f(y)) = E_{H}^{\oplus k}(\beta \cdot x,y) \quad \quad \quad \forall x,y \in H^{\oplus k}. 
\]
In particular,
$
E_M(\beta \cdot x, \beta \cdot y) = E_{H}^{\oplus k}(\beta \cdot x,y)$ for all $x,y \in M$, hence
$
E_M(x, y) = E_{H}^{\oplus k}(x,\beta^{-1} \cdot y)$ for $ x,y \in M$. 
As $E_{H}^{\oplus k}(x, \beta^{-1} \cdot y) = E_{H}^{\oplus k}(\beta^{-1} \cdot x, y)$ for each $x,y \in H^{\oplus k}$, the lemma follows.  
\end{proof}

\section{Powers of abelian varieties isomorphic to products of Jacobians} \label{section:polarizations}

The goal of this section is to prove the following theorem, which will be used in the proof of Theorems \ref{theorem:maintheorem1:hyperelliptic} and \ref{thm:cor:maintheorem2}.

\begin{theorem}\label{cor:iso-of-products-of-curves-main-main}
Let $g\in \Z_{\geq 1}$ and let $Z\subset \mathcal A_g$ be a subvariety of the moduli space of principally polarized abelian varieties of dimension $g$ with the following properties:
\begin{itemize}
    \item there is a point $[(A_0,\lambda_0)]\in Z$ such that $A_0\cong E_0\times B_0$ (as polarized abelian varieties), where $B_0$ is a principally polarized abelian variety of dimension $g-1$ and $E_0$ is an elliptic curve with transcendental $j$-invariant; 
    \item a very general point $[(A,\lambda)]\in Z$ satisfies $\End(A)=\Z$.  
\end{itemize} 
If for some very general point  $[(A,\lambda)]\in Z$ and some integers $k,n \geq 1$, there are some smooth projective connected curves $C_1,\dots ,C_n$ of positive genus and an isomorphism $\prod_{i = 1}^nJC_i \cong A^k$ of unpolarized abelian varieties, then $k = n$ and for each $i$ we have an isomorphism $(JC_i,\Theta_{C_i}) \cong (A,\lambda)$ of polarized abelian varieties.
\end{theorem}

\subsection{Applications of Theorem \ref{cor:iso-of-products-of-curves-main-main}} 
Before we turn to the proof of Theorem \ref{cor:iso-of-products-of-curves-main-main}, we show that it implies Theorem \ref{cor:iso-of-products-of-curves-main} stated in the introduction. 

\begin{proof}[Proof of Theorem \ref{cor:iso-of-products-of-curves-main}]
Let $\bar{Z}\subset \mathcal A_g$ be the closure in $\mathcal A_g$ of the image of $Z$ under the Torelli map $\ca M_g \to \ca A_g$. 
By assumption, $\bar{Z}$ contains the hyperelliptic Torelli locus. Since the Jacobian $JX$ of a very general hyperelliptic curve $X$ satisfies $\End(JX)=\Z$, we conclude via specialization that $\End(A)=\Z$ for any very general point $[(A,\lambda)]\in \bar{Z}$.
    Moreover, there are hyperelliptic compact type curves $X_0$ with $JX_0\cong E_0\times B_0$ for an elliptic curve $E_0$ with transcendental $j$-invariant.
    Since $\bar{Z}$ is closed in $\mathcal A_g$ and contains the hyperelliptic locus, the point $[JX_0,\Theta_{X_0}]$ is contained in $\bar{Z}$ and so Theorem \ref{cor:iso-of-products-of-curves-main-main} applies to the subvariety $\bar{Z}$ of $\mathcal A_g$. 
    Hence, for each $i$ we have $(JC_i,\Theta_{C_i}) \cong (JX,\Theta_X)$ as polarized abelian varieties and this implies $C_i\cong X$ by the Torelli theorem.
    This concludes the proof.
\end{proof}

Another consequence of Theorem \ref{cor:iso-of-products-of-curves-main-main} is as follows.

\begin{corollary} \label{cor:cubicthreefold-isomorphism}
Let $Y$ be a very general cubic threefold and $k$ a positive integer. There exist no smooth projective curves $C_1, \dotsc, C_n$ such that $JC_1 \times \cdots \times JC_n \cong (J^3Y)^k$. 
\end{corollary}
\begin{proof} 
By \cite[Theorem (0.1)]{collino-fanofundamentalgroup}, we can degenerate $Y$ into a singular cubic threefold $Y_0$ such that $J^3Y_0 = JX$ is the Jacobian of a very general hyperelliptic curve $X$ of genus five.
In particular, the closure of the locus of intermediate Jacobians of cubics inside $\mathcal A_5$ contains the locus of Jacobians of hyperelliptic curves.
We can then argue as in the proof of Theorem \ref{cor:iso-of-products-of-curves-main} to deduce from Theorem \ref{cor:iso-of-products-of-curves-main-main} that $(J^3Y,\Theta_Y)$ is isomorphic to the Jacobian of a curve, which contradicts the main result of  \cite{clemensgriffiths-cubicthreefold}.
This concludes the proof of the corollary.
\end{proof}

The remaining part of Section \ref{section:polarizations} will be devoted to a proof of Theorem \ref{cor:iso-of-products-of-curves-main-main}.

\subsection{Special subvarieties and powers of abelian varieties isomorphic to Jacobians} \label{sec:specialsubvarieties}

\begin{proof}[Proof of Theorem \ref{cor:iso-of-products-of-curves-main-main}]
Recall that  $Z\subset \mathcal A_g$ is a subvariety such that
\begin{itemize}
    \item there is a point $[(A_0,\lambda_0)]\in Z$ such that $A_0\cong E_0\times B_0$ and some elliptic curve $E_0$ with transcendental $j$-invariant; 
    \item a very general point $[(A,\lambda)]\in Z$ satisfies $\End(A)=\Z$.  
\end{itemize} 
We assume that  for some $k,n \geq 1$, there are some smooth projective connected curves $C_1,\dots ,C_n$ of positive genus and an isomorphism $\prod_{i = 1}^nJC_i \cong A^k$ of unpolarized abelian varieties.
By Lemma \ref{lem:factor of A}, we reduce to the case $n=1$ and get an isomorphism 
$$
f\colon JC\stackrel{\sim}\longrightarrow A^k
$$
for some $k\geq 1$, where $C\coloneqq C_1$.
We aim to prove that $k=1$.
Since $\End(A)=\Z$, this already implies $(JC,\Theta_C)\cong (A,\lambda)$ as polarized varieties, because $A$ carries only one principal polarization, since $\NS(A)=\Z$.

The above isomorphism $f \colon JC \to A^k$ provides $A^k$ with an indecomposable principal polarization, say $\mu$. 
Since $\End(A) = \Z$, one has $\mu = \lambda_\alpha$, the polarization on $A^k$ associated to a positive definite symmetric unimodular bilinear form $\alpha$ on $\Z^k$, see Lemma \ref{lemma:polarizationcorrespondence}.
By abuse of notation, we will denote the principal polarization $\lambda_\alpha$ by $\alpha$ for simplicity, hence write $(A^k,\alpha) \coloneqq (A^k, \lambda_\alpha)$.  

As the principally polarized abelian variety $(A^k,\alpha)$ is indecomposable, the positive definite integral quadratic space $(\Z^k,\alpha)$ is indecomposable, see Lemma \ref{lemma:polarizationcorrespondence-4}. 
Kneser's classification of indecomposable integral inner product spaces of rank at most  $16$ implies then that
$k = 1$ or $k \geq 8$, see  \cite[p.\ 28, Remark 1]{milnorhusemoller} and \cite{kneser-quadratisch}.

As explained above we only need to show $k=1$ and so we assume for a contradiction that $k\geq 8$.
By our assumptions there is a degeneration of $A$ to $A_0\cong E_0\times B_0$, where $E_0$ is an elliptic curve with transcendental $j$-invariant.
This yields an isomorphism $(E_0\times B_0)^k \cong JC_0$ for some compact type degeneration $C_0$ of $C$. If $C_{0,1}, \dotsc, C_{0,n}$ are the non-rational irreducible components of $C_0$, then we can write 
$$
(E_0\times B_0)^k = E_0^k\times B_0^k \cong JC_{0,1} \times \cdots \times JC_{0,n}. 
$$
This is an isomorphism of principally polarized abelian varieties, where the polarization on $E_0^k$ (resp.\ $B_0^k$) is the one induced by $\alpha$ and the principal polarization of $E_0$ (resp.\ $B_0$), see Lemma \ref{lemma:polarizationcorrespondence}
(this step uses that $A_0\cong E_0\times B_0$ as principally polarized abelian varieties.)  
By Lemma \ref{lemma:polarizationcorrespondence-4}, the principally polarized abelian variety $(E_0^k,\alpha)$ is an indecomposable principally polarized abelian variety because $(\Z^k,\alpha)$ is an indecomposable integral inner product space. 
By uniqueness of the decomposition of any principally polarized abelian variety into a product of indecomposable principally polarized abelian subvarieties \cite{clemensgriffiths-cubicthreefold, debarre-produits}, it follows that for some $i$, there is an isomorphism of principally polarized abelian varieties $(JC_{0,i},\Theta_{C_{0,i}}) \cong (E_0^k, \alpha)$. 

To simplify notation we write $E\coloneqq E_0$ and $C\coloneqq C_{0,i}$ and get an isomorphism $(E^k,\alpha) \cong (JC,\Theta_C)$ of principally polarized abelian varieties where $E$ is a very general elliptic curve. 
By Theorem \ref{theorem:verygeneral-powers-elliptic-isogenies} in Appendix \ref{appendix:B}, we conclude that $k \leq 11$. 

For $8 \leq k \leq 11$, the only indecomposable integral inner product space is, by Kneser's classification, given by the $E_8$-lattice, see \cite[p.\ 28, Remark 1]{milnorhusemoller} and \cite{kneser-quadratisch}. 
 Hence we are reduced to the case $k=8$ and we have $(E^8, \alpha) \cong (JC,\Theta_C)$ for a smooth projective curve $C$ of genus eight, with $\alpha$ induced by the $E_8$-lattice. 
 In particular, the automorphism group of $(E^8, \alpha)$ is isomorphic to $W(E_8)$, the Weyl group of type $E_8$, hence $\va{\Aut(E,\alpha)} = \va{W(E_8)} = 4! \cdot 6! \cdot  8!= 696\,729\,600$ by \cite[Section 2.12]{humphreys-reflectiongroups}. 
 Since $(E^8,\alpha) \cong (JC,\Theta_C)$ as principally polarized abelian varieties, the Torelli theorem implies $\Aut(C) = W(E_8)$ or $\Aut(C) \times \langle \pm{1} \rangle = W(E_8)$. 
 This is absurd: the genus $g(C)$ of $C$ is equal to eight, hence $\va{\Aut(C)} \leq 84(g(C)-1)  = 84 \cdot 7=588$. 
This contradiction concludes the proof of the theorem.
\end{proof}

\section{Modules and lattice theory} \label{section:linear-algebra}

The goal of this section is to prove Lemmas \ref{lemma:matricessaturatedsubmodules} and \ref{lemma:appendix-k=arbitrary-new} below. 
We consider unimodular symplectic lattices
$M$ and $H$ such that $M \subset H^{\oplus k}$ as well as four matrices $\alpha_i \in \rm{M}_k(\Z)$  for $i=1,2,3,4$.
In Lemma \ref{lemma:matricessaturatedsubmodules} we give sufficient conditions that guarantee the inclusions $\alpha_i H^{\oplus k} \subset M$, and in Lemma \ref{lemma:appendix-k=arbitrary-new} we give sufficient conditions for the inclusions $\alpha_i H^{\oplus k} \subset M$ to be an equality. When combined, these lemmas provide a key technical step in the proof of our main theorem, which shall be provided in the next section, see Section \ref{section:maintheorem}. 

\subsection{Preliminary lemmas} 

We start by collecting three basic lemmas for future reference. 

\begin{lemma} \label{lemmma:images-freemodules:linear-algebra}
    Let $\Lambda$ be a free $\Z$-module of positive rank. Let $n \in \Z_{\geq 1}$ and let $\alpha_1, \alpha_2 \in \rm{M}_n(\Z)$ be matrices with non-zero determinant. The following assertions are equivalent:
    \begin{enumerate}
        \item \label{itemlemma:1} We have $\alpha_1 \cdot \Lambda^{\oplus n} = \alpha_2 \cdot \Lambda^{\oplus n}$, for the natural actions of the $\alpha_i$ on $\Lambda^{\oplus n}$. 
        \item \label{itemlemma:2} There exists an invertible matrix $\gamma \in \GL_n(\Z)$ such that $\alpha_2 = \alpha_1 \gamma \in \rm{M}_n(\Z)$. 
    \end{enumerate} 
\end{lemma}
\begin{proof}
Clearly, \eqref{itemlemma:2} implies \eqref{itemlemma:1}. To prove the other implication, assume that \eqref{itemlemma:1} holds. 
Define $\gamma = \alpha_1^{-1}\alpha_2 \in \GL_{n}(\Q)$. One readily shows that $\gamma \in \GL_n(\Z)$. 
    \end{proof}

\begin{lemma} \label{lemma:saturated-matrices-embedding}
Let $M$ be a free $\Z$-module of finite rank and $N \subset M$ a saturated submodule. For $k \geq 1$, let $\alpha \in \rm{M}_k(\Z)$ be a matrix with non-zero determinant. Then $N^{\oplus k} \cap (\alpha \cdot M^{\oplus k}) = \alpha \cdot N^{\oplus k}$. 
\end{lemma}
\begin{proof}
We may assume that $N \neq M$. In particular, $M/N$ is a non-zero free $\Z$-module. Moreover, the matrix $\alpha \in \rm{M}_k(\Z)$ induces an endomorphism
\begin{align} \label{align:endomorphism-embedding-alpha}
\alpha \colon (M/N)^{\oplus k} \to (M/N)^{\oplus k}. 
\end{align}
As $\alpha$ has non-zero determinant, and $(M/N)^{\oplus k}$ is torsion-free, the endomorphism \eqref{align:endomorphism-embedding-alpha} is injective. 
Therefore, $N^{\oplus k} \cap (\alpha \cdot M^{\oplus k}) \subset \alpha \cdot N^{\oplus k}$. The other inclusion is clear. 
\end{proof}

\begin{lemma} \label{lemma:basic}
    Let $\Lambda$ be a $\Z$-module. Let $M_1,M_2$ and $N$ be submodules of $\Lambda$ with $M_1 \subset M_2$. Suppose that $M_1 + N = M_2 + N$, and that $M_1 \cap N = M_2 \cap N$. Then $M_1 = M_2 \subset \Lambda$. 
\end{lemma}
\begin{proof}
It suffices to prove that $M_2 \subset M_1$. Let $x_2 \in M_2$. The hypotheses imply that there exists an element $z \in N$ such that $x_2 = x_1 + z$ for some $x_1 \in M_1$. As $M_1 \subset M_2$, we have $x_1 \in M_2$, hence $x_2 - x_1 = z \in M_2 \cap N = M_1 \cap N$. Therefore, $z \in M_1$, so that $x_2 = x_1 + z \in M_1$. 
\end{proof}

\subsection{Matrices and saturated submodules} \label{subsec:matrix-saturated}
 Consider finite free $\Z$-modules $M$ and $H$ such that $M \subset H^{\oplus k}$. For matrices $\alpha_i \in \rm{M}_k(\Z)$ ($i = 1,2,3,4$), we would like to know whether $\alpha_i H^{\oplus k} \subset M$, assuming that this holds in certain subquotients of $H^{\oplus k}$. 
The goal of this section is to provide some sufficient conditions. 
The main result in this direction is Lemma \ref{lemma:matricessaturatedsubmodules} below. 

\begin{lemma} \label{lemma:preliminary-three}
    Let $H$ be a free $\Z$-module of finite rank. Let $W \subset H$ be a submodule, and let $V_i \subset H$ be saturated submodules for $i = 1,2$, such that $V_1 \cap V_2 = 0$ and $V_1 \oplus V_2 \subsetneq W$. Let $M \subset H^{\oplus k}$ be a submodule such that  for each $i \in \set{1,2}$ and some $\alpha_i \in \rm{M}_k(\Z)$ with non-zero determinant, we have
    \begin{align} \label{align:fundamental-assumption}
    M \cap W^{\oplus k} \equiv \alpha_i \cdot W^{\oplus k} \mod V_i^{\oplus k} .
    \end{align} 
    Then the following holds: \begin{enumerate}
        \item \label{item:firstitem:WH}There exists $\gamma \in \GL_k(\Z)$ such that $\alpha_2 = \alpha_1 \gamma$.
        \item \label{item:seconditem:WH}We have 
        \begin{align}
            \label{align:equality-submodules-k}
            M \cap W^{\oplus k} = \alpha_1 \cdot W^{\oplus k} = \alpha_2 \cdot W^{\oplus k}. 
        \end{align}
    \end{enumerate}
\end{lemma}

\begin{proof}
Replacing $M$ by $M \cap W^{\oplus k}$ and $H$ by $W$, we may assume that $M \subset W^{\oplus k} = H^{\oplus k}$. 
Notice that
\begin{align} \label{align:noticethat-align}
M \equiv \alpha_1 \cdot W^{\oplus k} \equiv \alpha_2 \cdot W^{\oplus k} \quad \mod V_1^{\oplus k} \oplus V_2^{\oplus k}.
\end{align}
As $(W/(V_1\oplus V_2))^{\oplus k}$ is torsion-free, item \eqref{item:firstitem:WH} follows from Lemma \ref{lemmma:images-freemodules:linear-algebra} because $V_1 \oplus V_2 \neq W$. 

Next, let us prove item \eqref{item:seconditem:WH}. Note that since $\alpha_2 = \alpha_1\gamma$ for some $\gamma \in \GL_k(\Z)$ by item \eqref{item:firstitem:WH}, we may and do assume that $\alpha_2 = \alpha_1$. 
We claim that
\begin{align} \label{align:M-inclusion}
M 
\subset \alpha \cdot W^{\oplus k}, \quad \quad \text{where} \quad \quad \alpha \coloneqq \alpha_1 = \alpha_2.  
\end{align}
To prove this, let $x \in M \subset W^{\oplus k}$. By \eqref{align:fundamental-assumption},  
we can write 
$
x = \alpha \cdot w_1 + v_1 = \alpha \cdot w_2 + v_2 \in W^{\oplus k}$ with $w_1, w_2 \in W^{\oplus k}$ and $ v_i \in V_i^{\oplus k}$. 
Therefore, we have
\[
\alpha \cdot \left( w_1 - w_2 \right) = - v_{1} + v_{2} \in \alpha \cdot W^{\oplus k} \cap \left(V_1^{\oplus k} \oplus V_2^{\oplus k} \right) = \alpha \left( V_1^{\oplus k} \oplus V_2^{\oplus k}\right) = \alpha V_1^{\oplus k} \oplus \alpha V_2^{\oplus k},
\] 
where we used Lemma \ref{lemma:saturated-matrices-embedding}, which applies because $V_1^{\oplus k} \oplus V_2^{\oplus k} $ is saturated in $H^{\oplus k}$ by assumption.
In particular, $v_1 \in \alpha V_1^{\oplus k}$ and $v_2 \in \alpha V_2^{\oplus k}$. Thus, we have $
x = \alpha w_1 + v_1 \in \alpha W^{\oplus k}$. This proves the inclusion \eqref{align:M-inclusion}. 

Furthermore, we claim that 
\begin{align} \label{align:alphaU1M}
\alpha \cdot V_1^{\oplus k} \subset M \cap V_1^{\oplus k}. 
\end{align}
To prove this, notice that  
$
\alpha \cdot V_1^{\oplus k} \subset \alpha \cdot W^{\oplus k} \equiv M \cap W^{\oplus k} \bmod  V_2^{\oplus k}$, where we use the assumption $V_1\subset W$ for the inclusion and \eqref{align:fundamental-assumption} for the congruence.
As $\alpha \cdot V_1^{\oplus k} \subset V_1^{\oplus k}$, we obtain
\begin{align} \label{align:alphaU1M:congruence}
\alpha \cdot V_1^{\oplus k} \subset M \cap V_1^{\oplus k} \quad \mod V_2^{\oplus k}. 
\end{align}
Notice that \eqref{align:alphaU1M:congruence} implies \eqref{align:alphaU1M}, because $V_1 \cap V_2 = 0$. Our claim is proved.  
In a similar way (or by symmetry), one proves that $\alpha \cdot V_2^{\oplus k} \subset M \cap V_2^{\oplus k} \subset M$. 

As $V_1 \cap V_2 = 0$, it follows that $\alpha\cdot V_1^{\oplus k} \oplus \alpha \cdot V_2^{\oplus k} \subset M$. Via Lemma \ref{lemma:saturated-matrices-embedding}, we thus obtain:
\begin{align}
    \label{align:MU1U2:i}
\alpha \cdot W^{\oplus k} \cap \left(V_1^{\oplus k} \oplus V_2^{\oplus k} \right) =
\alpha \cdot \left(V_1^{\oplus k} \oplus V_2^{\oplus k} \right)   
\subset M \cap \left(V_1^{\oplus k} \oplus V_2^{\oplus k} \right). 
\end{align}
By \eqref{align:M-inclusion}, we have $M \subset \alpha \cdot W^{\oplus k}$, so that 
\begin{align}
\label{align:MU1U2:ii}
M \cap \left(V_1^{\oplus k} \oplus V_2^{\oplus k} \right) \subset \alpha \cdot W^{\oplus k} \cap \left(V_1^{\oplus k} \oplus V_2^{\oplus k} \right).
\end{align}
Combining \eqref{align:MU1U2:i} and \eqref{align:MU1U2:ii}, we see that
\begin{align}
\label{align:MU1U2:iii}
M \cap \left(V_1^{\oplus k} \oplus V_2^{\oplus k} \right) = 
\alpha \cdot W^{\oplus k} \cap \left( V_1^{\oplus k} \oplus V_2^{\oplus k} \right). 
\end{align}
We are now in position to apply Lemma \ref{lemma:basic} to the following $\Z$-modules: let $\Lambda \coloneqq W^{\oplus k}$, $M_1 \coloneqq M \subset \Lambda$, $M_2 \coloneqq \alpha \cdot W^{\oplus k} \subset \Lambda$, and $N \coloneqq V_1^{\oplus k} \oplus V_2^{\oplus k}$. Observe that \eqref{align:noticethat-align}, \eqref{align:M-inclusion} and \eqref{align:MU1U2:iii} imply respectively that $M_1 + N = M_2 + N$, $M_1 \subset M_2$, and $M_1 \cap N = M_2 \cap N$. Therefore, by Lemma \ref{lemma:basic}, we have $M_1 = M_2$. That is, $
M = \alpha \cdot W^{\oplus k}$, and the lemma follows. 
\end{proof}

\begin{lemma} \label{lemma:matricessaturatedsubmodules}
Let $H$ be a free $\Z$-module of finite rank. Assume that, for each $i \in \set{1,2,3,4}$, there exists a sequence of free submodules 
$$
V_i \subset W_i  \subset H,
$$
such that the following properties are satisfied:
\begin{enumerate} 
\item \label{item:embeddingTFcokernel} The natural map $V_1 \oplus V_2 \oplus V_3 \oplus V_4 \to H$ is an embedding with torsion-free cokernel.  
\item \label{item:H1H2H} We have 
$(W_1 \cap W_2) + (W_3 \cap W_4) = H$.
\item \label{item:nonemptysaturatedintersection} The intersections $$W_1 \cap {W_2}, \quad \quad  {W_3} \cap {W_4}, \quad \quad \text{and} \quad \quad W_1 \cap W_2 \cap W_3 \cap W_4$$ are non-zero and saturated in $H$, and $W_1 \cap W_2$ (resp.\ $W_3 \cap W_4$, resp.\ $W_1\cap W_2 \cap W_3 \cap W_4$) strictly contains $V_1 \oplus V_2$ (resp.\ $V_3 \oplus V_4$, resp.\ $V_1 \oplus V_2 \oplus V_3 \oplus V_4$).  
\end{enumerate}
Let $k\geq 1$ and let $M \subset H^{\oplus k}$ be a submodule that satisfies the following condition: for each $i \in \set{1,2,3,4}$ there exists a matrix $\alpha_i \in \rm{M}_k(\Z)$ with non-zero determinant, such that 
\begin{align} \label{hypothesis:MVik}
M \cap W_i^{\oplus k} \equiv \alpha_i \cdot W_i^{\oplus k} \mod V_i^{\oplus k}. 
\end{align}
Then for each $i,j \in \set{1,2,3,4}$, 
we have:
\begin{align} \label{align:MV1V2-equality:app} 
\alpha_i \cdot H^{\oplus k} = \alpha_j \cdot H^{\oplus k} \subset M.  
\end{align} 
\end{lemma}

\begin{proof}  
We claim that the following equality holds:
\begin{align} \label{claim:alphaX12}
\alpha_1 \cdot \left( (W_1 \cap W_2 )/ ( V_1 \oplus V_2) \right)^{\oplus k} = \alpha_2 \cdot \left( (W_1 \cap W_2 )/ ( V_1 \oplus V_2) \right)^{\oplus k}. 
\end{align}
To prove this, first observe that $(W_1 \cap W_2 )/ ( V_1 \oplus V_2)$ is saturated in $W_i/(V_1 \oplus V_2)$ for $i = 1,2$, and that both modules are non-zero and torsion-free. Indeed, by condition \eqref{item:embeddingTFcokernel}, $V_1\oplus V_2$ is saturated in $H$ and hence in any submodule that contains it; by condition \eqref{item:nonemptysaturatedintersection}, we have $V_1\oplus V_2\subsetneq W_1 \cap W_2 \subset W_i$ and hence the quotients $(W_1 \cap W_2 )/ ( V_1 \oplus V_2)$ and $W_i/(V_1 \oplus V_2)$ are non-zero and torsion-free. The saturation of $(W_1 \cap W_2 )/ ( V_1 \oplus V_2) \subset W_i/(V_1 \oplus V_2)$ follows from the saturation of $W_1 \cap W_2 \subset W_i$, which holds because of condition \eqref{item:nonemptysaturatedintersection}. 
In view of Lemma \ref{lemma:saturated-matrices-embedding}, we deduce that 
\begin{align} \label{align:added-1}
    \alpha_1 \cdot \left( W_1^{\oplus k} \cap W_2^{\oplus k} \right) &\equiv \left(\alpha_1 \cdot W_1^{\oplus k}\right) \cap W_2^{\oplus k}  \quad \quad  \mod V_1^{\oplus k} \oplus V_2^{\oplus k},\\
    \label{align:added-2}
    W_1^{\oplus k}  \cap \left( \alpha_2 \cdot W_2^{\oplus k}\right)   &\equiv \alpha_2 \cdot \left( W_1^{\oplus k} \cap W_2^{\oplus k} \right) \quad \quad  \mod V_1^{\oplus k} \oplus V_2^{\oplus k}.
\end{align}
Moreover, because of \eqref{hypothesis:MVik}, we have:
\begin{align} \label{align:added-3}
\begin{split}
    \left(\alpha_1 \cdot W_1^{\oplus k}\right) \cap W_2^{\oplus k}  
&\equiv 
\left( M \cap W_1^{\oplus k} \right) \cap W_2^{\oplus k}
\equiv
W_1^{\oplus k}  \cap \left( M \cap W_2^{\oplus k}\right) 
\\
&
\equiv W_1^{\oplus k}  \cap \left( \alpha_2 \cdot W_2^{\oplus k}\right) \quad \quad
\mod V_1^{\oplus k} \oplus V_2^{\oplus k}.
\end{split}
\end{align}
Taken together, \eqref{align:added-1}, \eqref{align:added-2} and \eqref{align:added-3} imply \eqref{claim:alphaX12}, proving the claim. 

From Lemma \ref{lemmma:images-freemodules:linear-algebra} and equation \eqref{claim:alphaX12}, we conclude that there exists $\gamma_{12} \in \GL_k(\Z)$ such that $\alpha_2 = \alpha_1 \gamma_{12} \in \rm{M}_k(\Z)$. By symmetry, there exists a matrix $\gamma_{34} \in \GL_k(\Z)$ such that $\alpha_4 = \alpha_3\gamma_{34}$. Replacing $\alpha_2$ by $\alpha_2 \gamma_{12}^{-1} = \alpha_1$ and $\alpha_4$ by $\alpha_4\gamma_{34}^{-1} = \alpha_3$, we may (and will) assume that $\alpha_1 = \alpha_2$ and $\alpha_3 = \alpha_4$. 
Consider, for $i = 1,2,3,4$, the inclusions
\[
V_1 \oplus V_2 \oplus V_3 \oplus V_4 \subset W_1 \cap W_2 \cap W_3 \cap W_4 \subset W_i. 
\]
Both inclusions are saturated by conditions \eqref{item:embeddingTFcokernel} and \eqref{item:nonemptysaturatedintersection}. Consequently, for each $i$, the inclusion
\[
\left(W_1 \cap W_2 \cap W_3 \cap W_4 \right)/\left(V_1 \oplus V_2 \oplus V_3 \oplus V_4 \right) \subset W_i/\left(V_1 \oplus V_2 \oplus V_3 \oplus V_4 \right)
\]
is a saturated embedding of free $\Z$-modules, which are non-zero in view of condition \eqref{item:nonemptysaturatedintersection}. 

We obtain the following congruences:
\begin{align*}
\alpha_1 \cdot (W_1 \cap W_2 \cap W_3 \cap W_4)^{\oplus k} &\equiv \left( \alpha_1 \cdot W_1^{\oplus k} \right) \cap \left(W_1 \cap W_2 \cap W_3 \cap W_4\right)^{\oplus k} 
\\
& \equiv \left(M \cap  W_1^{\oplus k}\right) \cap W_2^{\oplus k} \cap W_3^{\oplus k} \cap W_4^{\oplus k}   \\
&\equiv W_1^{\oplus k} \cap W_2^{\oplus k} \cap \left(M\cap W_3^{\oplus k}\right) \cap W_4^{\oplus k} \\
& \equiv W_1^{\oplus k} \cap W_2^{\oplus k} \cap \left(\alpha_3 \cdot W_3^{\oplus k}\right) \cap W_4^{\oplus k} \\
& \equiv  \alpha_3 \cdot \left(W_1 \cap W_2 \cap W_3 \cap W_4\right)^{\oplus k} \  \mod V_1^{\oplus k} \oplus V_2^{\oplus k} \oplus V_3^{\oplus k} \oplus V_4^{\oplus k}. 
\end{align*}
Here, the first congruence follows from Lemma \ref{lemma:saturated-matrices-embedding}, the second from equation \eqref{hypothesis:MVik}, the third congruence is clear, the fourth congruence follows from \eqref{hypothesis:MVik} again and the final congruence from Lemma \ref{lemma:saturated-matrices-embedding} again. 

Therefore, by Lemma \ref{lemmma:images-freemodules:linear-algebra} using that $(W_1 \cap W_2 \cap W_3 \cap W_4)/(V_1 \oplus V_2 \oplus V_3 \oplus V_4)$ is a non-zero free $\Z$-module, there exists an invertible matrix $\gamma \in \GL_k(\Z)$ such that $\alpha_3 = \alpha_1 \gamma$. Let $\alpha = \alpha_1$. 
Replacing $\alpha_3$ by $\alpha_3\gamma^{-1} = \alpha_1$, we may and do assume that 
\begin{align} \label{equalities:alpha}
\alpha = \alpha_1 = \alpha_2 = \alpha_3 = \alpha_4. 
\end{align}
By \eqref{hypothesis:MVik} and \eqref{equalities:alpha}, we get:
\begin{align} \label{equalities:MmodW}
M \cap W_i^{\oplus k} &\equiv \alpha \cdot W_i^{\oplus k} \mod V_i^{\oplus k} \quad \quad  \forall i \in \set{1,2}. 
\end{align}
Combining \eqref{equalities:MmodW} with Lemma \ref{lemma:saturated-matrices-embedding}, using that, for $i = 1,2$, $(W_1 \cap W_2)/(V_1 \oplus V_2) \subset W_i/(V_1 \oplus V_2)$ is a saturated inclusion of free $\Z$-modules (see conditions \eqref{item:embeddingTFcokernel} and \eqref{item:nonemptysaturatedintersection}), we obtain:
\begin{align} \label{equalities:MmodW:congruence2}
M \cap \left(W_1 \cap W_2\right)^{\oplus k} &\equiv \alpha \cdot \left(W_1 \cap W_2\right)^{\oplus k} \mod V_i^{\oplus k} \quad \forall i \in \set{1,2}. 
\end{align}
By Lemma \ref{lemma:preliminary-three} and the fact that $V_1 \cap V_2 = 0$, it follows from \eqref{equalities:MmodW:congruence2} that
\[
M \cap \left(W_1 \cap W_2 \right)^{\oplus k} = \alpha \cdot \left(W_1 \cap W_2 \right)^{\oplus k}. 
\]
By symmetry, we obtain $
M \cap \left(W_3 \cap W_4 \right)^{\oplus k} = \alpha \cdot \left(W_3 \cap W_4 \right)^{\oplus k} 
$. 

Finally, as we have an equality $H^{\oplus k} = (W_1 \cap W_2 )^{\oplus k}  + (W_3 \cap W_4 )^{\oplus k}$, see condition \eqref{item:H1H2H}, we conclude that
\begin{align*}
\alpha \cdot H^{\oplus k} = \alpha \cdot 
\left(W_1 \cap W_2 \right)^{\oplus k} + \alpha \cdot \left(W_3 \cap W_4 \right)^{\oplus k} 
= 
M \cap \left(W_1 \cap W_2 \right)^{\oplus k} + M \cap \left(W_3 \cap W_4 \right)^{\oplus k} \subset M.
\end{align*}
In other words, the inclusion \eqref{align:MV1V2-equality:app} is proved, and we are done. 
\end{proof}

\subsection{Sublattices of powers of unimodular symplectic lattices} \label{subsection:sublattices-of-powers}

We continue by investigating sublattices $M \subset H^{\oplus k}$ of some power of a unimodular symplectic lattice $H$. This section is independent of Section \ref{subsec:matrix-saturated}. The main result of this section is Lemma \ref{lemma:appendix-k=arbitrary-new}. 

\begin{lemma}\label{lem:isotropic-subspace}
Let $(H,E_H)$ be a unimodular symplectic lattice of positive finite rank. 
Let $U\subset H$ be a saturated isotropic subspace.
Then there is a saturated subspace $U'\subset H$ such that $U\oplus U'\subset H$ is a unimodular subspace.
\end{lemma}
\begin{proof}
We argue by induction on the rank of $U$.
If $U$ has rank one, then $U=\langle u\rangle $.
The map $E_H(u,-):H\to \Z$ is surjective because $u$ is indivisible and $E_H$ is unimodular.
Hence there is a class $u'\in H$ with $E_H(u,u')=1$, as we want.

 If $U$ has rank $r\geq 2$, then we pick a saturated subspace $U_1\subset U$ of rank $r-1$ and apply the induction hypothesis to $U_1$ to get a unimodular subspace $U_1\oplus U_1'\subset H$.
 Any unimodular subspace of a unimodular lattice admits a unimodular complement.
 Hence, we can extend any symplectic basis of $U_1\oplus U_1'$ to a symplectic basis of $H$.
 It follows that there is a symplectic basis
 $$
\{  e_1,\dots ,e_g,f_1,\dots ,f_g \}
$$
of $H$ with $U=\langle e_1,\dots ,e_{r-1}\rangle\oplus \langle u\rangle  $ for some primitive element
$
u\in \langle e_r,\dots ,e_{g},f_r,\dots ,f_g\rangle .
$
Applying the $r=1$ case of the lemma to the subspace $ \langle u\rangle\subset \langle e_r,\dots ,e_{g},f_r,\dots ,f_g\rangle $, we get
 an element $u'\in  \langle e_r,\dots ,e_{g},f_r,\dots ,f_g\rangle $ such that $$ 
\langle u,u'\rangle\subset  \langle e_r,\dots ,e_{g},f_r,\dots ,f_g\rangle $$
is a unimodular subspace.  
This subspace has a complement and hence up to change of basis we can assume $u=e_r$ and $u'=f_r$.
At this point the lemma is clear.
\end{proof}

\begin{lemma} \label{lemma:appendix:preliminary-lemma}
Let $(H,E_H)$ be a unimodular symplectic lattice of rank $2g$. 
Let $k\geq 1$ and consider the induced unimodular symplectic lattice $(H^{\oplus k},E_{H}^{\oplus k})$. 
Let $M\subset H^{\oplus k}$ be a sublattice such that there is a matrix $\beta\in \rm{M}_k(\Z)$ with positive determinant such that the intersection form $$
E_M(-,-)\coloneqq E_H^{\oplus k}(\beta^{-1}-,-)
$$ 
is unimodular and integral on $M$. Assume that there is a matrix $\alpha\in \rm{M}_k(\Z)$ with nonzero determinant such that $\alpha H^{\oplus k}\subset M$. Let $[H^{\oplus k} \colon M]$ be the index of $M$ in $H^{\oplus k}$. Then 
\[
\left[H^{\oplus k} \colon M \right] = \det(\beta)^g  
 \mid \det(\alpha)^{2g}. 
\]
\end{lemma}
\begin{proof}
The inclusions $\alpha H^{\oplus k} \subset M \subset H^{\oplus k}$, together with the fact that the index of $\alpha H^{\oplus k}$ in $H^{\oplus k}$ equals $\det(\alpha)^{2g}$, show that the index of $M$ in $H^{\oplus k}$ divides $\det(\alpha)^{2g}$. Thus, it suffices to prove that $[H^{\oplus k} \colon M] = \det(\beta)^g$. 

To prove this, note that $E_M$ and $\psi$ induce isomorphisms $M \cong M^\vee$ and $H \cong H^\vee$. In particular, we can dualize the inclusion $M \subset H^{\oplus k}$ to obtain an inclusion $H^{\oplus k} \subset M$, and, with respect to this embedding, the index of $M$ in $H^{\oplus k}$ equals the index of $H^{\oplus k}$ in $M$. 
The fact that $E_H^{\oplus k}(\beta^{-1}\cdot  -, -)$ restricts to the unimodular pairing $E_M$ on $M$ implies that $E_M$ restricts to $E_H^{\oplus k}(\beta \cdot -, -)$ on $H^{\oplus k} \subset M$. This means precisely that the composition
\begin{align}\label{appendix:composition-M-H:preliminary}
\xymatrix{
H^{\oplus k} \ar@{^{(}->}[r] & M \ar[r]^-{\cong} & M^\vee \ar@{^{(}->}[r] & (H^{\oplus k})^\vee = (H^\vee)^{\oplus k} 
}
\end{align}
is given by the map
\[
x \mapsto \left(y \mapsto E_H^{\oplus k}(\beta \cdot x, y) \right). 
\]
Consequently, by identifying $M$ and $H^{\oplus k}$ with their respective duals $M^\vee$ and $(H^{\oplus k})^\vee$ (via $E_M$ and $E_H^{\oplus k}$) in \eqref{appendix:composition-M-H:preliminary}, it follows that the composition
\[
\xymatrix{
H^{\oplus k} \ar@{^{(}->}[r] & M \ar@{^{(}->}[r] & H^{\oplus k}
}
\]
is given by multiplication by the matrix $\beta$. In particular, if $[H^{\oplus k} \colon M]$ denotes the index of $M$ in $H^{\oplus k}$ and $[ H^{\oplus k} \colon \beta \cdot H^{\oplus k}]$ denotes the index of $\beta \cdot H^{\oplus k}$ in $H^{\oplus k}$, then
\[
\left[ H^{\oplus k} \colon M \right]^2= \left[ H^{\oplus k} \colon \beta \cdot H^{\oplus k}\right] = \det(\beta)^{2g}.
\]
As $\det(\beta)$ is positive, we get $[H^{\oplus k} \colon M] = \det(\beta)^g$ as desired. 
\end{proof}

\begin{lemma} \label{lemma:linearalgebra:determinantsdivision}
    Let $\Lambda$ be a free $\Z$-module and put $\Lambda_\Q = \Lambda \otimes \Q$. 
    Let $\psi \colon \Lambda_\Q \times \Lambda_\Q \to \Q$ be a bilinear form that takes integral values on $\Lambda \times \Lambda$. 
    Let $\alpha_i, \beta_i \in \rm{M}_k(\Z)$ $(i = 1,2)$ be matrices with non-zero determinant. Let $e_1, \dotsc, e_k, f_1, \dotsc, f_k \in \Lambda$. For $\gamma = (\gamma_{ij}) \in \rm{M}_k(\Q)$, define $\gamma e_i = \sum_{j} \gamma_{ji}e_j$ and $\gamma f_i = \sum_{j} \gamma_{ji} f_j$, and assume $\psi(\beta^{-1}_1\alpha_1 e_i, \beta_2^{-1}\alpha_2 f_j) = \delta_{ij}$ for all $i,j \in \set{1,\dotsc, k}$. 
    
    Then $\det(\alpha_1)\det(\alpha_2) \mid \det(\beta_1)\det(\beta_2)$. 
\end{lemma}

\begin{proof}
    Let $\set{h_1, \dotsc, h_g}$ be the canonical basis of $\Z^{\oplus k}$, and define a bilinear form 
\[
\Psi \colon \Q^{\oplus k} \times \Q^{\oplus k} \to \Q \quad \quad  \text{ by } \quad\quad   \Psi(h_i, h_j) = \psi(e_i, f_j).  
\]
Then $\Psi$ takes integral values on $\Z^{\oplus k} \times \Z^{\oplus k}$. Let $Q = (Q_{ij})$ be the $k \times k$-matrix with entries $Q_{ij} = \Psi(h_i,h_j)$. Then $Q \in \rm{M}_k(\Z)$, and for $x,y \in \Q^{\oplus k}$, we have
\[
\Psi(x,y) = x^\top Q y. 
\]
Moreover, $\Psi( \beta^{-1}_1 \alpha_1 h_i, \beta_2^{-1}\alpha_2 h_j) = \psi(\beta^{-1}_1\alpha_1 e_i, \beta_2^{-1}\alpha_2 f_j) = \delta_{ij}$. Therefore,
\[
\delta_{ij} = 
\Psi( \beta^{-1}_1 \alpha_1 h_i, \beta_2^{-1}\alpha_2 h_j) = 
(\beta^{-1}_1\alpha_1 h_i)^\top \cdot Q  \cdot (\beta_2^{-1}\alpha_2 h_j) = 
h_i^\top \left( \alpha_1^\top \beta_1^{-\top} Q \beta_2^{-1}\alpha_2 \right) h_j. 
\]
Hence, the matrix $\alpha_1^\top \beta_1^{-\top} Q \beta_2^{-1}\alpha_2\in \rm{M}_k(\Q)$ is the identity matrix. In particular, 
\[
\det(\alpha_1^\top \beta_1^{-\top} Q \beta_2^{-1}\alpha_2) = \det(\alpha_1)\det(\alpha_2)\det(\beta_1)^{-1}\det(\beta_2)^{-1}\det(Q) = 1,
\]
which implies that 
$
\det(\alpha_1)\det(\alpha_2) \det(Q) = \det(\beta_1)\det(\beta_2)$. We have $\det(Q) \in \Z$ because $Q \in \rm{M}_k(\Z)$, and the lemma follows. 
\end{proof}

We come now to the main result of Section \ref{subsection:sublattices-of-powers}, which is the following lemma. 
A proof of a simpler version of this lemma is contained in the appendix, 
see Proposition \ref{proposition:marcucci:improved}. 

Let $(H,E_H)$ be a unimodular symplectic lattice of rank $2g$.
Let $k\geq 1$ and consider the unimodular lattice $(H,E_H)^{\oplus k}=(H^{\oplus k},E_H^{\oplus k})$ with symplectic basis 
$$
\set{ e_{1i},\dots ,e_{gi},f_{1i},\dots ,f_{gi}}, \ \ \ \ i=1,\dots ,k.
$$
Let 
$$
H'_1\coloneqq \bigoplus_{i=1}^k \langle e_{1i},\dots ,e_{gi},f_{2i},f_{3i},\dots ,f_{gi}\rangle, \quad \quad 
H'_2\coloneqq \bigoplus_{i=1}^k \langle e_{1i},\dots ,e_{gi},f_{1i},f_{3i},\dots ,f_{gi}\rangle ,
$$
and define
$$
V_1 \coloneqq\bigoplus_{i = 1}^k \langle e_{1i} \rangle, \quad \quad 
V_2 \coloneqq \bigoplus_{i = 1}^k \langle e_{2i} \rangle. 
$$

\begin{lemma} \label{lemma:appendix-k=arbitrary-new}
In the above notation, let $M\subset H^{\oplus k}$ be a sublattice such that there is a matrix $\beta\in \rm{M}_k(\Z)$ with positive determinant such that the intersection form $$
E_M(-,-)\coloneqq E_H^{\oplus k}(\beta^{-1}-,-)
$$ 
is unimodular and integral on $M$. Assume that there are matrices $\alpha_1, \alpha_2 \in \rm{M}_k(\Z)$ with nonzero determinant such that $\alpha_i H^{\oplus k}\subset M$ for $i = 1,2$, and such that \begin{align}\label{align:appendix:mcapH'modV:k=arbitrary}M\cap H'_i\equiv \alpha_i H'_i \mod V_i, \quad \quad i = 1,2.\end{align}
Then there exists $\gamma \in \GL_k(\Z)$ such that $\alpha_2 = \alpha_1 \gamma$, and we have $(\det \alpha_i)^2= \det \beta$ for $i = 1,2$. Furthermore, we have $M=\alpha_i H^{\oplus k}$ for $i = 1,2$. 
\end{lemma}

\begin{proof}
As a first step, we aim to show that there exists $\gamma \in \GL_k(\Z)$ such that $\alpha_2 = \alpha_1 \gamma$. To prove this, notice that by Lemma \ref{lemma:saturated-matrices-embedding}, we have, for $i = 1,2$: 
\begin{align}\label{equality:sixone}
\left(\alpha_i \cdot H_1'\right) \cap H_2'
=\alpha_i \cdot \left( H_1' \cap H_2' \right) =
H_1'  \cap \left( \alpha_i \cdot H_2'\right). 
\end{align}
Moreover, $V_1 \oplus V_2 \subset H_1' \cap H_2'$. Modulo $V_1 \oplus V_2$, the equality \eqref{equality:sixone} combined with \eqref{align:appendix:mcapH'modV:k=arbitrary} gives 
\begin{align} \label{equation:alpha1alpha2}
\begin{split}
\alpha_1 \cdot \left( H_1' \cap H_2' \right) &\equiv \left(\alpha_1 \cdot H_1'\right) \cap H_2'  
\equiv \left( M \cap H_1'\right) \cap H_2' 
\equiv H_1'  \cap \left( M \cap H_2'\right) 
\\
 &\equiv H_1'  \cap \left( \alpha_2 \cdot H_2'\right)    \equiv \alpha_2 \cdot \left( H_1' \cap H_2' \right) \quad \quad  \mod V_1 \oplus V_2.  
\end{split}
\end{align}
Remark that $V_1 \oplus V_2 \subset H_1' \cap H_2'$ is saturated in $H_1' \cap H_2'$. Thus, by Lemma \ref{lemmma:images-freemodules:linear-algebra}, equation \eqref{equation:alpha1alpha2} implies that there exists $\gamma \in \GL_k(\Z)$ such that $\alpha_2 = \alpha_1\gamma$, as we want.

As a second step, we let
    $$
    U\coloneqq \alpha \cdot \left(V_1 \oplus V_2 \right) \subset M, \quad \quad \alpha \coloneqq \alpha_1 \in \rm{M}_k(\Z),
    $$
and claim that $U$ is saturated in $M$. To prove this, let $x \in M$ and assume that $rx \in U$ for some $r \in \Z_{\geq 1}$. We need to show that $x \in U$. As $rx \in U$, we have $x \in V_1 \oplus V_2$. In particular, $x \in H_i'$ for $i = 1,2$, hence $x \in M \cap  H_1' \cap H_2'$. By \eqref{align:appendix:mcapH'modV:k=arbitrary}, this means that we can write 
\begin{align*}
    x & = \alpha h_1' + v_1, \quad \quad \text{for some} \quad \quad h_1' \in H_1'\quad \text{and} \quad v_1 \in V_1, \\
    x & = \alpha h_2' + v_2, \quad \quad \text{for some} \quad \quad h_2' \in H_2' \quad \text{and} \quad v_2 \in V_2.
\end{align*}
As $x \in V_1 \oplus V_2$, we get that $h_i' \in V_1 \oplus V_2$ and hence 
$$
\alpha(h_1' - h_2') =  -v_1 + v_2 \in  \langle e_1, e_2 \rangle^{\oplus k}=V_1 \oplus V_2.
$$
As the action of $\alpha$ respects the decomposition $ V_1 \oplus V_2$, this implies $v_1 \in \alpha V_1$ and $v_2 \in \alpha V_2$. Therefore, 
$
x  = \alpha h_1' + v_1 \in \alpha \langle e_1, e_2 \rangle^{\oplus k} = U,
$
proving our claim that $U$ is saturated in $M$. 

Next, observe that the subspace $U$ of $M$ is isotropic. As we have just proved that $U \subset M$ is saturated in $M$, we conclude from Lemma \ref{lem:isotropic-subspace} that there are classes
$
g_{1i},g_{2i}\in M
$
such that
$$
 U \oplus \bigoplus_{i=1}^k \langle g_{1i},g_{2i} \rangle 
$$
is a unimodular sublattice of $M$.
More precisely, we can assume that the $g_{ai}$ ($a \in \set{1,2}$) are chosen such that
$
E_M(\alpha e_{ai}, g_{bj})=\delta_{ab}\delta_{ij} . 
$
Hence
\begin{align} \label{align:appendix-deltaij}
E_H^{\oplus k}(\beta^{-1}\alpha e_{ai}, g_{bj})=\delta_{ab}\delta_{ij} . 
\end{align}

Since the $e_{rs},f_{rs}$ form a basis of $H^{\oplus k}$, we can write
$$
g_{ij}=\sum_{r,s} \left(a_{ijrs}e_{rs}+b_{ijrs}f_{rs} \right) 
$$
for uniquely determined $a_{ijrs},b_{ijrs} \in \Z$.
By \eqref{align:appendix-deltaij}, we have 
$$
E_H^{\oplus k}(\beta^{-1}\alpha e_{1i}, g_{2j})=0\ \ \ \text{for all $i,j$.}
$$   
Since $\alpha$ and $\beta$ have nonzero determinants, the classes $\beta^{-1}\alpha e_{1i}$ with $i=1,\dots ,k$ span rationally the space $\oplus_{i=1}^k\langle e_{1i}\rangle $ and so we find that $g_{2j}$ contains the basis element $f_{1i}$ trivially, that is, we have $b_{2j1s}=0$ for all $j,s=1,\dots,k$. This implies 
$
g_{2j}\in M\cap H'_1$ for all $j = 1, \dotsc, k$. Hence, by \eqref{align:appendix:mcapH'modV:k=arbitrary}, there exist 
classes $g_{2j}'\in H'_1$ such that, for all $j$, we have 
\begin{align} \label{align:appendix:delta_ij_ap}
g_{2j}=\alpha g'_{2j} + e_j', \quad  \text{for some} \quad e_j'  \in V_1. 
\end{align}
By \eqref{align:appendix-deltaij}, we have
\begin{align}\label{eq:EH=deltaij-new-ii}
E_H^{\oplus k}(\beta^{-1}\alpha e_{2i},g_{2j})= \delta_{ij} .
\end{align}
Combining \eqref{align:appendix:delta_ij_ap} and \eqref{eq:EH=deltaij-new-ii}, and noticing that $E_H^{\oplus k}(\beta^{-1}\alpha e_{2i}, e_j') = 0$, 
we get
\begin{align*}
\delta_{ij} = E_H^{\oplus k}(\beta^{-1}\alpha e_{2i},g_{2j}) = 
E_H^{\oplus k}(\beta^{-1}\alpha e_{2i},g_{2j} - e'_j)= 
E_H^{\oplus k}(\beta^{-1}\alpha e_{2i}, \alpha g_{2j}'
)
\end{align*}
for all $i,j \in \set{1, \dotsc, k}$. 
By Lemma \ref{lemma:linearalgebra:determinantsdivision}, this implies that $\det(\alpha)^2 \mid \det(\beta)$. 

Conversely, by Lemma \ref{lemma:appendix:preliminary-lemma}, we have that $\det(\beta)^g  
 \mid \det(\alpha)^{2g}$. We conclude that $\det(\alpha)^{2g} = \det(\beta)^g$. As $\det(\beta)$ is positive, it follows that $\det(\beta) = \det(\alpha)^2$. Finally, by Lemma \ref{lemma:appendix:preliminary-lemma} again, we have that $[H^{\oplus k} \colon M] = \det(\beta)^g$. Consequently, $[H^{\oplus k} \colon M] = \det(\alpha)^{2g} = [H^{\oplus k} \colon \alpha H^{\oplus k}]$. As $\alpha H^{\oplus k} \subset M$, we must have $M = \alpha H^{\oplus k}$. This finishes the proof of the lemma. 
\end{proof}

\section{Proof of the main theorem} \label{section:maintheorem}

The goal of this section is to prove Theorem \ref{theorem:maintheorem1:hyperelliptic}, the main result of this paper.
The idea is to degenerate the given isogeny $JC \to (JX)^{ k}$ in four different directions; these are provided by Lemma \ref{lemma:hyperelliptic-degeneration}. 
For each such a degeneration $JC_{0i} \to (JX_{0i})^{k}$ ($i \in \set{1,2,3,4}$), we consider the induced map on the compact quotients, or equivalently the map $H^1(\wt C_{0i},\Z) \to H^1(\wt X_{0i},\Z)^{\oplus k}$ between the cohomology groups of the normalizations  $\wt C_{0i}$ and $\wt X_{0i}$ of $C_{0i}$ and $X_{0i}$. We want to show that, for each $i$, the image of this map is $\alpha_i \cdot H^1(\wt X_{0i},\Z)^{\oplus k} \subset H^1(\wt X_{0i},\Z)^{\oplus k}$ for some $\alpha_i \in \rm{M}_k(\Z)$. 
We then combine these pieces of information obtained in the different degenerations to conclude that $\alpha_j = \alpha_i \gamma_{ij}$ for some $\gamma_{ij} \in \GL_k(\Z)$, and that the image of $H^1(C,\Z) \to H^1(X,\Z)^{\oplus k}$ is given by $ \alpha_i \cdot H^1(X,\Z)^{\oplus k}$. In this last step, we use the linear algebra and lattice theory worked out in Section \ref{section:linear-algebra}. More precisely, this is where Lemmas \ref{lemma:matricessaturatedsubmodules} and \ref{lemma:appendix-k=arbitrary-new} enter the picture.

\subsection{Degeneration of the isogeny: moving the extension class} In this section, we prove:  
    \begin{proposition} \label{proposition:degenerate-isogeny:moving-extension}
    Let $\Delta$ be a connected normal complex analytic space, and let $
p \colon \ca X \to \Delta
$ be a family of semi-stable curves
of arithmetic genus $g \geq 2$. 
Assume there is a family of semi-stable curves $q \colon \mathcal C\to \Delta$ together with, for some $k\geq 1$, an isogeny
\[
\psi \colon (J\mathcal X)^k \coloneqq 
J\mathcal X\times_\Delta\dots \times_\Delta J\mathcal X\longrightarrow J\mathcal C
\]
of semi-abelian schemes over $\Delta$, where $J\mathcal X$ and $J\mathcal C$ denote the respective relative Jacobians over $\Delta$, such that the following conditions are satisfied. 
    \begin{enumerate}[label=(\Roman*)] 
    \item \label{item:prop-constant-normalization:1}
    For all $s \in \Delta$, $X_s= p^{-1}(s)$ is an irreducible one-nodal hyperelliptic curve. 
    \item \label{item:prop-constant-normalization:2}  
    The moduli map $\Delta \to \overline{\ca M}_g$ induced by the family $p \colon \ca X \to \Delta$ is generically finite onto its image, where $\overline{\ca M}_g$ is the moduli stack of stable curves of genus $g$. 
\item \label{item:prop-constant-normalization:3}  
   For any simultaneous normalization $\nu \colon \wt{\ca X} \to \ca X$ of $p$, see Proposition \ref{prop:simultaneous-normalization}, the map from $ \Delta$ to the moduli stack of smooth hyperelliptic curves of genus $g-1$ induced by the family of hyperelliptic curves $p \circ \nu \colon \wt{\ca X} \to \Delta$ is dominant with positive dimensional generic fibre.
\end{enumerate} Then for any general point $0 \in \Delta$, there is an isomorphism of abelian varieties \[J\wt C_{0} \cong (J\wt X_{0})^k,\] and the composition $(J\wt X_{0})^k \xrightarrow{\wt \psi_0} J\wt C_{0} \cong (J\wt X_{0})^k$  
is given by a matrix $\alpha \in \rm{M}_k(\Z)$. Here, $\wt C_0$ and $\wt X_0$ are the normalizations of $C_0$ and $X_0$, and $\wt \psi_0 \colon (J \wt X_0)^k \to J\wt C_0$ is the isogeny induced by $\psi$.  
    \end{proposition}
    
    \begin{proof}
Let $\ca H_{g-1}$ be the moduli stack of smooth hyperelliptic curves of genus $g-1$. 
Let $$f \colon \Delta \longrightarrow \ca H_{g-1}$$ be the map induced by a simultaneous normalization of $p \colon \ca X \to \Delta$, see Proposition \ref{prop:simultaneous-normalization}. 
Then for $0$ in a dense open subset of $\Delta$, the closed analytic subset $$H'_0 \coloneqq f^{-1}(f(0)) \subset  \Delta$$ is positive dimensional by condition \ref{item:prop-constant-normalization:3}. 
For each $u \in H'_0$, the normalization $\wt X_u$ of $X_u$ is isomorphic to $\wt X_0$. 
Moreover, as the map $\Delta \to \overline{\ca M}_g$ is generically finite onto its image (see condition \ref{item:prop-constant-normalization:2}), the same holds for the composition $H'_0 \hookrightarrow \Delta \to \overline{\ca M}_g$. 
In particular, for general $0 \in \Delta$ as above, there exists a connected normal complex analytic space $H_0$ and a dominant generically finite morphism 
\begin{align} \label{align:dominantHH'}
H_0 \longrightarrow H'_0
\end{align}
with the following properties, where $p|_{H_0} \colon \ca X|_{H_0} \to H_0$ denotes the pull-back of $p$ along \eqref{align:dominantHH'}, and $\widetilde{\mathcal X|_{H_0}} \to \ca X|_{H_0}$ the simultaneous normalization of $p|_{H_0}$ (see Proposition \ref{prop:simultaneous-normalization}): 
there is an isomorphism
$$
\widetilde{\mathcal X|_{H_0}}\cong  \widetilde X_0\times H_0
$$
of families of curves over $H_0$, and a non-constant morphism $H_0 \to \widetilde X_0$, $u\mapsto x_u$, such that for any $u\in H_0$, the fibre $X_u$ is obtained from $\widetilde X_u$ by gluing the pair of points $(x_u,\iota(x_u))$, where $\iota$ denotes the hyperelliptic involution on $\widetilde X_0$. 

Let $\widetilde C_{0,i}$ for $i=1,\dots ,n$ be the non-rational irreducible components of the normalization $\widetilde C_0$ of the fibre $C_0=q^{-1}(0)$.  
We apply Proposition \ref{prop:moving} to the family $p|_{H_0} \colon \ca X|_{H_0} \to H_0$. This yields an integer $N \in \Z_{\geq 1}$ such that for each $i \in \set{1, \dotsc, n}$, there is a non-constant morphism
\[
\wt X_0 \longrightarrow
N \cdot \left(
\wt C_{0,i} - \wt C_{0,i}
\right).
\]
By Proposition \ref{proposition:hyperelliptic-gauss}, this implies that the curve $\wt C_{0,i}$ is hyperelliptic for each $i$. Moreover, the assumption that $g = g(X_t) \geq 4$ implies that $g(\wt X_0) \geq 3$ for the genus $g(\wt X_0)$ of the $\wt X_0$. 

By condition \ref{item:prop-constant-normalization:3}, the curve $\wt X_0$ is a very general hyperelliptic curve of genus $g-1 \geq 3$. Thus, by Theorem \ref{proposition:matrixalphai}, the fact that there exists an isogeny
$$
\wt \psi_0 \colon (J\wt X_0)^k
\longrightarrow
J \wt C_{0,1} \times \cdots \times J\wt C_{0,n}$$ implies that for each $i \in \set{1, \dotsc, n}$ we have an isomorphism $\wt C_{0,i} \cong \wt X_0$, and that the composition 
$
(J\wt X_0)^k \to
J \wt C_{0,1} \times \cdots \times J\wt C_{0,n} \cong (J\wt X_0)^k
$
is given by a matrix in $ \rm{M}_k(\Z)$. 
\end{proof}

\subsection{Degeneration of the isogeny: one-dimensional base} \label{subsection:degenerate-one-dim}

Let $D \subset \C$ be a sufficiently small open disc around $0 \in \C$.  
Let $$p \colon \mathcal X\longrightarrow D \quad \text{and} \quad q \colon \ca C \longrightarrow D$$ be families of semi-stable curves over $D$, cf.\ Section \ref{section:conventions}. Let $\psi$ be an isogeny
\begin{align} \label{algin:phiisogeny}
\psi \colon (J\mathcal X)^k \coloneqq 
J\mathcal X\times_{D}\dots \times_{D} J\mathcal X\longrightarrow J\mathcal C
\end{align}
of semi-abelian schemes over $D$, where $J\mathcal X$ and $J\mathcal C$ denote the relative Jacobians over $D$.  
Assume that:
\begin{enumerate}
\item \label{item:degenerate-isogeny-one-dim-base:0} 
for each $s \in D^\ast = D - \set{0}$, the fibre $X_s = p^{-1}(s)$ is a curve of compact type;
\item \label{item:degenerate-isogeny-one-dim-base:3}there exists an isomorphism $J\wt C_{0} \cong (J\wt X_{0})^k$ and a matrix $\alpha \in \rm{M}_k(\Z)$ such that the composition 
$$
    (J\wt X_{0})^k \stackrel{\tilde \psi_0}\longrightarrow J\wt C_{0} \cong (J\wt X_{0})^k
$$ 
    is given by the multiplication by $\alpha$.  
\end{enumerate} 
Recall that, possibly up to shrinking $D$ around $0$, the natural map $H^1(\ca X,\Z) \to H^1(X_0,\Z)$ is an isomorphism (cf.\ Proposition \ref{proposition:fibrewise-retraction}). 
For such sufficiently small $D$, let $t \in D^\ast$, and consider the natural embeddings
    \begin{align} \label{align:injections}
H^1(X_0,\Z) \subset H^1(X_t,\Z)\ \ \ \text{and}\ \ \ H^1(C_t,\Z) \subset H^1(X_t,\Z)^{\oplus k}
\end{align}
induced by the specialization map 
\begin{align}\label{align:def:spX}
sp_{\ca X} \colon H^1(X_0,\Z) \stackrel{\sim}\longleftarrow H^1(\ca X,\Z) \longrightarrow H^1(X_t,\Z)
\end{align} 
and the isogeny $\psi_t \colon (JX_t)^k \to JC_t$. 
 
    \begin{proposition} \label{proposition:congruence-one-dim-base}
    Under the above assumptions, up to possibly shrinking $D$ around $0$ and using the maps in \eqref{align:injections} and \eqref{align:def:spX}, the following identity holds in $H^1(X_t,\Z)^{\oplus k}$  modulo
 $ W_0H^1(X_0,\Z)^{\oplus k} $:
\begin{align} \label{align:degenerate-isogeny:alphacongruence}
H^1(C_t,\Z) \cap H^1(X_{0},\Z)^{\oplus k} \equiv \alpha \cdot H^1(X_{0},\Z)^{\oplus k} \mod  W_0H^1(X_{0},\Z)^{\oplus k},
\end{align}
where
 $ W_0H^1(X_0,\Z) = W_0H^1(X_0,\Q) \cap H^1(X_0,\Z) $ and $\alpha$ is the matrix from \eqref{item:degenerate-isogeny-one-dim-base:3} above.
    \end{proposition}

Before we prove Proposition \ref{proposition:congruence-one-dim-base}, which is the main result of this section, we need the following lemma. To state it, for $s \in D$, let 
\[
\psi_s^\ast \colon H^1(C_s,\Z) \longrightarrow H^1(X_s,\Z)^{\oplus k}
\]
denote the morphism on cohomology induced by the isogeny $\psi_s \colon (JX_s)^k \to JC_s$. 
\begin{lemma} \label{lemma:reduction2:one-dim} 
Up to possibly shrinking $D$ around $0$, we have the following equality:
\begin{align} \label{equation:invariant-image-proofmaintheorem}
\begin{split} 
\Ima&\left(H^1(C_{t},\Z) \stackrel{\psi_t^\ast}\longrightarrow H^1(X_t,\Z)^{\oplus k} \right) \cap \Ima\left( H^1(X_{0},\Z)^{\oplus k} \stackrel{sp_{\ca X}^{\oplus k}}\longrightarrow H^1(X_{t},\Z)^{\oplus k}\right)\\
&=\Ima\left( H^1(C_{0}, \Z) \stackrel{\psi_0^\ast}\longrightarrow H^1(X_{0},\Z)^{\oplus k} \xrightarrow{sp_{\ca X}^{\oplus k}} H^1(X_{t},\Z)^{\oplus k} \right).
\end{split}
\end{align}
Here, $sp_{\ca X}$ is the specialization map defined in \eqref{align:def:spX}.
\end{lemma}

Assuming Lemma \ref{lemma:reduction2:one-dim}, we can prove Proposition \ref{proposition:congruence-one-dim-base} as follows. 

\begin{proof}[Proof of Proposition \ref{proposition:congruence-one-dim-base}]
By Proposition \ref{proposition:fibrewise-retraction}, we may shrink $D$ around $0$ so that the natural map $H^1(\ca X,\Z) \to H^1(X_0,\Z)$ is an isomorphism. By assumption \eqref{item:degenerate-isogeny-one-dim-base:3},
there is an isomorphism $H^1(\wt C_0,\Z) \cong H^1(\wt X_0,\Z)^{\oplus k}$ such that the following diagram commutes, and its rows are exact: 
\begin{align*}
\xymatrix{
0 \ar[r] & W_0H^1(C_{0},\Z)\ar@{^{(}->}[d]^{\psi_0^\ast}  \ar[r] & H^1(C_{0},\Z) \ar[r] \ar@{^{(}->}[d]^{\psi_0^\ast} &H^1(\wt C_0,\Z) \cong  H^1(\wt {X}_{0},\Z)^{\oplus k} \ar@{^{(}->}[d]^{\alpha} \ar[r] & 0 \\
0 \ar[r] & W_0H^1(X_{0},\Z)^{\oplus k} \ar[r] & H^1(X_{0},\Z)^{\oplus k} \ar[r] & H^1(\wt {X}_{0},\Z)^{\oplus k} \ar[r] & 0. 
}
\end{align*}
Here, $W_0H^1(C_0,\Z) = W_0H^1(C_0,\Q) \cap H^1(C_0,\Z)$.
By the above commutative diagram with exact rows, we get
$$
\Ima \left(
H^1(C_{0},\Z) \stackrel{\psi_0^\ast}\longrightarrow  H^1(X_{0},\Z)^{\oplus k}  \longrightarrow H^1(\wt X_{0},\Z)^{\oplus k}
    \right) =\alpha \cdot H^1(\wt X_{0},\Z)^{\oplus k}. 
$$
Consequently, we have 
\begin{align} \label{equation:final-maintheorem}
\Ima\left(
H^1(C_{0},\Z) \stackrel{\psi_0^\ast}\longrightarrow H^1(X_{0},\Z)^{\oplus k}
\right) \equiv \alpha \cdot H^1(X_{0},\Z)^{\oplus k} \mod W_0H^1(X_{0},\Z)^{\oplus k}. 
\end{align}
Combining \eqref{equation:invariant-image-proofmaintheorem} with \eqref{equation:final-maintheorem}, we obtain: 
\begin{align*}
\Ima&\left(H^1(C_{t},\Z) \stackrel{\psi_t^\ast}\longrightarrow H^1(X_t,\Z)^{\oplus k} \right) \cap \Ima\left( H^1(X_{0},\Z)^{\oplus k} \stackrel{sp_{\ca X}^{\oplus k}}\longrightarrow H^1(X_{t},\Z)^{\oplus k}\right) \\
& = \Ima\left( H^1(C_{0}, \Z) \stackrel{\psi_0^\ast}\longrightarrow H^1(X_{0},\Z)^{\oplus k} \stackrel{sp_{\ca X}^{\oplus k}}\longrightarrow  H^1(X_{t},\Z)^{\oplus k} \right) \\
& \equiv \alpha \cdot \Ima\left( H^1(X_{0},\Z)\stackrel{sp_{\ca X}}\longrightarrow H^1(X_{t},\Z)\right)^{\oplus k}\quad \mod \quad  W_0H^1(X_{0},\Z)^{\oplus k}. 
\end{align*}
Thus, \eqref{align:degenerate-isogeny:alphacongruence} holds, and hence the proposition is proved. 
\end{proof}
    
    It remains to prove Lemma \ref{lemma:reduction2:one-dim}. 

\begin{proof}[Proof of Lemma \ref{lemma:reduction2:one-dim}]
By the existence of the isogeny \eqref{algin:phiisogeny}, the restriction 
$
J\ca C|_{D^\ast} \to D^\ast
$
is an abelian scheme over $D^\ast$.  
In particular, 
$
q|_{D^\ast} \colon \ca C|_{D^\ast} \to D^\ast
$
is a family of compact type curves. Thus, 
$
R^1(q|_{D^\ast})_\ast\Z
$ is a local system on $D^\ast$, and in fact a sub-local system of $(R^1(p|_{D^\ast})\Z)^{\oplus k}$. Let
$
S \in \Aut(H^1(C_t,\Z))$ and $ T \in \Aut(H^1(X_t,\Z))$ be generators of the monodromy groups attached to $
R^1(q|_{D^\ast})_\ast\Z
$ and $
R^1(p|_{D^\ast})_\ast\Z$, such that 
\begin{align*} 
    T^{\oplus k}|_{H^1(C_t,\Z)} = S \in \Aut(H^1(C_t,\Z))
\end{align*}
with respect to the inclusion $
H^1(C_t,\Z) \subset H^1(X_t,\Z)^{\oplus k}$ given by $\psi_t^\ast$. Then, we have: 
\begin{align}\label{equation:monodromy-compatibility-maintheorem}
\begin{split}
\Ima&\left(H^1(C_{t},\Z)\stackrel{\psi_t^\ast}\longrightarrow  H^1(X_{t},\Z)^{\oplus k}\right) \cap \left(H^1(X_{t},\Z)^{\oplus k}\right)^{T^{\oplus k}}\\
& = 
\Ima\left(
H^1(C_{t},\Z)^{S}
\longrightarrow H^1(C_t,\Z) \stackrel{\psi_t^\ast}\longrightarrow 
H^1(X_{t},\Z)^{\oplus k}\right). 
\end{split}
\end{align} 
By Proposition \ref{proposition:fibrewise-retraction}, we may shrink $D$ around $0$ so that pulling back along the inclusions $X_{0} \hookrightarrow \ca X$ and $C_0 \hookrightarrow \ca C$ yields isomorphisms $H^1(\ca X,\Z) \cong H^1(X_0,\Z)$ and $H^1(\ca C,\Z) \cong H^1(C_0,\Z)$. 
Define $sp_{\ca X} \colon H^1(X_0,\Z) \to H^1(\ca X,\Z)$ as in \eqref{align:def:spX}, and define in a similar way
$$
sp_{\ca C} \colon H^1(C_0,\Z) \longrightarrow H^1(C_t,\Z).
$$ 
Thus, $sp_{\ca C}$ is the composition of the inverse of the restriction map $H^1(\ca C,\Z) \xrightarrow{\sim} H^1(C_0,\Z)$ with the restriction map $H^1(\ca C,\Z) \to H^1(C_t,\Z)$. 
Then Lemma \ref{lemma:basechange-surjective-monodromy} applied to $p \colon \ca X \to D$ and to $q\colon \ca C \to D$ implies that, possibly after further shrinking $D$ around $0$, we have:
\begin{align}\label{equation:monodromydegeneration-cohomology:final}
\begin{split}
\Ima\left( sp_{\ca X}^{\oplus k} \colon H^1(X_0,\Z)^{\oplus k} \to H^1(X_{t},\Z)^{\oplus k}\right) &= (H^1(X_{t},\Z)^{\oplus k})^{T^{\oplus k}} \subset H^1(X_t,\Z)^{\oplus k}, \\
\Ima\left( sp_{\ca C} \colon  H^1(C_0,\Z) \to H^1(C_{t},\Z)\right) &= H^1(C_{t},\Z)^{S} \subset H^1(C_t,\Z).
\end{split}
\end{align}

\begin{claim} \label{claim:leray}
We have $\psi_t^\ast \circ sp_{\ca C} = sp_{\ca X}^{\oplus k} \circ \psi_0^\ast$ as maps $H^1(C_0,\Z) \to H^1(X_t,\Z)^{\oplus k}$.
\end{claim}

\begin{proof}[Proof of Claim \ref{claim:leray}]
The Leray spectral sequence provides canonical morphisms $H^1(\ca X,\Z)  \to H^0(D, R^1p_\ast\Z)$ and $H^1(\ca C,\Z)
\to H^0(D, R^1q_\ast\Z)$. These make the following diagram commute:
\begin{align*} 
\begin{split}
\xymatrix{
H^1(C_{0},\Z)\ar[d]^-{\psi_0^\ast} 
& H^1(\ca C,\Z) \ar[l]_-{\sim}
\ar[r]
&H^0(D,R^1g_\ast\Z) \ar[d]\ar[r]
&H^1(C_t,\Z) \ar[d]^-{\psi_t^\ast} \\
H^1(X_{0},\Z)^{\oplus k} & H^1(\ca X,\Z)^{\oplus k} \ar[l]_-{\sim} \ar[r]&H^0(D,R^1f_\ast\Z)^{\oplus k}\ar[r]  & H^1(X_t,\Z)^{\oplus k}. 
}
\end{split}
\end{align*} 
As the specialization maps $sp_{ \ca C}$ and $sp_{\ca X}^{\oplus k}$ are obtained by following the horizontal arrows in this diagram from left to right, it follows that $\psi_t^\ast \circ sp_{\ca C} = sp_{\ca X}^{\oplus k} \circ \psi_0^\ast$ as desired.  
\end{proof}
We can finish the proof of Lemma \ref{lemma:reduction2:one-dim}. It suffices to prove the following sequence of equalities: 
\begin{align*} 
\begin{split} 
\Ima&\left(H^1(C_{t},\Z) \stackrel{\psi_t^\ast}\longrightarrow H^1(X_t,\Z)^{\oplus k} \right) \cap \Ima\left( H^1(X_{0},\Z)^{\oplus k} \stackrel{sp_{\ca X}^{\oplus k}}\longrightarrow H^1(X_{t},\Z)^{\oplus k}\right)\\
& = \Ima\left(H^1(C_{t},\Z) \stackrel{\psi_t^\ast}\longrightarrow  H^1(X_t,\Z)^{\oplus k} \right) \cap \Ima\left( \left(H^1(X_{t},\Z)^{\oplus k}\right)^{T^{\oplus k}} \longrightarrow H^1(X_{t},\Z)^{\oplus k}\right) \\
& = 
\Ima\left(
H^1(C_t,\Z)^{S} \longrightarrow H^1(C_t,\Z) \stackrel{\psi_t^\ast}\longrightarrow H^1(X_t,\Z)^{\oplus k}
\right) \\
& = 
\Ima \left(
H^1(C_{0},\Z) \stackrel{sp_{\ca C}}\longrightarrow H^1(C_t,\Z) \stackrel{\psi_t^\ast}\longrightarrow H^1(X_t,\Z)^{\oplus k}
\right) \\
&=\Ima\left( H^1(C_{0}, \Z) \stackrel{\psi_0^\ast}\longrightarrow H^1(X_{0},\Z)^{\oplus k} \stackrel{sp_{\ca X}^{\oplus k}}\longrightarrow H^1(X_{t},\Z)^{\oplus k} \right).
\end{split}
\end{align*}
The first equality follows from \eqref{equation:monodromydegeneration-cohomology:final}, the second equality from \eqref{equation:monodromy-compatibility-maintheorem}, the third equality from \eqref{equation:monodromydegeneration-cohomology:final} again, and the last equality from Claim \ref{claim:leray}. 
\end{proof}

\subsection{Symplectic bases adapted to paths}

Let $U$ be a complex analytic space and let $t \in U$ be a sufficiently general point. Let $n$ and $g$ be positive integers with $n \leq g$, and for $i \in \set{1, \dotsc, n}$, let $\Delta_i \subset U$ be an effective divisor. For each $i$, consider a sufficiently small disc $D_i \subset U$ that intersects the divisor $\Delta \coloneqq  \cup_i\Delta_i$ transversally in a general point $0_i \in \Delta_i$. 

Let $p \colon \ca X \to U$ be a family of nodal curves of arithmetic genus $g$ over $U$. Assume that $p$ is smooth over $U - \Delta$. For each $i$, shrink the disc $D_i$ around $0_i$ so that the natural map 
\[
H^1(\ca X|_{D_i},\Z) \longrightarrow H^1(X_{0_i},\Z)
\]
is an isomorphism (see Proposition \ref{proposition:fibrewise-retraction}). Consider a point $t_i \in D_i - \set{0_i}$ sufficiently close to $0_i$, and let 
\begin{align} \label{path}
    \rho_i \colon [0,1] \longrightarrow U - \Delta
\end{align} be a path from $t$ to $t_i$. The path $\rho_i$ together with a single counter-clockwise loop on $D_i$ induces a loop on $U - \Delta$ and we let $T_i \in \Aut(H^1(X_t,\Z))$ denote the associated monodromy operator. Moreover, the path $\rho_i$ induces a canonical isomorphism
$
H^1(X_{t_i},\Z) \xrightarrow{\sim} H^1(X_t,\Z)$, and we let $sp_{\ca X}^i \colon H^1(X_{0_i},\Z) \to H^1(X_t,\Z)$ denote the composition
\begin{align}\label{align:specializationmap-composition}
sp_{\ca X}^i \colon 
H^1(X_{0_i},\Z) \stackrel{\sim}\longleftarrow H^1(\ca X|_{D_i},\Z) \longrightarrow H^1(X_{t_i},\Z) \stackrel{\sim}\longrightarrow H^1(X_t,\Z). 
\end{align}
Define $W_0H^1(X_{0_i},\Z) = W_0H^1(X_{0_i},\Q) \cap H^1(X_{0_i},\Z)$.  

\begin{definition} \label{definition:monodromyadapted} 
Consider the above notation. 
We say that a symplectic basis
\begin{align}\label{align:definition:adapted-symplectic}
\set{ \delta_1, \dotsc, \delta_g; \gamma_1, \dotsc, \gamma_g } \subset  H^1(X_t,\Z)
\end{align} 
is \emph{adapted to the paths} $\rho_1, \dotsc, \rho_n \colon [0,1] \to U-\Delta$ defined in \eqref{path} if, for each $i \in \{1,\dots ,n\}$, we have:
\begin{align*}
\Ima\left(
H^1(X_{0_i},\Z) \stackrel{sp_{\ca X}^i}\longrightarrow H^1(X_t,\Z)
\right)  &= H^1(X_t,\Z)^{T_i} = \langle \delta_1, \dotsc, \delta_g; \gamma_1, \dotsc, \gamma_{i-1}, \widehat\gamma_i, \gamma_{i+1}, \dotsc, \gamma_g \rangle, \\
\Ima\left(W_0H^1(X_{0_i},\Z) \stackrel{sp_{\ca X}^i}\longrightarrow H^1(X_t,\Z)\right)
&= 
 \Z \cdot \delta_i = \left\langle \delta_i   \right \rangle. 
\end{align*}
\end{definition}

Continue with the above notation. Let $U'$ be a complex analytic space with a surjective generically finite morphism 
$$
\pi \colon U' \longrightarrow U
$$ 
of complex analytic spaces.
Let $t', 0_i' \in U'$ be points on $U'$ $(i = 1, \dotsc, n)$ with $\pi(t') = t$ and $\pi(0_i') = 0_i$ for $i = 1, \dotsc, n$. 
Let $\Delta_i' = \pi^{-1}(\Delta_i)$ 
and $\Delta' = \cup_i \Delta_i' = \pi^{-1}(\Delta)$. 

As $0_i \in \Delta_i$ is a general point, the map $\pi$ looks analytically locally at $0_i$ like the product of a ramified cover of a disc with the identity on a ball of dimension $\dim(U)-1$.

Let $D_i' \subset U'$ be the unique connected component of $\pi^{-1}(D_i)$ that contains $0_i'$.
Up to shrinking $D_i$, $D_i'$ is a disc and $D_i' \to D_i$ is a finite cover of discs, totally ramified at $0_i'$ and \'etale outside of $0_i'$. 
Let $t_i'\in D_i'$ be a lift of $t_i\in D_i$. 
Note that $t_i'$ is automatically in the \'etale locus of $\pi$.
Since $t\in U$ is general, we may assume $t$ is in the \'etale locus of $\pi$ as well.
Up to a small deformation of the path $\rho_i$ 
which does not change its homotopy class nor its beginning and endpoints, we can assume that $\rho_i$ lies also in the \'etale locus of $\pi$.
Under these assumptions, there is a unique path  $$\rho_i' \colon [0,1] \longrightarrow U' - \Delta' \quad \quad (i = 1, \dotsc, n)$$ from $t'$ to $t_i'$ that lifts $\rho_i$.
We then let $T_i' \in \Aut(H^1(X'_{t'},\Z)$ be the monodromy operator induced by the path $\rho_i'$ and the pointed disc $(D_i', 0_i')$. 
 
 \begin{lemma} \label{lemma:adapted-monodromy-basechange}
 Consider the above notation and assumptions. Let $\set{\delta_1, \dotsc, \delta_g;\gamma_1, \dotsc, \gamma_g} \subset  H^1(X_t,\Z)$ be a symplectic basis adapted to the paths $\rho_1, \dotsc, \rho_n$, see Definition \ref{definition:monodromyadapted}. Then the image of $\set{\delta_1, \dotsc, \delta_g;\gamma_1, \dotsc, \gamma_g}$ under the canonical isomorphism
\[
H^1(X_t,\Z) \cong H^1(X'_{t'},\Z),
\]
is a symplectic basis of $H^1(X'_{t'},\Z)$ adapted to the paths $\rho_1', \dotsc, \rho_n'$.  
\end{lemma}

\begin{proof}
Notice that $D_i' \to D_i$ is a finite cover of discs, of the form $z \mapsto z^{m_i}$ for some $m_i \in \Z_{\geq 1}$. 
If we identify $H^1(X'_{t'},\Z)$ with $ H^1(X_t,\Z)$, then the monodromy operator $T_i'$ satisfies $T_i' = T_i^{m_i}$, as elements of $\Aut(H^1(X_t,\Z))$. 
The first thing to show is that
\begin{align}\label{align:invariant-Timi}
 H^1(X_t,\Z)^{T_i^{p}} = \langle \delta_1, \dotsc, \delta_g; \gamma_1, \dotsc, \gamma_{i-1}, \widehat\gamma_i, \gamma_{i+1}, \dotsc, \gamma_g \rangle
\end{align}
for $p = m_i$, knowing that it holds for $p = 1$. To prove this, let $x \in H^1(X_t,\Z)$, and write 
\[
x = a \cdot \gamma_i + y, \quad a \in \Z, \quad y \in  \langle \delta_1, \dotsc, \delta_g; \gamma_1, \dotsc, \gamma_{i-1}, \widehat\gamma_i, \gamma_{i+1}, \dotsc, \gamma_g \rangle = H^1(X_t,\Z)^{T_i}. 
\]
We must show that $T_i^{m_i}(x) = x$ if and only if $a = 0$. 
In other words, we must prove that $T_i^{p}(\gamma_i) \neq \gamma_i$
for $p = m_i$, knowing that it holds for $p = 1$. But this is clear: if $T_i^{m_i}(\gamma_i) = \gamma_i$, then $T_i^{m_i}$ acts trivially on $H^1(X_t,\Z)$, which is absurd since $T_i$ does not act trivially on $H^1(X_t,\Z)$. We conclude that, for each $i$, \eqref{align:invariant-Timi} holds for $p = m_i$. 

Finally, if $f \colon H^1(X_t,\Z) \xrightarrow{\sim}  H^1(X'_{t'},\Z)$ is the isomorphism induced by the canonical isomorphism $X_{t'}'\cong X_t$, then $f$ identifies the images of $W_0H^1(X_{0_i},\Z)$ and $W_0H^1(X'_{0_i'},\Z)$. 
\end{proof}

\subsection{Extending the isogeny}
To prove Theorem \ref{theorem:maintheorem1:hyperelliptic}, we would like to apply Proposition \ref{proposition:degenerate-isogeny:moving-extension}. To do so, we will need the following lemma.  

Let $n$ and $g$ be integers with $g\geq 2$ and $1 \leq n \leq g$. 
Consider a normal algebraic variety $U$ of dimension $2g-1$, irreducible divisors $\Delta_i \subset U$ for $i \in \set{1, \dotsc, n}$, and a family of stable genus $g$ hyperelliptic curves \begin{align} \label{align:family-X-U-semi-stable}
p \colon \ca X \longrightarrow U
\end{align} satisfying the conditions of Lemma \ref{lemma:hyperelliptic-degeneration}.

\begin{lemma} \label{lemma:extending-isogeny:existence:3rd}
 Consider the above notation and let $t\in U$ be a very general point.  
Assume that, for the  fibre $X_t$ of \eqref{align:family-X-U-semi-stable}, there is a smooth projective curve $C$ and an isogeny $\varphi \colon (JX_t)^k \to JC$.   
Then up to replacing $U$ by a normal variety with surjective generically finite morphism $\pi \colon U' \to U$, $p$ by its pull-back along $\pi$, $t$ by a point in $U'$ that lies over it, and $\Delta_i$ by an irreducible divisor in $U'$ that dominates it, there exists a family of stable curves $$q \colon \ca C \longrightarrow U \quad \quad \text{ such that } \quad \quad C_t = q^{-1}(t) = C,$$ together with an isogeny 
    \begin{align} \label{isogeny:semi-abelian:U:2nd} 
\psi \colon (J\mathcal X)^k \coloneqq 
J\mathcal X\times_{U}\dots \times_{U} J\mathcal X\longrightarrow J\mathcal C
\end{align}
of semi-abelian schemes over $U$ that extends the given isogeny $\varphi \colon (JX_t)^k \to JC$, and such that the following conditions are satisfied:
\begin{enumerate}[label=(\roman*)] 
    \item \label{item:condition(ii)}
     For each $i \in \set{1,\dotsc, n}$, the family $p \colon \ca X|_{\wt \Delta_i} \to \wt \Delta_i$ obtained by pulling back $p \colon \ca X \to U$ along the normalization $\wt \Delta_i \to \Delta_i \subset U$ satisfies conditions \ref{item:prop-constant-normalization:1}--\ref{item:prop-constant-normalization:3} in Proposition \ref{proposition:degenerate-isogeny:moving-extension}. 
    \item \label{item:condition(iii)} 
For general $0_i \in \Delta_i$ $(i \in \set{1, \dotsc, n})$, there is a disc $D_i \subset U$ that intersects $\Delta$ transversally in $0_i \in \Delta_i$, so that for general $t_i \in D_i - \set{0_i}$, there is a path $\rho_i \colon [0,1] \to U - \Delta$ from $t$ to $t_i$ with the following property.
There is a symplectic basis $\set{ \delta_1, \dotsc, \delta_g; \gamma_1, \dotsc, \gamma_g} \subset  H^1(X_t,\Z)$ which is adapted to the paths $\rho_1, \dotsc, \rho_n$, see Definition \ref{definition:monodromyadapted}. 
\end{enumerate}
\end{lemma}

\begin{proof} 
We need to spread out the curve $C$ and the isogeny $\varphi \colon (JX_t)^k \to JC$, and our plan is to do this after a suitable base change $U' \to U$. 

\begin{claim} \label{claim:claimm1}
    There is a normal variety $U'$, a generically finite surjective map $\pi \colon U'  \to U$  
    and a family of stable curves \begin{align} \label{align:construct-family-q'}q' \colon \ca C' \longrightarrow U'\end{align} 
    such that for a point $t' \in U'$ lying over $t \in U$, the fibre $C'_{t'} = (q')^{-1}(t')$ is isomorphic to $C$ and the following holds.
    If $p' \colon \ca X' \to U'$ is the pull-back of the family \eqref{align:family-X-U-semi-stable} along $\pi$, then 
   there is a dense open subset $V'\subset U'$ such that the family of curves \eqref{align:construct-family-q'} is smooth over $V'$, and the isogeny $\varphi \colon (JX)^k \to JC$ extends to an isogeny of abelian schemes 
\begin{align} \label{align:isogeny-abelian-V'}
(J\ca X')^k|_{V'} \longrightarrow J\ca C'|_{V'}.
\end{align} 
\end{claim}
\begin{proof}[Proof of Claim \ref{claim:claimm1}] 
This follows from standard spreading out arguments and the properness of the stack $\overline{\ca M}_{kg}$ of stable genus $kg$ curves.  
\end{proof}

\noindent
\begin{claim} \label{claim:claimm2}
Let $U'$, $q' \colon \ca C' \to U'$ and $p' \colon \ca X' \to U'$ be as in Claim \ref{claim:claimm1}. The isogeny \eqref{align:isogeny-abelian-V'}, which is an isogeny of abelian schemes over the open subset $V' \subset U'$, extends to an isogeny 
\[
\psi' \colon (J\mathcal X')^k = 
J\mathcal X'\times_{U'}\dots \times_{U'} J\mathcal X'\longrightarrow J\mathcal C'
\]
of semi-abelian schemes over $U'$. 
\end{claim}
\begin{proof}[Proof of Claim \ref{claim:claimm2}]
As $U'$ is normal, this follows from \cite[Chapter I, Proposition 2.7]{faltings-chai}. 
\end{proof}
Finally, to finish the proof of Lemma \ref{lemma:extending-isogeny:existence:3rd}, it remains to prove that properties \ref{item:condition(ii)} and \ref{item:condition(iii)} hold. As for property \ref{item:condition(ii)}, this holds by Lemma \ref{lemma:hyperelliptic-degeneration} and by the fact that it is stable under base change. Property \ref{item:condition(iii)} follows from Lemmas \ref{lemma:hyperelliptic-degeneration}, \ref{lemma:conditionssatisfied} and \ref{lemma:adapted-monodromy-basechange}. 
\end{proof}

\subsection{Proof of the main theorem}

We are now in position to prove our main theorem. 

\begin{proof}[Proof of Theorem \ref{theorem:maintheorem1:hyperelliptic}] To prove the theorem, we begin with the following reduction step. 
\begin{claim} \label{claim:claim1}
    Theorem \ref{theorem:maintheorem1:hyperelliptic} is implied by the following statement:
    
    \vspace{2mm}
 \emph{   \hypertarget{star}{$(\ast)$} Let $k$ be a positive integer. 
 If, for a very general hyperelliptic curve $X$ of genus $g \geq 4$, there exists a smooth projective curve $C$ and an isogeny $(JX)^k \to JC$, then $k = 1$ and $C \cong X$.}
\end{claim}

\begin{proof}[Proof of Claim \ref{claim:claim1}]
Let $g\geq 4$ and let  $Z\subset \mathcal M_g$ be an irreducible closed subvariety that contains the hyperelliptic locus.
Let $[X]\in Z$ be a very general point, corresponding to a smooth curve $X$ of genus $g$. As mentioned in the introduction, because $JX$ is simple, Theorem \ref{theorem:maintheorem1:hyperelliptic} readily reduces to the case $n=1$: there is an isogeny $(JX)^k \to JC$ between $(JX)^k$ and the Jacobian $JC$ of a smooth projective connected curve $C$. We need to show, under the assumption that \hyperlink{star}{$(\ast)$} holds, that $k = 1$ and $C \cong X$. We specialize $X$ to a very general hyperelliptic curve $Y$. 
This yields a specialization of $C$ into a compact type curve $D$, and an isogeny $(JY)^k \to JD = JD_1 \times \cdots \times JD_n$, where the $D_i$ are the non-rational irreducible components of $D$. 
As $JY$ is simple, there is, for each $i$, an integer $k_i \leq k$ and an isogeny $\varphi_i \colon (JY)^{k_i} \to JD_i$. 
Then \hyperlink{star}{$(\ast)$} implies $k_i = 1$ and $D_i \cong Y$ for each $i$, hence $JD \cong (JY)^k$. Lemma \ref{lemma:isomorphismcriterion} implies $JC \cong (JX)^k$, and then Theorem \ref{cor:iso-of-products-of-curves-main} implies $k=1$ and $C \cong X$ as wanted. 
\end{proof}

Our goal is to prove \hyperlink{star}{$(\ast)$}. 
Thus, let $k \geq 1$ be an integer, and assume that for a very general hyperelliptic curve $X$ of genus $g \geq 4$, there exists a smooth projective curve $C$ and an isogeny $\varphi \colon (JX)^k \to JC$. 
We aim to show that $k=1$ and $C\cong X$.
The strategy is to spread out the isogeny $\varphi$ to an isogeny of families.
To this end we apply Lemma \ref{lemma:extending-isogeny:existence:3rd} and we get
 a normal algebraic variety $U$ with irreducible divisors $\Delta_i \subset U$ $(i = 1,2,3,4)$, families of stable curves
\[
    p \colon \ca X \longrightarrow U \quad \text{and} \quad q \colon \ca C \longrightarrow U
    \]
 with smooth general fibres, and an isogeny
\[
\psi \colon (J\ca X)^k \longrightarrow J\ca C
\]
of semi-abelian schemes over $U$, such that for some $t \in U$, we have $X_t = X, C_t = C$ and $\psi$ restricts to the given isogeny $\varphi \colon (JX)^k \to JC$, and such that all the conditions in Lemma \ref{lemma:extending-isogeny:existence:3rd} are satisfied. 
In particular, $p \colon \ca X \to U$ is smooth over the complement $U - \Delta$ of the divisor $\Delta = \cup_i \Delta_i$, and the morphism $U - \Delta \to \ca H_g$ induced by $p$ is dominant. At this point, in order to prove Theorem \ref{theorem:maintheorem1:hyperelliptic}, it suffices to show that $k = 1$ and $C_t \cong X_t$ (see Claim \ref{claim:claim1}).

\begin{claim} \label{claim:claim2}
For $i \in \set{1,2,3,4}$ and general $0_i \in \Delta_i$, there is an isomorphism of abelian varieties
\[
J\wt C_{0_i} \cong (J\wt X_{0_i})^k,
\]
and the composition $(J\wt X_{0_i})^k \xrightarrow{\wt \psi_{0_i}} J\wt C_{0_i} \cong (J\wt X_{0_i})^k$ is given by a matrix \begin{align} \label{align:matrix:alphai}\alpha_i \in \rm{M}_k(\Z).\end{align} Here, the curves $\wt C_{0_i}$ and $\wt X_{0_i}$ are the normalizations of $C_{0_i}$ and $X_{0_i}$ respectively, and the morphism $\wt \psi_{0_i} \colon (J \wt X_{0_i})^k \to J\wt C_{0_i}$ is the isogeny induced by $\psi$.     
\end{claim}

\begin{proof}[Proof of Claim \ref{claim:claim2}]
For $i \in \set{1,\dotsc, 4}$, the family $p \colon \ca X|_{\wt \Delta_i} \to \wt \Delta_i$ obtained by pulling back $p \colon \ca X \to U$ along the normalization $\wt \Delta_i \to \Delta_i \subset U$ of the divisor $\Delta_i$ satisfies conditions \ref{item:prop-constant-normalization:1}--\ref{item:prop-constant-normalization:3} in Proposition \ref{proposition:degenerate-isogeny:moving-extension}, see Lemma \ref{lemma:extending-isogeny:existence:3rd}. Therefore, the claim follows from Proposition \ref{proposition:degenerate-isogeny:moving-extension}. 
\end{proof}
By construction (see Lemma \ref{lemma:extending-isogeny:existence:3rd}), for each $i  \in \set{1,2,3,4}$ and general $0_i \in \Delta_i$, the fibre $X_{0_i}$ is an irreducible one-nodal hyperelliptic curve, and there exists a disc $D_i \subset U$ intersecting $\Delta$ transversally in $0_i \in \Delta_i$, a path $\rho_i$ from $t \in U$ to a point $t_i \in D_i - \set{0_i}$, and a symplectic basis 
\[
\set{\delta_1, \dotsc, \delta_g; \gamma_1, \dotsc, \gamma_g} \subset H^1(X_t,\Z)
\]
which is adapted to the paths $\rho_1, \dotsc, \rho_4$ in the sense of Definition \ref{definition:monodromyadapted}. For $i \in \set{1,2,3,4}$, shrink $D$ around $0$ so that  the path $\rho_i$ induces a well-defined specialization map $$sp_{\ca X}^i \colon H^1(X_{t_i},\Z) \to H^1(X_{t},\Z),$$ see \eqref{align:specializationmap-composition}. Consider the group $H^1(X_{0_i},\Z)$ as a submodule
$
H^1(X_{0_i},\Z) \subset H^1(X_t,\Z)
$ via $sp_{\ca X}^i$ and consider $H^1(C_t,\Z)$ as a submodule
$
H^1(C_t,\Z) \subset H^1(X_t,\Z)^{\oplus k} 
$ of $H^1(X_t,\Z)^{\oplus k}$ via the map $\psi_t^\ast \colon H^1(C_t,\Z) \to H^1(X_t,\Z)^{\oplus k}$ induced by the isogeny $\psi_t \colon (JX_t)^k \to JC_t$. Define 
$$W_0H^1(X_{0_i},\Z) = W_0H^1(X_{0_i},\Q) \cap H^1(X_{0_i},\Z).$$ 
Thus, $W_0H^1(X_{0_i},\Z)$ is the integral part of the zeroth piece of the weight filtration on $H^1(X_{0_i},\Q)$. 

\begin{claim} \label{claim:claim3}
 For each $i \in \set{1,2,3,4}$, consider the matrix $\alpha_i \in \rm{M}_k(\Z)$ of Claim \ref{claim:claim2}
 , see equation \eqref{align:matrix:alphai}. We have:
\begin{align} \label{align:degenerate-isogeny:alphacongruence:new-proof}
H^1(C_t,\Z) \cap H^1(X_{0_i},\Z)^{\oplus k} \equiv \alpha_i \cdot H^1(X_{0_i},\Z)^{\oplus k} \mod  W_0H^1(X_{0_i},\Z)^{\oplus k}. 
\end{align}   
\end{claim}
\begin{proof}[Proof of Claim \ref{claim:claim3}] 
By Proposition \ref{proposition:congruence-one-dim-base}, the claim follows from Claim \ref{claim:claim2}. 
\end{proof}
We can finish the proof of Theorem \ref{theorem:maintheorem1:hyperelliptic}. By Claim \ref{claim:claim3}, 
we know that \eqref{align:degenerate-isogeny:alphacongruence:new-proof} holds for the matrices $\alpha_i \in \rm{M}_k(\Z)$ of Claim \ref{claim:claim2}. 
Moreover, by Lemma \ref{lemma:extending-isogeny:existence:3rd}, 
the monodromy operators $T_1, \dotsc, T_4 \in \Aut(H^1(X_t,\Z))$ induced by the paths $\rho_i$ and the discs $D_i$ satisfy the property that 
\begin{align*}
\Ima\left(
H^1(X_{0_i},\Z) \stackrel{sp_{\ca X}^i}\longrightarrow H^1(X_t,\Z)
\right)  &= H^1(X_t,\Z)^{T_i} = \langle \delta_1, \dotsc, \delta_g; \gamma_1, \dotsc, \gamma_{i-1}, \widehat\gamma_i, \gamma_{i+1}, \dotsc, \gamma_g \rangle, \\
\Ima\left(W_0H^1(X_{0_i},\Z) \stackrel{sp_{\ca X}^i}\longrightarrow H^1(X_t,\Z)\right)
&= 
 \Z \cdot \delta_i = \left\langle \delta_i   \right \rangle. 
\end{align*}
Thus, by Lemma \ref{lemma:matricessaturatedsubmodules}, we have $\alpha_i \cdot  H^1(X_{t},\Z)^{\oplus k} = \alpha_j \cdot  H^1(X_{t},\Z)^{\oplus k} \subset H^1(C_t,\Z) \subset H^1(X_t,\Z)^{\oplus k}$. In view of Lemma \ref{lemmma:images-freemodules:linear-algebra}, there exist invertible matrices $\gamma_{ij} \in \GL_k(\Z)$ for $i, j \in \set{1,2,3,4}$, such that $\alpha_j = \alpha_i \gamma_{ij}$ for each $i,j$. Moreover, as $\alpha_j \cdot  H^1(X_{t},\Z)^{\oplus k} \subset H^1(C_t,\Z)$, Lemma \ref{lemma:appendix-k=arbitrary-new} implies
\[
H^1(C_t,\Z) = \alpha_i \cdot  H^1(X_{t},\Z)^{\oplus k} \subset  H^1(X_{t},\Z)^{\oplus k} \quad \quad \forall i = 1,2,3,4. 
\]
Consequently, by Lemma \ref{lemma:isogenyVHS}, there exists an isomorphism of abelian schemes
$$J\ca C \cong (J\ca X)^k$$
over $U$. By Theorem \ref{cor:iso-of-products-of-curves-main}, it follows that $k = 1$ and $C_t \cong X_t$.

We have proven that \hyperlink{star}{$(\ast)$} holds. By Claim \ref{claim:claim1}, we are done.
\end{proof}

\section{Abelian varieties with no power isogenous to a Jacobian} \label{section:maintheorem2:IHCapplication}
The goal of this section is to prove Theorem \ref{thm:cor:maintheorem2} and Corollaries \ref{cor:coleman-oort},  \ref{corollary:IHCconsequence} and \ref{cor:IHCconsequence-2} stated in the introduction.  

\begin{proof}[Proof of Theorem \ref{thm:cor:maintheorem2}] 

First, we deal with the case where $A = J^3Y$ is the intermediate Jacobian of a very general cubic threefold $Y$. We claim that there is no integer $k \geq 1$ for which there exists an isogeny between $(J^3Y)^k$ and a product of Jacobians. To prove this, assume that such an integer and such an isogeny exist. As $Y$ is very general, $J^3Y$ is simple; in particular, we may assume that there exists a curve $C$ and an isogeny $\phi \colon JC \to (J^3Y)^k$. 
Degenerate $Y$ into a singular cubic $Y_0$ such that $J^3Y_0 = JX$ is the Jacobian of a very general hyperelliptic curve $X$, cf.\ \cite{collino-fanofundamentalgroup}. This leads to a degeneration of $C$ into a compact type curve $D$ and an isogeny $\phi_0 \colon JD \to (JX)^k$. By Theorem \ref{theorem:maintheorem1:hyperelliptic}, this implies 
$JD \cong (JX)^k$, hence by Lemma \ref{lemma:isomorphismcriterion}, we get $JC \cong (J^3Y)^k$, which contradicts Corollary \ref{cor:cubicthreefold-isomorphism}. 

It remains to show that if $A$ is a very general principally polarized abelian variety of dimension $g \geq 4$, then there is no integer $k \geq 1$ for which there exists an isogeny between $A^k$ and a product of Jacobians. To arrive at a contradiction, we may assume that, for some $k \in \Z_{\geq 1}$, there exists an isogeny $\phi \colon JC \to A^k$ for some smooth projective curve $C$. Specialize $A$ to the Jacobian $A_0 = JX$ of a very general hyperelliptic curve of genus $g$. The curve $C$ specializes to a compact type curve $D$, hence we obtain an isogeny $\phi_0 \colon JD  \to (JX)^k$. 
By Theorem \ref{theorem:maintheorem1:hyperelliptic}, we have $JD \cong (JX)^k$, 
hence there exists an isomorphism of abelian varieties $JC \cong A^k$ by Lemma \ref{lemma:isomorphismcriterion}. 
This implies by Theorem \ref{cor:iso-of-products-of-curves-main-main} that $A$ is isomorphic as a polarized abelian variety to the Jacobian of a curve, which is absurd by dimension reasons, because $g\geq 4$. The theorem follows.  
\end{proof}

\begin{proof}[Proof of Corollary \ref{cor:coleman-oort}] 
The fact that $Z\subset \mathcal A_g$ is special is  well-known to experts; we include an argument in Lemma \ref{lemma:special-lemma} in Appendix \ref{appendix:B}. 
The fact that $Z\subset \mathcal A_g$ satisfies the Coleman--Oort conjecture follows from the fact that for a very general principally polarized abelian variety $A$ of dimension $g \geq 4$, the $k$-th power $A^k$ is not isogenous to a Jacobian of a curve, see Theorem \ref{thm:cor:maintheorem2}. 
 \end{proof}
  
\begin{proof}[Proof of Corollary \ref{corollary:IHCconsequence}]
Let $A$ be either the intermediate jacobian of a very general cubic threefold or a very general principally polarized abelian variety of dimension at least $4$.
Let $A_1$ be an abelian variety isogenous to a power of $A$ and let $A_2$ be an abelian variety with $\Hom(A,A_2)=0$.

Suppose that there are smooth projective curves $C_1,\dots ,C_n$ and an isomorphism
 \begin{align} \label{equation:isomorphism:maintheorem2}
 A_1 \times A_2 \cong JC_1 \times \cdots \times JC_n .
 \end{align}
 We claim that there exists a non-empty subset $I \subset \set{1,\dotsc, n}$ such that $A_1 \cong \prod_{i \in I}JC_i$. Indeed, the product polarization on $\prod_{i =1}^nJC_i$ and the isomorphism \eqref{equation:isomorphism:maintheorem2} equip $A_1 \times A_2$ with a principal polarization, call it $\lambda$. We have $\NS(A_1 \times A_2) = \NS(A_1) \times \NS(A_2)$ because $\Hom(A_1, A_2) = 0$. Hence $\lambda = \lambda_1 \times \lambda_2$ for principal polarizations $\lambda_i$ on $A_i$. By \cite[Corollary 3.23]{clemensgriffiths-cubicthreefold} (see also \cite{debarre-produits}), the decomposition of a principally polarized abelian variety into a product of principally polarized abelian subvarieties is unique. Therefore, $(A_1, \lambda_1)$ is isomorphic to $\prod_{i =1}^nJC_i$ for some non-empty subset $I \subset \set{1,\dotsc, n}$, proving the claim. 
 
 Since $A_1$ is isogenous to a power of $A$, we find that $\prod_{i =1}^nJC_i$ is isogenous to a power of $A$, which
 contradicts Theorem \ref{thm:cor:maintheorem2}.
 This concludes the proof of Corollary \ref{corollary:IHCconsequence}.
\end{proof}

\begin{proof}[Proof of Corollary \ref{cor:IHCconsequence-2}]
    This is a direct consequence of Corollary \ref{corollary:IHCconsequence}.
\end{proof}

\appendix

\section{Remark on the degeneration method} 

\label{appendix}

In an influential work, Bardelli and Pirola \cite{bardellipirola-curvesofgenusg} proved that for a very general curve $X$ of genus $g\geq 4$, the Jacobian $JX$ is not isogenous to a Jacobian of any smooth curve $C$ with $C\not \cong X$.
Their argument is based on the following idea, which also played an important role in consecutive papers (including this paper). 
If $f \colon JC\to JX$ is an isogeny, then the goal is to show that there is an integer $n$ such that the image of $f^\ast \colon H^1(JX,\Z)\to H^1(JC,\Z)$ satisfies
\begin{align}\label{eq:im(f*)=nH1}
    \im(f^\ast)=n\cdot H^1(JC,\Z)\subset H^1(JC,\Z).
    \end{align}
Indeed, \eqref{eq:im(f*)=nH1} implies the existence of an isomorphism $JX\cong JC$ with respect to which the isogeny $f$ is given by multiplication by $n$. As $X$ is very general, this isomorphism has to respect the polarizations, and thereby $C \cong X$ because of the Torelli theorem for curves. 
 
    To prove \eqref{eq:im(f*)=nH1}, 
    a degeneration argument to nodal curves is used to show that there are suitable monodromy operators $T_1,T_2$ corresponding to two different Picard--Lefschetz degenerations of $X$, such that \eqref{eq:im(f*)=nH1} holds when intersected with the $T_i$-invariant subspaces for $i=1,2$.
    To get \eqref{eq:im(f*)=nH1} from this, it is then used without proof in \cite[Proposition 4.1.3]{bardellipirola-curvesofgenusg} that 
    \begin{align} \label{eq:imf*=imT1+imT2}     \im(f^\ast)=\im(f^\ast)^{T_1}+\im(f^\ast)^{T_2} .
    \end{align}
    The argument is formalized in \cite[Proposition 3.6]{marcucci-genusofcurves}, but also here, the identity \eqref{eq:imf*=imT1+imT2} (which is the identity $H_1+H_2=\mathcal H_x$ in the notation of loc.\ cit.)\ is assumed implicitly in the proof. 
The same identity is used 
in the proof of Theorem \ref{theorem:C-hyperell-g=3} in \cite{naranjopirola2018}; see also Remark \ref{remark:naranjo-pirola-advances} in Section \ref{section:genusgcurveshyperellipticjacobians}.

It turns out that the version in \cite[Proposition 3.6]{marcucci-genusofcurves}
is incorrect, see Proposition \ref{proposition:counterexample} below. 
Nonetheless, it is possible to prove \eqref{eq:imf*=imT1+imT2} under the additional assumption that $f^\ast\Theta_{X}$ is a multiple of the theta divisor $\Theta_{C}$ of $JC$ (which is the case in \cite{bardellipirola-curvesofgenusg, naranjopirola2018}). 
For this, one can use the lattice theoretic results provided in Section \ref{section:linear-algebra}, see in particular Lemmas \ref{lemma:matricessaturatedsubmodules} and \ref{lemma:appendix-k=arbitrary-new}. 
These statements and their proofs greatly simplify if one is only interested in the $k=1$ case.  
For convenience of the reader, we state and prove the precise statement that one needs to prove \eqref{eq:imf*=imT1+imT2}  in the aforementioned applications in \cite{bardellipirola-curvesofgenusg,naranjopirola2018}
in Proposition \ref{proposition:marcucci:improved} below.

\begin{proposition} \label{proposition:marcucci:improved}
Let $H\subset G$ be free $\Z$-modules of the same finite rank and let $E_G$ and $E_H$ be unimodular symplectic forms on $G$ and $H$, respectively. 
Let $\set{\delta_1,\dotsc, \delta_g; \gamma_1,\dotsc, \gamma_g} \subset G$ be a symplectic basis of $G$.
Let $T_i\in \Aut(G)$ for $i=1,2$ be linear automorphisms with $T_i(H)\subset H$ for all $i=1,2$. 
Assume the following conditions. 
\begin{enumerate}[label=(\roman*)]
\item
\label{item:marcucci:improved:A}
For each $i$, there is a positive integer $n_i$ such that $H^{T_i}=n_iG^{T_i}$. 
    \item \label{item:marcucci:improved:B} 
    The $T_i$-invariant subspaces of $G$ are given by 
    $$
G^{T_1}=\langle\delta_1,\dots ,\delta_g,\gamma_2,\dots ,\gamma_g \rangle \quad \quad \text{and} \quad \quad G^{T_2}=\langle\delta_1,\dots ,\delta_g,\gamma_1,\gamma_3,\dots ,\gamma_g \rangle.
    $$
 \item \label{item:marcucci:improved:C} The form $E_G$ on $G$ restricts to a multiple of the form $E_H$ on $H\subset G$.
\end{enumerate}
Then $n:=n_1=n_2$ and $H = nG \subset G$. 
\end{proposition}

Before we prove Proposition \ref{proposition:marcucci:improved}, we consider the following result, due to Marcucci.  

\begin{proposition}[Marcucci]
\label{proposition:marcucci}
Let $H\subset G$ be free $\Z$-modules of the same finite rank. Let $T_i\in \Aut(G)$ for $i=1,2$ be linear automorphisms with $T_i(H)\subset H$ for all $i=1,2$. 
Assume that 
\begin{enumerate}
\item $G^{T_1} + G^{T_2} = G$;
\label{item:marcucci:1}
\item 
\label{item:marcucci:2}
$G^{T_1} \cap G^{T_2} \neq 0$;
\item 
\label{item:marcucci:3}for each $i$, there is a positive integer $n_i$ such that $H^{T_i}=n_iG^{T_i}$. 
\end{enumerate}
Then $n:=n_1 = n_2$ and $n G \subset H \subset G$. 
\end{proposition}
\begin{proof}
    See the proof of \cite[Proposition 3.6]{marcucci-genusofcurves} (or see the more general Lemma \ref{lemma:matricessaturatedsubmodules}). 
\end{proof}

Next, we prove Proposition \ref{proposition:marcucci:improved}. 

\begin{proof}[Proof of Proposition \ref{proposition:marcucci:improved}]
Note that $n_1 = n_2$ and $n_1G \subset H$ by Proposition \ref{proposition:marcucci}. By Lemma \ref{lemma:appendix-k=arbitrary-new}, this implies $nG_1 = H$ as we want. For convenience of the reader, let us sketch the proof in this particular situation. 
Let $n \coloneqq n_1$ and $U \coloneqq \langle ne_1, ne_2 \rangle \subset G$. Then $U$ is a saturated, isotropic subspace of $H$; by Lemma \ref{lem:isotropic-subspace}, there are $g_1, g_2 \in G$ such that $U \oplus \langle g_1, g_2 \rangle \subset H$ is a unimodular sublattice of $H$. By condition \ref{item:marcucci:improved:C}, there exists $m \in \Z_{\geq 1}$ such that $E_H(x,y) = E_G(m^{-1}x,y)$ for $x,y \in H$. Thus, $[G \colon H] = m^g$ (cf.\ Lemma \ref{lemma:appendix:preliminary-lemma}), and $g_2 \in \langle \delta_1, \dotsc, \delta_g; \gamma_2, \dotsc, \gamma_g\rangle = G^{T_1}$. 

In view of condition \ref{item:marcucci:improved:A}, there exists $g_2' \in H_1$ such that $g_2 = n g_2'$. As we have an equality $E_H(ne_2,g_2)=E_G(m^{-1}n e_2, ng_2') = 1$, we must have $n^2 \mid m$. Moreover, the inclusion $nG \subset H$ shows that $ m^g = [G \colon H]  \mid [G \colon nG]= n^{2g}$. Thus, $m =n^2$, and the equality $H = nG$ follows. 
\end{proof}

\subsection{Counterexample}

In this section we show that, in contrast to what is claimed in \cite[Proposition 3.6]{marcucci-genusofcurves}, the equality $nG = H$ does not follow in general from the conditions \eqref{item:marcucci:1}--\eqref{item:marcucci:3} of Proposition \ref{proposition:marcucci}.  
We try to keep our notation as compatible as possible with loc.\ cit.

\begin{notation} \label{notation:marcucci}
    Let $G$ be a free $\Z$-module of rank four, with basis $\set{e_1, e_2, f_1, f_2} \subset G$. We equip $G$ with the symplectic form 
\[
\left(-, - \right) \colon G\times G \longrightarrow \Z 
\]
that has $\set{e_1, e_2, f_1, f_2}$ as symplectic basis. In other words, for $i,j \in \set{1,2}$, we have $$(e_i, f_j) = \delta_{ij}, \quad (e_i, e_j) = 0 = (f_i, f_j), \quad (y,x) = - (x,y) \quad \forall x,y\in G.$$  
Let $k$ and $n$ be positive integers such that \begin{align}\label{align:1<k<n}
k \mid n, \quad 1 < k \leq n.
\end{align}
We define a submodule $H \subset G$ as follows:
\[
H = \left\langle ne_1, ne_2, nf_1, \frac{n}{k}f_2 + \frac{n}{k}f_1 \right\rangle \subset \left\langle e_1, e_2, f_1, f_2 \right\rangle = G. 
\]
\end{notation}
Notice that
\[
nG = \langle ne_1, ne_2, nf_1, nf_2\rangle 
\subset 
\left\langle ne_1, ne_2, nf_1, \frac{n}{k}f_2 + \frac{n}{k}f_1 \right\rangle = H. 
\]
\begin{lemma} \label{lemma:strict-inequality}
    The inequality $nG \subset H$ is strict. More precisely, the quotient $H/nG$ is a finite cyclic group of order $k > 1$. 
\end{lemma}
\begin{proof}
    This follows from \eqref{align:1<k<n}. 
\end{proof}

Consider the free abelian group $\Z^{\oplus 2}$. We are going to define an action of $\Z^{\oplus 2}$ on $G$. For $i = 1,2$, define an operator $T_i \colon G \to G$ by $T_i(x) = x + (x, ke_i)ke_i$. 
Then $T_1$ and $T_2$ are automorphisms of $G$. As $(e_1, e_2) = 0$, we obtain an action of $\Z^{\oplus 2}$ on $G$:
\begin{align} \label{align:groupaction-marcucci}
\Z^{\oplus 2} \longrightarrow \Aut(G), \quad b_i \mapsto T_i,
\end{align}
where $b_1 = (1,0) \in \Z^{\oplus 2}$ and $b_2 = (0,1) \in \Z^{\oplus 2}$. 
\begin{lemma} \label{lemma:preserves}
The action of $\Z^{\oplus 2}$ on $G$ preserves the submodule $H \subset G$. 
\end{lemma}
\begin{proof}
We have $T_1(ne_1) = ne_1$ and $T_1(ne_2) = ne_2$. Moreover, we have $
    T_1(nf_1)= nf_1 - k^2 ne_1$ and $
    T_1\left( \frac{n}{k}f_2 + \frac{n}{k}f_1  \right) =
    \frac{n}{k}f_2 + \frac{n}{k}f_1 - kn e_1$. In a similar way, $T_2(ne_1) = ne_1$, $T_2(ne_2) = ne_2$, $ T_2(nf_1) = nf_1$, and $T_2\left( \frac{n}{k}f_2 + \frac{n}{k}f_1  \right) = 
    \frac{n}{k}f_2 + \frac{n}{k}f_1 - kn e_2$. 
These are all elements of $H$. 
\end{proof}
Next, we would like to calculate $H^{T_i} = H \cap G^{T_i}$ for $i = 1,2$.
\begin{lemma} \label{lemma:invariant-Ti-H}
    We have $H^{T_i} = n \cdot G^{T_i}$ for $i = 1,2$. 
\end{lemma}
\begin{proof}
Notice that 
$
G^{T_i} = \set{x \in G \mid T_i(x) = x} = \set{x \in G \mid (x,e_i) = 0}. 
$
Hence, \begin{align} \label{align:invariantGTi}
G^{T_1} = \langle e_1, e_2, f_2 \rangle \quad \quad \text{and}\quad \quad 
G^{T_2} = \langle e_1, e_2, f_1 \rangle.
\end{align}
Therefore, $H^{T_1} = H \cap G^{T_1} = H \cap \langle e_1, e_2, f_2 \rangle = \left \langle ne_1, ne_2, nf_2 \right \rangle = n G^{T_1}$. Similarly, we have $H^{T_2} = H \cap G^{T_2} = H \cap \langle e_1, e_2, f_1 \rangle= \left \langle ne_1, ne_2, nf_1 \right \rangle = nG^{T_2}$. 
\end{proof}

From the previous results, we deduce the following result, which shows that \cite[Proposition 3.6]{marcucci-genusofcurves} fails in the generality stated. For a corrected version with stronger hypotheses, which seems to suffice for most of the applications, see Proposition \ref{proposition:marcucci:improved}. 

\begin{proposition} \label{proposition:counterexample}
There exists a connected and path-connected pointed topological space $(X,x)$, two local systems of free abelian groups of rank four $\ca H$ and $\ca G$ on $X$, an injective morphism of local systems $\ca H \hookrightarrow \ca G$, and two elements $\gamma_i \in \pi_1(X,x)$ $(i = 1,2)$, such that the following holds. If the monodromy representations attached to $\ca H$ and $\ca G$ are denoted by 
\[
\rho \colon \pi_1(X,x) \longrightarrow \Aut(\ca H_x) \quad \text{and} \quad \sigma \colon \pi_1(X,x) \longrightarrow \Aut(\ca G_x), 
\]
and if $G_i \subset \ca G_x$ and $H_i \subset \ca H_x$ are defined as
\begin{align} \label{align:GiHi}
\begin{split}
G_i \coloneqq \rm{Inv}(\sigma(\gamma_i)) = \set{a \in \ca G_x \mid \sigma(\gamma_i)(a) = a}, \\
H_i \coloneqq \rm{Inv}(\rho(\gamma_i)) = \set{a \in \ca H_x \mid \rho(\gamma_i)(a) = a}, 
\end{split}
\end{align}
then the following conditions are satisfied: 
\begin{enumerate}
    \item \label{item:LS-one} $G_1 + G_2 = \ca G_x$;
    \item \label{item:LS-two} $G_1 \cap G_2 \neq 0$; 
    \item \label{item:LS-three} there exists an integer $n \in \Z_{\geq 1}$ such that $H_i = nG_i$ for $i  =1,2$;
    \item \label{item:LS-four} with respect to the integer $n$ in condition \eqref{item:LS-three}, we have
    \[
    \ca H_x \neq n \ca G_x
    \]
    as submodules of $\ca G_x$. 
\end{enumerate}
\end{proposition}
\begin{proof}
Let $D^\ast  = \set{z \in \C \mid 0 < \va{z} < 1}$, and define $X = D^\ast \times D^\ast$. Let $x \in X$ be any point. Then $\pi_1(X,x) = \Z^{\oplus 2}$. 

Define two free $\Z$-modules of rank four $H \subset G$ as in Notation \ref{notation:marcucci}. Let $\Z^{\oplus 2}$ act on $G$ as in \eqref{align:groupaction-marcucci}. In particular, the action of $\Z^{\oplus 2}$ on $G$ restricts to an action of $\Z^{\oplus 2}$ on $H \subset G$, see Lemma \ref{lemma:preserves}. This yields two representations $\rho \colon \pi_1(X,x) \to \Aut(H)$ and $\sigma \colon \pi_1(X,x) \to \Aut(G)$. 

Let $\ca H$ and $\ca G$ be the local systems on $X$ attached to the representation $\rho \colon \pi_1(X,x) \to \Aut(H)$ and $\sigma \colon \pi_1(X,x) \to \Aut(G)$, respectively. Define $G_i$ and $H_i$ as in \eqref{align:GiHi}. 

We claim that conditions \eqref{item:LS-one}--\eqref{item:LS-four} are verified. Conditions \eqref{item:LS-one} and \eqref{item:LS-two} follow from \eqref{align:invariantGTi}. Condition \eqref{item:LS-three} follows from Lemma \ref{lemma:invariant-Ti-H}. Condition \eqref{item:LS-four} follows from Lemma \ref{lemma:strict-inequality}. 
\end{proof} 

\section{Jacobians isogenous to a power of an elliptic curve} \label{appendix:B} 

In \cite{luzuo-shimuracurves}, Lu and Zuo prove that for a very general elliptic curve $E$, no power $E^g$ with $g \geq 12$ is in the Hecke orbit of the Jacobian of a smooth projective connected curve of genus g (see \cite[Theorem A]{luzuo-shimuracurves} and Remark \ref{remark:be-aware} below). The goal of this appendix is to show that the methods of \cite{luzuo-shimuracurves} in fact imply the following stronger result.

\begin{theorem} \label{theorem:verygeneral-powers-elliptic-isogenies}
For an elliptic curve $E$ with transcendental $j$-invariant, the following holds: \begin{enumerate}
    \item \label{item:1:corollary:luzuo:vg}There exists no integer $g \geq 12$ such that $E^g$ is isogenous to the Jacobian of a smooth projective connected curve. \label{item:thm:powers-of-elliptic:1}
    \item \label{item:2:corollary:luzuo:vg}There exists no integer $g \geq 5$ such that $E^g$ is isogenous to the Jacobian of a smooth projective connected hyperelliptic curve. 
    \end{enumerate}
\end{theorem} 
 
 \begin{remark} \label{remark:be-aware}
At most places in \cite{luzuo-shimuracurves} the term ``isogenous'' means ``to lie in the same Hecke orbit", see \cite[paragraph below Definition 2.12]{luzuo-shimuracurves} and \cite[Lemma 2.13]{luzuo-shimuracurves}. 
If the points in $\ca A_g$ associated to two principally polarized abelian varieties $(A, \lambda_A)$ and $(B,\lambda_B)$ lie in the same Hecke orbit, then $A$ and $B$ are isogenous, but the converse is not necessarily true. In fact, one can show that the following are equivalent:
\begin{enumerate}
    \item The associated moduli points $[A],[B] \in \ca A_g$ lie in the same Hecke orbit, that is, admit lifts $x,y \in \bb H_g$ to the Siegel upper half space $\bb H_g$ that lie in the same $\rm{GSp}_{2g}(\Q)_+$-orbit.
    \item There is an isomorphism of  
    rational Hodge structures $H^1(A,\Q) \cong H^1(B,\Q)$ that preserves the polarizations up to a positive rational multiple.
    \item There is an isogeny $\phi \colon A \to B$ such that $\phi^\ast(\lambda_B) = n \cdot \lambda_A$
    for some $n \in \Z_{\geq 1}$. 
\end{enumerate} 
\end{remark}

\begin{remark} \label{rem:Hecke-orbit-A^k}
One can show that if $(A,\lambda)$ is a very general principally polarized abelian variety of dimension $g$, then for any integer $k \geq 1$ and any principal polarization $\mu$ on the $k$-th power $A^k$ of $A$, the moduli points of $(A^k, \mu)$ and $(A^k, \lambda^k)$ have isomorphic polarized
rational Hodge structures
(where $\lambda^k$ denotes the product polarization on $A^k$ associated to $\lambda$). 
To prove this, notice that by Lemma \ref{lemma:polarizationcorrespondence}, $\mu $ coincides with the polarization $\lambda_\alpha$ associated to a unimodular positive definite symmetric matrix $\alpha \in \GL_k(\Z)$. 
We thus need to show that the principally polarized abelian varieties $(A^k, \lambda_\alpha)$ and $(A^k, \lambda^k)$ have isomorphic polarized
rational Hodge structures. 
By a suitable analogue of Lemma \ref{lemma:polarizationcorrespondence}, that considers rational equivalence classes of polarizations on $A^k$, this comes down to proving that for each matrix $\alpha \in \GL_k(\Z)$ as above there exists a matrix $\gamma \in \GL_k(\Q)$ such that 
$\gamma \alpha \gamma^t$ is the identity matrix, where $\gamma^t$ denotes the transpose of $\gamma$. 
This turns out to be true, and can be deduced from the results in \cite[Chapter V, Sections 1.3.5 and 1.3.6]{serre-course-in-arithmetic}.
\end{remark}

\begin{remark}
We used Theorem \ref{theorem:verygeneral-powers-elliptic-isogenies} in the proof of Theorem \ref{cor:iso-of-products-of-curves-main-main}, which in turn is used in the proofs of Theorems \ref{theorem:maintheorem1:hyperelliptic} and \ref{thm:cor:maintheorem2}.
In fact, for these applications, one only needs the special case of item \eqref{item:1:corollary:luzuo:vg} in Theorem \ref{theorem:verygeneral-powers-elliptic-isogenies} in which the isogeny $E^g\to JC$ is an isomorphism of unpolarized abelian varieties, and this special case could alternatively be deduced directly from Remark \ref{rem:Hecke-orbit-A^k} and \cite[Theorem A]{luzuo-shimuracurves}. 
We decided to give the proof of Theorem \ref{theorem:verygeneral-powers-elliptic-isogenies} (instead of spelling out the details of the result alluded to in Remark \ref{rem:Hecke-orbit-A^k}), since Theorem \ref{theorem:verygeneral-powers-elliptic-isogenies} seems of independent interest, and naturally complements Theorem \ref{thm:cor:maintheorem2}.
\end{remark}

\subsection{Special subvarieties in moduli spaces of abelian varieties}\label{section:specialsub} 
We need to gather some results on special subvarieties in $\ca A_g$. Let us start by recalling the definition. 
For integers $g \geq 1$ and $n \geq 1$, let $\ca A_{g,[n]}$ be the moduli space of principally polarized abelian varieties of dimension $g$ with level $n$ structure (which is a scheme if $n \geq 3$ and an algebraic stack in general).  
A closed subvariety $Z \subset \ca A_{g,[n]}$ is called a \emph{special subvariety} if it is a Hodge locus of the $\Q$-variation of Hodge structure $R^1\mf h_\ast\Q$, where $\mf h \colon \ca X_{g,[n]} \to \ca A_{g,[n]}$ denotes the universal family, see \cite[Definition 3.7]{moonenoort-specialsubvarieties}.
Special subvarieties of $\ca A_{g,\delta,[n]}$ for some polarization type $\delta$ are defined similarly. 
By way of example, we have:
\begin{lemma} \label{lemma:special-lemma}
Let $n \geq 3$ be an integer and let $h, k, g$ be positive integers with $g = hk$. Let $Z \subset \ca A_{g,[n]}$ be a subvariety such that the general element of $Z$ is isogenous to the $k$-th power of a general polarized abelian variety of dimension $h$. Then $Z \subset \ca A_{g,[n]}$ is a special subvariety of dimension $h(h+1)/2$.
\end{lemma}
\begin{proof} 
By \cite[Remark 3.13]{moonenoort-specialsubvarieties}, this easily reduces to the case where $Z \subset \ca A_{g,[n]}$ is the image of the diagonal embedding $\ca A_{h,[n]} \hookrightarrow \ca A_{hk,[n]} = \ca A_{g,[n]}$, in which case the result is clear.  
\end{proof}

We are grateful to Kang Zuo for pointing us to item \eqref{item:smoothness-special} of the following lemma.

\begin{lemma} \label{lemma:moonen-deligne}
    Let $g \geq 1, n \geq 3$ be integers. Let $Z \subset \ca A_{g,[n]}$ be a special subvariety. Then the following assertions are true. 
    \begin{enumerate}
        \item \label{item:smoothness-special}There exists an integer $m \geq 3$ with $n \mid m$ and a smooth special subvariety $Y \subset \ca A_{g,[m]}$ such that $Y$ is an irreducible component of the preimage of $Z$ under the natural finite \'etale map $\ca A_{g,[m]} \to \ca A_{g,[n]}$. 
        \item \label{item:unipotent-monodromy}If $Z \subset \ca A_{g,[n]}$ is one-dimensional, then there are $m \geq 3$ and $Y \subset \ca A_{g,[m]}$ as in item \eqref{item:smoothness-special} with the additional property that the pull-back of the universal weight one $\Q$-local system on $\ca A_{g,[m]}$ to the smooth one-dimensional subvariety $Y \subset \ca A_{g,[m]}$ has unipotent local monodromy around each point of $\ol{Y} - Y$, where $\ol{Y}$ is the smooth projective model of $Y$.  
    \end{enumerate}
\end{lemma}
\begin{proof}
This is well-known; item \eqref{item:smoothness-special} follows e.g.\ from \cite[Lemma 3.3]{moonen-linearity-I} and item \eqref{item:unipotent-monodromy} from item \eqref{item:smoothness-special} together with \cite[Proposition 1.6(2)]{yau-zhang-toroidal}.  
\end{proof}

\subsection{Arakelov inequality for a family of abelian varieties over a curve}

To prove Theorem \ref{theorem:verygeneral-powers-elliptic-isogenies}, the idea is to apply the following result, due to Lu and Zuo \cite{luzuo-shimuracurves}. To state it, we need to introduce the following notation. Let $\ol C$ be a smooth projective connected curve, and let $C \subset \ol C$ be an open subscheme.  
Let $h \colon A \to C$ be a family of abelian varieties. 
Define $\Delta_{\ol C} = \ol C - C$ and assume that the local monodromy of $R^1h_\ast \Q$ around each point of $\Delta_{\ol C}$ is unipotent. Consider the Deligne extension $(R^1h_\ast \Q \otimes_{\Q} \OO_C)_{\rm{ext}}$ of the vector bundle $R^1h_\ast \Q \otimes_{\Q} \OO_C$, see \cite[Proposition 5.2, page 91]{MR417174}, 
which is a vector bundle on $\ol C$ that extends $R^1h_\ast \Q \otimes_{\Q} \OO_C$. The Hodge filtration
\[
0 \subset E^{1,0} \subset 
R^1h_\ast \Q \otimes_{\Q} \OO_C \quad \quad (E^{1,0} = h_\ast \Omega^1_{A/C})
\]
extends, in view of the nilpotent orbit theorem (see \cite[Theorem 2.1]{cattani-kaplan-degenerating}), to a filtration
\begin{align} \label{align:filtration-deligne-extension}
0 \subset E_{\ol C}^{1,0} \subset (R^1h_\ast \Q \otimes_{\Q} \OO_C)_{\rm{ext}}.
\end{align}
    By \cite{faltings-arakelov}, we have that 
\begin{align}\label{align:faltings:prelim}
\deg(E^{1,0}_{\ol C}) \leq \frac{g}{2} \cdot \deg \Omega^1_{\ol C}(\log \Delta_{\ol C}).
\end{align} 

\begin{theorem}[Lu--Zuo] \label{theorem:luzuo}
    Let $n \geq 3$ and $g \geq 5$ be integers. Let $C \subset \ca A_{g,[n]}$ be a smooth subvariety of dimension one with smooth projective model $C \subset \ol C$. Let $h \colon A \to C$ be the pull-back of the universal abelian scheme over $\ca A_{g,[n]}$ and assume that the local monodromy of $R^1h_\ast \Q$ around each point of $\Delta_{\ol C} = \ol C - C$ is unipotent. Let $E^{1,0}_{\ol C}$ be the vector bundle defined in \eqref{align:filtration-deligne-extension} above. Assume that $C$ is generically contained in the Torelli locus, and that $g \geq 12$ if $C$ is not contained in the hyperelliptic locus. Then \eqref{align:faltings:prelim} is a strict inequality.     
\end{theorem}
\begin{proof}
    See \cite[Theorem 1.4]{luzuo-shimuracurves} and its proof. Although the statement of \cite[Theorem 1.4]{luzuo-shimuracurves} only covers the $g \geq 12$ cases, its proof also deals with the $5 \leq g < 12$ cases under the additional hypothesis that the curve $C \subset \ca A_{g,[n]}$ generically contained in the hyperelliptic locus. 
\end{proof}

\begin{lemma} \label{lemma:arakelov-equality}
Let $g \geq 1, n \geq 3$ be integers. Let $C \subset \ca A_{g,[n]}$ be a smooth subvariety of dimension one and assume that the general element of $C$ is isogenous to the $g$-th power of an elliptic curve. Let $h \colon A \to C$ be the pull-back of the universal abelian scheme over $\ca A_{g,[n]}$ and assume that the local monodromy of $R^1h_\ast \Q$ around each point of $\Delta_{\ol C} = \ol C - C$ is unipotent. Let $E^{1,0}_{\ol C}$ be the vector bundle defined in \eqref{align:filtration-deligne-extension} above. Then \eqref{align:faltings:prelim} is an equality. 
\end{lemma}
\begin{proof}  
The vector bundle $E^{1,0}_{\ol C}$ is a direct summand of a vector bundle $E$ on $\ol C$ that has a natural Higgs bundle structure (see \cite{MR2074892,viehweg-zuo-characterization}) which by \cite{kollar-subadditivity} decomposes as a direct sum $E=F\oplus N$ of Higgs bundles such that $F^{1,0}_{\ol C} \coloneqq F \cap E^{1,0}_{\ol C}$ is ample and the Higgs field of $E$ vanishes on $N$. By Lemma \ref{lemma:special-lemma}, $C \subset \ca A_{g,[n]}$ is a one-dimensional special subvariety of $\ca A_{g,[n]}$. Therefore, by \cite{moller-shimura-teichmuller}, we have 
\begin{align} \label{align:maximal-higgs}
\deg F^{1,0}_{\ol C} = \frac{g_0}{2} \cdot \deg \Omega^1_{\ol C}(\Delta_{\ol C}), \quad \quad   \text{where}  \quad \quad g_0 \coloneqq \text{rank } F^{1,0}_{\ol C}. 
\end{align}
We claim that $\Delta_{\ol C}$ is non-empty. Indeed, the universal family of $\mathcal A_{g,[n]}$ restricted to $C$ is a non-isotrivial abelian scheme $h\colon A\to C$.
By assumption, this abelian scheme is (up to a finite surjective base change) isogenous to a self-fibre product of a family of elliptic curves.
If $C$ were proper, then one would conclude the properness of the moduli space of elliptic curves which is absurd. Hence, $\Delta_{\ol C}$ is non-empty, so that we can apply \cite[Theorem 0.2]{viehweg-zuo-characterization} to conclude that there exists an \'etale covering $\pi \colon C' \to C$ such that if $h' \colon A' = A \times_{C}C' \to C'$ denotes the pull-back of our family $h \colon A \to C$ along $\pi$, then $h' \colon A' \to C'$ is isogenous over $C'$ to a fibre product of the form $B' \times_{C'} \ca E' \times_{C'} \times \cdots \times_{C'} \ca E'$, where $B' \to C'$ is a constant family of $b$-dimensional abelian varieties over $Y'$ and $\ca E' \to C'$ is a non-isotrivial family of semi-stable elliptic curves over $Y'$. 
Here, $b = g-g_0 = g- \rm{rank}(F^{1,0}_{\ol C})$.  
Since the general fibre of $h \colon A \to C$ is isogenous to the $g$-th power of an elliptic curve and $h \colon A \to C$ is non-isotrivial, we see that $b = 0$ hence $g = g_0$. 
Thus, $N = 0$ and $E =F$, so that \eqref{align:maximal-higgs} implies  that $
\deg E^{1,0}_{\ol C} = (g/2) \cdot \deg \Omega^1_{\ol C}(\Delta_{\ol C})$.  
In other words, \eqref{align:faltings:prelim} is an equality, and we are done. 
\end{proof}

\subsection{Elliptic curves with no power isogenous to a Jacobian}  

\begin{proof}[Proof of Theorem \ref{theorem:verygeneral-powers-elliptic-isogenies}]
Let $E$ be an elliptic curve with transcendental $j$-invariant. Assume that, for some integer $g \geq 2$, we have an isogeny $\varphi \colon E^g \to JX$ where $X$ is a smooth projective connected curve. We must show that $g < 12$ and that $g < 5$ if $X$ is hyperelliptic. 

Since $E$ has transcendental $j$-invariant, the isogeny $\varphi \colon E^g \to JX$ spreads out to a one-dimensional family. More precisely, there is a one-dimensional variety $B$, a family of smooth projective connected curves $\ca X \to B$ whose fibres are hyperelliptic if $X$ is hyperelliptic, a non-isotrivial family of elliptic curves $\ca E \to B$ and an isogeny of abelian schemes $\varphi \colon \ca E^g \to J\ca X$ over $B$  that extends the isogeny $\varphi \colon E^g \to JC$. 
Up to replacing $B$ by an \'etale cover, we may assume that $J \ca X \to B$ is equipped with a level $n$ structure for some $n \geq 3$, so that it gives rise to a morphism $B \to \ca A_{g,[n]}$. Let $Z \subset \ca A_{g,[n]}$ denote the closure of the image of this map.  
By Lemma \ref{lemma:special-lemma}, $Z$ is a special subvariety of $\ca A_{g,[n]}$. 

By Lemma \ref{lemma:moonen-deligne}, there is an integer $m \geq 3$ with $n \mid m$ such that the following holds.
There is a one-dimensional smooth special subvariety $C \subset \ca A_{g,[m]}$ 
that dominates $Z$
such that the universal weight one $\Q$-local system on $C$ has unipotent monodromy at each point of $\Delta_{\ol C} = \ol C - C$, where $\ol C$ denotes the smooth projective model of $C$. 
Let $E^{1,0}_{\ol C}$ be the vector bundle defined in \eqref{align:filtration-deligne-extension} above. 
By Lemma \ref{lemma:arakelov-equality}, we have 
\begin{align} \label{align:arakelov-equality}
\deg(E^{1,0}_{\ol C}) = \frac{g}{2} \cdot \deg \Omega^1_{\ol C}(\log \Delta_{\ol C}).
\end{align}
Notice that $C \subset \ca A_{g,[m]}$ is generically contained in the Torelli locus. Thus, by Theorem \ref{theorem:luzuo}, the equality \eqref{align:arakelov-equality} implies that $g < 12$. Moreover, if the curve $X$ is hyperelliptic, then $C \subset \ca A_{g,[m]}$ is generically contained in the hyperelliptic Torelli locus. Therefore, we get $g < 5$ if $X$ is hyperelliptic, see Theorem \ref{theorem:luzuo}. This concludes the proof of the theorem.
\end{proof}

\printbibliography

@article {MR4332480,
    AUTHOR = {Chen, Ke and Lu, Xin and Zuo, Kang},
     TITLE = {On {CM} points away from the {T}orelli locus},
   JOURNAL = {J. Lond. Math. Soc. (2)},
  FJOURNAL = {Journal of the London Mathematical Society. Second Series},
    VOLUME = {104},
      YEAR = {2021},
    NUMBER = {3},
     PAGES = {1363--1383},
}

@article {yau-zhang-toroidal,
    AUTHOR = {Yau, Shing-Tung and Zhang, Yi},
     TITLE = {The geometry on smooth toroidal compactifications of {S}iegel
              varieties},
   JOURNAL = {Amer. J. Math.},
  FJOURNAL = {American Journal of Mathematics},
    VOLUME = {136},
      YEAR = {2014},
    NUMBER = {4},
     PAGES = {859--941},
}

@article {moonen-linearity-I,
    AUTHOR = {Moonen, Ben},
     TITLE = {Linearity properties of {S}himura varieties. {I}},
   JOURNAL = {J. Algebraic Geom.},
  FJOURNAL = {Journal of Algebraic Geometry},
    VOLUME = {7},
      YEAR = {1998},
    NUMBER = {3},
     PAGES = {539--567},
}

@article {moller-shimura-teichmuller,
    AUTHOR = {M{\"{o}}ller, Martin},
     TITLE = {Shimura and {T}eichm\"{u}ller curves},
   JOURNAL = {J. Mod. Dyn.},
  FJOURNAL = {Journal of Modern Dynamics},
    VOLUME = {5},
      YEAR = {2011},
    NUMBER = {1},
     PAGES = {1--32},
}

@book {serre-course-in-arithmetic,
    AUTHOR = {Serre, Jean-Pierre},
     TITLE = {A course in arithmetic},
    SERIES = {Graduate Texts in Mathematics, No. 7},
      NOTE = {Translated from the French},
 PUBLISHER = {Springer-Verlag, New York-Heidelberg},
      YEAR = {1973},
     PAGES = {viii+115},
   MRCLASS = {12-02 (10CXX 10DXX)},
  MRNUMBER = {344216},
}

@incollection {cattani-kaplan-degenerating,
    AUTHOR = {Cattani, Eduardo and Kaplan, Aroldo},
     TITLE = {Degenerating variations of {H}odge structure},
      NOTE = {Actes du Colloque de Th\'{e}orie de Hodge (Luminy, 1987)},
   JOURNAL = {Ast\'{e}risque},
  FJOURNAL = {Ast\'{e}risque},
    NUMBER = {179-180},
      YEAR = {1989},
     PAGES = {9, 67--96},
}

@article {MR2074892,
    AUTHOR = {Viehweg, Eckart and Zuo, Kang},
     TITLE = {Families over curves with a strictly maximal {H}iggs field},
   JOURNAL = {Asian J. Math.},
  FJOURNAL = {Asian Journal of Mathematics},
    VOLUME = {7},
      YEAR = {2003},
    NUMBER = {4},
     PAGES = {575--598},
}

@article {viehweg-zuo-characterization,
    AUTHOR = {Viehweg, Eckart and Zuo, Kang},
     TITLE = {A characterization of certain {S}himura curves in the moduli
              stack of abelian varieties},
   JOURNAL = {J. Differential Geom.},
  FJOURNAL = {Journal of Differential Geometry},
    VOLUME = {66},
      YEAR = {2004},
    NUMBER = {2},
     PAGES = {233--287},
}

@article {faltings-arakelov,
    AUTHOR = {Faltings, Gerd},
     TITLE = {Arakelov's theorem for abelian varieties},
   JOURNAL = {Invent. Math.},
  FJOURNAL = {Inventiones Mathematicae},
    VOLUME = {73},
      YEAR = {1983},
    NUMBER = {3},
     PAGES = {337--347},
}

@book {MR417174,
    AUTHOR = {Deligne, Pierre},
     TITLE = {\'{E}quations diff\'{e}rentielles \`a points singuliers r\'{e}guliers},
    SERIES = {Lecture Notes in Mathematics, Vol. 163},
 PUBLISHER = {Springer-Verlag, Berlin-New York},
      YEAR = {1970},
     PAGES = {iii+133},
}

@incollection {kollar-subadditivity,
    AUTHOR = {Koll\'{a}r, J\'{a}nos},
     TITLE = {Subadditivity of the {K}odaira dimension: fibers of general
              type},
 BOOKTITLE = {Algebraic geometry, {S}endai, 1985},
    SERIES = {Adv. Stud. Pure Math.},
    VOLUME = {10},
     PAGES = {361--398},
 PUBLISHER = {North-Holland, Amsterdam},
      YEAR = {1987},
}

@incollection {morrison-clemensschmid,
    AUTHOR = {Morrison, David},
     TITLE = {The {C}lemens-{S}chmid exact sequence and applications},
 BOOKTITLE = {Topics in transcendental algebraic geometry ({P}rinceton,
              {N}.{J}., 1981/1982)},
    SERIES = {Ann. of Math. Stud.},
    VOLUME = {106},
     PAGES = {101--119},
 PUBLISHER = {Princeton Univ. Press, Princeton, NJ},
      YEAR = {1984},
}

@article {chiang-lipman-simultaneous,
    AUTHOR = {Chiang-Hsieh, Hung-Jen and Lipman, Joseph},
     TITLE = {A numerical criterion for simultaneous normalization},
  JOURNAL = {Duke Mathematical Journal},
    VOLUME = {133},
      YEAR = {2006},
    NUMBER = {2},
     PAGES = {347--390},
}

@book {faltings-chai,
    AUTHOR = {Faltings, Gerd and Chai, Ching-Li},
     TITLE = {Degeneration of {A}belian {V}arieties},
    SERIES = {Ergebnisse der Mathematik und ihrer Grenzgebiete (3)},
    VOLUME = {22},
 PUBLISHER = {Springer-Verlag, Berlin},
      YEAR = {1990},
     PAGES = {xii+316},
}

@article {marcucci-naranjo-pirola,
    AUTHOR = {Marcucci, Valeria and Naranjo, Juan Carlos and Pirola, Gian
              Pietro},
     TITLE = {Isogenies of {J}acobians},
  JOURNAL = {Algebraic Geometry},
    VOLUME = {3},
      YEAR = {2016},
    NUMBER = {4},
     PAGES = {424--440},
}

@article {lange-sernesi,
    AUTHOR = {Lange, Herbert and Sernesi, Edoardo},
     TITLE = {Curves of genus {$g$} on an abelian variety of dimension
              {$g$}},
  JOURNAL = {Koninklijke Nederlandse Akademie van Wetenschappen.
              Indagationes Mathematicae. New Series},
    VOLUME = {13},
      YEAR = {2002},
    NUMBER = {4},
     PAGES = {523--535},
}

@article {masser-zannier,
    AUTHOR = {Masser, David and Zannier, Umberto},
     TITLE = {Abelian varieties isogenous to no {J}acobian},
  JOURNAL = {Annals of Mathematics. Second Series},
    VOLUME = {191},
      YEAR = {2020},
    NUMBER = {2},
     PAGES = {635--674},
}

@article {tsimerman,
    AUTHOR = {Tsimerman, Jacob},
     TITLE = {The existence of an abelian variety over {$\overline{\mathbb Q}$}
              isogenous to no {J}acobian},
  JOURNAL = {Annals of Mathematics. Second Series},
    VOLUME = {176},
      YEAR = {2012},
    NUMBER = {1},
     PAGES = {637--650},
}

@book {hatcher-AT,
    AUTHOR = {Hatcher, Allen},
     TITLE = {Algebraic {T}opology},
 PUBLISHER = {Cambridge University Press, Cambridge},
      YEAR = {2002},
     PAGES = {xii+544},
}

@article {kollar-plurigenera,
    AUTHOR = {Koll\'{a}r, J\'{a}nos},
     TITLE = {Shafarevich maps and plurigenera of algebraic varieties},
  JOURNAL = {Inventiones Mathematicae},
    VOLUME = {113},
      YEAR = {1993},
    NUMBER = {1},
     PAGES = {177--215},
}

@article {chai-oort,
    AUTHOR = {Chai, Ching-Li and Oort, Frans},
     TITLE = {Abelian varieties isogenous to a Jacobian},
  JOURNAL = {Annals of Mathematics. Second Series},
    NUMBER = {176},
      YEAR = {2012},
     PAGES = {589–-635},
}

@article {dejong-alterations,
    AUTHOR = {de Jong, Aise Johan},
     TITLE = {Smoothness, semi-stability and alterations},
  JOURNAL = {Institut des Hautes \'{E}tudes Scientifiques. Publications
              Math\'{e}matiques},
    NUMBER = {83},
      YEAR = {1996},
     PAGES = {51--93},
}

@book {voisin-II,
    AUTHOR = {Voisin, Claire},
     TITLE = {Hodge {T}heory and {C}omplex {A}lgebraic {G}eometry. {II}},
    SERIES = {Cambridge Studies in Advanced Mathematics},
    VOLUME = {77},
   EDITION = {English},
 PUBLISHER = {Cambridge University Press, Cambridge},
      YEAR = {2007},
     PAGES = {x+351},
}

@book {petersteenbrink,
    AUTHOR = {Peters, Chris and Steenbrink, Joseph},
     TITLE = {Mixed {H}odge {S}tructures},
    SERIES = {Ergebnisse der Mathematik und ihrer Grenzgebiete. 3. Folge. A
              Series of Modern Surveys in Mathematics},
    VOLUME = {52},
 PUBLISHER = {Springer-Verlag, Berlin},
      YEAR = {2008},
     PAGES = {xiv+470},
}

@Inbook{Milne1986,
author="Milne, James",
editor="Cornell, Gary
and Silverman, Joseph",
title="Abelian Varieties",
bookTitle="Arithmetic Geometry",
year="1986",
publisher="Springer New York",
pages="103--150",
}

@article {kneser-quadratisch,
    AUTHOR = {Kneser, Martin},
     TITLE = {Klassenzahlen definiter quadratischer {F}ormen},
  JOURNAL = {Archiv der Mathematik},
    VOLUME = {8},
      YEAR = {1957},
     PAGES = {241--250},
}

@article {schoen-hyperelliptic,
    AUTHOR = {Schoen, Chad},
     TITLE = {Bounds for rational points on twists of constant hyperelliptic
              curves},
  JOURNAL = {J.\ Reine Angew.\ Math.},
    VOLUME = {411},
      YEAR = {1990},
     PAGES = {196--204},
}

@incollection {collino-fanofundamentalgroup,
    AUTHOR = {Collino, Alberto},
     TITLE = {The fundamental group of the {F}ano surface. {I}},
 BOOKTITLE = {Algebraic threefolds ({V}arenna, 1981)},
    SERIES = {Lecture Notes in Math.},
    VOLUME = {947},
     PAGES = {209--218, 219--220},
 PUBLISHER = {Springer, Berlin-New York},
      YEAR = {1982},
}

@article {voisin-universalCHgroup,
    AUTHOR = {Voisin, Claire},
     TITLE = {On the universal {$\rm CH_0$} group of cubic hypersurfaces},
  JOURNAL = {Journal of the European Mathematical Society (JEMS)},
    VOLUME = {19},
      YEAR = {2017},
    NUMBER = {6},
     PAGES = {1619--1653},
}

@misc{voisin2022cycle,
      title={Cycle classes on abelian varieties and the geometry of the Abel-Jacobi map}, 
      author={Claire Voisin},
      year={2022},
      eprint={2212.03046},
      archivePrefix={arXiv},
      primaryClass={math.AG},
note={To appear in Pure and Applied Mathematics Quarterly}
}

@article {beckmann-degaayfortman,
    AUTHOR = {Beckmann, Thorsten and de Gaay Fortman, Olivier},
     TITLE = {Integral {F}ourier transforms and the integral {H}odge
              conjecture for one-cycles on abelian varieties},
  JOURNAL = {Compositio Mathematica},
    VOLUME = {159},
      YEAR = {2023},
    NUMBER = {6},
     PAGES = {1188--1213},
}

@article {naranjopirola1994,
    AUTHOR = {Naranjo, Juan Carlos and Pirola, Gian Pietro},
     TITLE = {On the genus of curves in the generic {P}rym variety},
  JOURNAL = {Koninklijke Nederlandse Akademie van Wetenschappen.
              Indagationes Mathematicae. New Series},
    VOLUME = {5},
      YEAR = {1994},
    NUMBER = {1},
     PAGES = {101--105},
}

@article {naranjopirola2018,
    AUTHOR = {Naranjo, Juan Carlos and Pirola, Gian Pietro},
     TITLE = {Hyperelliptic {J}acobians and isogenies},
  JOURNAL = {Advances in Mathematics},
    VOLUME = {335},
      YEAR = {2018},
     PAGES = {896--909},
}

@article {lazarsfeld2023measures,
    AUTHOR = {Robert Lazarsfeld and Olivier Martin},
     TITLE = {Measures of association between algebraic varieties, II: {S}elf-correspondences},
  JOURNAL = {Épijournal de Géométrie Algébrique},
    VOLUME = {10},
      YEAR = {2023},
}

@book {chai-siegelcompactifications,
    AUTHOR = {Chai, Ching-Li},
     TITLE = {Compactification of {S}iegel moduli schemes},
    SERIES = {London Mathematical Society Lecture Note Series},
    VOLUME = {107},
 PUBLISHER = {Cambridge University Press, Cambridge},
      YEAR = {1985},
     PAGES = {xvi+326},
}

@article {alexeev-compactifiedjacobians,
    AUTHOR = {Alexeev, Valery},
     TITLE = {Compactified {J}acobians and {T}orelli map},
  JOURNAL = {Kyoto University. Research Institute for Mathematical
              Sciences. Publications},
    VOLUME = {40},
      YEAR = {2004},
    NUMBER = {4},
     PAGES = {1241--1265},
}

@book {humphreys-reflectiongroups,
    AUTHOR = {Humphreys, James},
     TITLE = {Reflection groups and {C}oxeter groups},
    SERIES = {Cambridge Studies in Advanced Mathematics},
    VOLUME = {29},
 PUBLISHER = {Cambridge University Press, Cambridge},
      YEAR = {1990},
     PAGES = {xii+204},
}

@incollection {moonenoort-specialsubvarieties,
    AUTHOR = {Moonen, Ben and Oort, Frans},
     TITLE = {The {T}orelli locus and special subvarieties},
 BOOKTITLE = {Handbook of {M}oduli. {V}ol. {II}},
    SERIES = {Adv. Lect. Math. (ALM)},
    VOLUME = {25},
     PAGES = {549--594},
 PUBLISHER = {Int. Press, Somerville, MA},
      YEAR = {2013},
}

@book {milnorhusemoller,
    AUTHOR = {Milnor, John and Husemoller, Dale},
     TITLE = {Symmetric {B}ilinear {F}orms},
    SERIES = {Ergebnisse der Mathematik und ihrer Grenzgebiete, Band 73},
 PUBLISHER = {Springer-Verlag, New York-Heidelberg},
      YEAR = {1973},
     PAGES = {viii+147},
   MRCLASS = {15A63 (10C05 57D65)},
  MRNUMBER = {0506372},
MRREVIEWER = {Louis H. Kauffman},
}

@article {marcucci-genusofcurves,
    AUTHOR = {Marcucci, Valeria Ornella},
     TITLE = {On the genus of curves in a {J}acobian variety},
  JOURNAL = {Annali della Scuola Normale Superiore di Pisa. Classe di
              Scienze. Serie V},
    VOLUME = {12},
      YEAR = {2013},
    NUMBER = {3},
     PAGES = {735--754},
}

@article {bardellipirola-curvesofgenusg,
    AUTHOR = {Bardelli, Fabio and Pirola, Gian Pietro},
     TITLE = {Curves of genus {$g$} lying on a {$g$}-dimensional {J}acobian
              variety},
  JOURNAL = {Inventiones Mathematicae},
    VOLUME = {95},
      YEAR = {1989},
    NUMBER = {2},
     PAGES = {263--276},
}

@incollection {carlson-extensions,
    AUTHOR = {Carlson, James},
     TITLE = {Extensions of mixed {H}odge structures},
 BOOKTITLE = {Journ\'{e}es de {G}\'{e}ometrie {A}lg\'{e}brique d'{A}ngers, {J}uillet
              1979/{A}lgebraic {G}eometry, {A}ngers, 1979},
     PAGES = {107--127},
 PUBLISHER = {Sijthoff \& Noordhoff, Alphen aan den Rijn---Germantown, Md.},
      YEAR = {1980},
}

@article {clemensgriffiths-cubicthreefold,
    AUTHOR = {Clemens, Herbert and Griffiths, Phillip},
     TITLE = {The intermediate {J}acobian of the cubic threefold},
  JOURNAL = {Annals of Mathematics. Second Series},
    VOLUME = {95},
      YEAR = {1972},
     PAGES = {281--356},
}

@article {luzuo-shimuracurves,
    AUTHOR = {Lu, Xin and Zuo, Kang},
     TITLE = {The {O}ort conjecture on {S}himura curves in the {T}orelli
              locus of curves},
  JOURNAL = {Journal de Math\'{e}matiques Pures et Appliqu\'{e}es. Neuvi\`eme S\'{e}rie},
    VOLUME = {123},
      YEAR = {2019},
     PAGES = {41--77},
}

@article {debarre-produits,
    AUTHOR = {Debarre, Olivier},
     TITLE = {Polarisations sur les vari\'{e}t\'{e}s ab\'{e}liennes produits},
   JOURNAL = {C. R. Acad. Sci. Paris S\'{e}r. I Math.},
  FJOURNAL = {Comptes Rendus de l'Acad\'{e}mie des Sciences. S\'{e}rie I.
              Math\'{e}matique},
    VOLUME = {323},
      YEAR = {1996},
    NUMBER = {6},
     PAGES = {631--635},
}

@book {birkenhakelange-complexAVs,
    AUTHOR = {Birkenhake, Christina and Lange, Herbert},
     TITLE = {Complex {A}belian {V}arieties},
    SERIES = {Grundlehren der mathematischen Wissenschaften [Fundamental
              Principles of Mathematical Sciences]},
    VOLUME = {302},
   EDITION = {Second},
 PUBLISHER = {Springer-Verlag, Berlin},
      YEAR = {2004},
     PAGES = {xii+635},
}

\end{document}